\documentclass[11pt]{amsart}
\usepackage{amsxtra}
\usepackage{amssymb}
\usepackage{mathtools}
\usepackage{tikz-cd}
\usetikzlibrary{chains}
\usetikzlibrary{cd,shapes}
\usetikzlibrary{arrows}
\usepackage[mathscr]{eucal}
\input xy \xyoption{matrix} \xyoption{arrow}\xyoption{frame}

\theoremstyle{plain}
\usepackage{color}
\newcommand \red {\textcolor{black}}

\tikzset{node distance=2em, ch/.style={circle,draw,on chain,inner sep=2pt},chj/.style={ch,join},every path/.style={shorten >=4pt,shorten <=4pt},line width=1pt,baseline=-1ex}

\let\dlabel=\alabel
\let\ulabel=\mlabel

\newcommand{\dnode}[2][chj]{%
\node[#1,label={below:\dlabel{#2}}] {};
}

\newcommand{\dnodea}[3][chj]{%
\dnode[#1,label={above:\ulabel{#2}}]{#3}
}

\newcommand{\dnodeanj}[2]{%
\dnodea[ch]{#1}{#2}
}

\newcommand{\QLeftarrow}{%
\begingroup
\tikz
\draw[shorten >=0pt,shorten <=0pt] (0,3pt) -- ++(-1em,0) (0,1pt) -- ++(-1em-1pt,0) (0,-1pt) -- ++(-1em-1pt,0) (0,-3pt) -- ++(-1em,0) (-1em+1pt,5pt) to[out=-105,in=45] (-1em-2pt,0) to[out=-45,in=105] (-1em+1pt,-5pt);
\endgroup
}
\tikzstyle{mutable}=[align=center, draw,inner sep=0pt,shape=ellipse,minimum height=0.4cm, minimum width=0.4cm]
\tikzstyle{frozen}=[inner sep=0.5mm,rectangle,draw]

\newcommand{\cleqn}{\setcounter{equation}{0}}
\newcommand{\clth}{\setcounter{theorem}{0}}
\newcommand {\sectionnew}[1]{\section{#1}\cleqn\clth}
\newcommand{\nn}{\hfill\nonumber}
\newtheorem{theorem}{Theorem}[section]
\newtheorem{lemma}[theorem]{Lemma}
\newtheorem{definition-theorem}[theorem]{Definition-Theorem}
\newtheorem{proposition}[theorem]{Proposition}
\newtheorem{corollary}[theorem]{Corollary}
\newtheorem{definition}[theorem]{Definition}
\newtheorem{example}[theorem]{Example}
\newtheorem{remark}[theorem]{Remark}
\newtheorem{conjecture}[theorem]{Conjecture}
\newtheorem{notation}[theorem]{Notation}

\newtheorem*{maintheorem*}{Main Theorem}
\newtheorem*{theorem*}{Theorem}
\newtheorem*{theoremA*}{Theorem A}
\newtheorem*{theoremB*}{Theorem B}
\theoremstyle{definition}

\newtheorem*{definition*}{Definition}
\newcommand \bth[1] { \begin{theorem}\label{t#1} }
\newcommand \ble[1] { \begin{lemma}\label{l#1} }

\newcommand \bpr[1] { \begin{proposition}\label{p#1} }
\newcommand \bco[1] { \begin{corollary}\label{c#1} }
\newcommand \bde[1] { \begin{definition}\label{d#1}\rm }
\newcommand \bex[1] { \begin{example}\label{e#1}\rm }
\newcommand \bre[1] { \begin{remark}\label{r#1}\rm }
\newcommand \bcj[1] { \begin{conjecture}\label{j#1}\rm }

\newcommand \bnota[1] { \begin{notation}\label{n#1}\rm }
\renewcommand {\eth} { \end{theorem} }
\newcommand {\ele} { \end{lemma} }

\newcommand {\epr} { \end{proposition} }
\newcommand {\eco} { \end{corollary} }
\newcommand {\ede} { \end{definition} }
\newcommand {\eex} { \end{example} }
\newcommand {\ere} { \end{remark} }
\newcommand {\ecj} { \end{conjecture} }

\newcommand {\enota} { \end{notation} }
\newcommand \thref[1]{Theorem \ref{t#1}}
\newcommand \leref[1]{Lemma \ref{l#1}}
\newcommand \prref[1]{Proposition \ref{p#1}}
\newcommand \coref[1]{Corollary \ref{c#1}}

\newcommand \deref[1]{Definition \ref{d#1}}
\newcommand \exref[1]{Example \ref{e#1}}
\newcommand \reref[1]{Remark \ref{r#1}}

\def \Cset {{\mathbb C}}
\def \DD {{\mathbb D}}
\def \KK {{\mathbb K}}

\def \Zset {{\mathbb Z}}
\def \Nset {{\mathbb N}}
\def \Qset {{\mathbb Q}}

\def \iA {{\mathcal{A}}}             
\def \iAA {{\mathcal{A}}^{1/2}}
\def \iFF {{\mathcal{F}}^{1/2}}
\def \AA  {{\boldsymbol{\mathsf A}}}    
\def \UU {{\boldsymbol{\mathsf U}}}
\def \FF {{\mathcal{F}}}
\def \EE {{\mathcal{E}}}

\def \TT {{\mathcal{T}}}

\def \OO {{\mathcal{O}}}
\def \Ointg  {{ \mathcal{O}_{\mathrm{int}}({\mathfrak{g}})}}

\def \JJ {{\mathcal{J}}}

\def \Sbb {{\mathcal{S}}}

\def \TT {{\mathcal{T}}} 

\def \BB {\bf{B}}

\def \rb {{\bf{r}}}

\def \tb {{\bf{t}}}
\def \ex {{\bf{ex}}}
\def \Scr {{\mathscr{S}}}

\def \inv {{\bf{inv}}}
\def \de {\delta}
\def \al {\alpha}
\def \be {\beta}
\def \vpi {\varpi}
\def \la {\lambda}

\def \Om {\Omega}
\def \ga {\gamma}
\def \de {\delta}
\def \Ga {\Gamma}

\def \sig {\sigma}

\def \vp {\varphi}
\def \ep {\epsilon}
\def \De {\Delta}
\def \rb { {\bf{r}} }

\def \mt  {\mapsto}
\def \lra {\longrightarrow}

\def \hra {\hookrightarrow}

\def \sy  {\ast}                         

\def \rcor {\rangle}
\def \lcor {\langle}

\def \ol {\overline}
\def \wt {\widetilde}

\def \id { {\mathrm{id}} }

\DeclareMathOperator \rk { {\mathrm{rk}} }

\def \g  {\mathfrak{g}}   

\def \sl {\mathfrak{sl}} 
\def \h  {\mathfrak{h}}

\def \n  {\mathfrak{n}}

\def \b  {\mathfrak{b}}

\def \sl {\mathfrak{sl}}

\DeclareMathOperator \wtt { {\mathrm{wt}}}
\DeclareMathOperator \sign {{\mathrm{sign}}}
\DeclareMathOperator \Aut { {\mathrm{Aut}} }

\DeclareMathOperator \charr { {\mathrm{char}} }

\DeclareMathOperator \diag { {\mathrm{diag}} }

\DeclareMathOperator \Hom { {\mathrm{Hom}} }

\DeclareMathOperator \lt  { {\mathrm{lt}} }
\DeclareMathOperator \lc  { {\mathrm{lc}} }

\DeclareMathOperator \Fract { {\mathrm{Fract}} }

\renewcommand \Im { {\mathrm{Im}} }

\DeclareMathOperator*{\moplus}{\text{\raisebox{0.25ex}{\scalebox{0.8}{$\bigoplus$}}}}
\newcommand\kx{\KK^*}

\newcommand\HH{{\mathcal{H}}}
\newcommand\xh{X(\HH)}
\DeclareMathOperator \Spec {Spec}

\newcommand \Znn {\Zset_{\ge 0}}

\DeclareMathOperator{\up}{up} 
\DeclareMathOperator{\low}{low}

\def\lab{{\boldsymbol \lambda}}

\def\nub{{\boldsymbol \nu}}

\begin{document}
\title[Integral quantum cluster structures]
{Integral quantum cluster structures}
\author[K. R. Goodearl]{K. R. Goodearl}
\address{
Department of Mathematics \\
University of California\\
Santa Barbara, CA 93106 \\
U.S.A.
}
\email{goodearl@math.ucsb.edu}
\thanks{The research of K.R.G was partially supported by NSF grant DMS-1601184, and that of 
M.T.Y. by NSF grant DMS-1601862 and Bulgarian Science Fund grant H02/15.}
\author[M. T. Yakimov]{M. T. Yakimov}
\address{
Department of Mathematics \\
Louisiana State University \\
Baton Rouge, LA 70803 \\
U.S.A.
}
\email{yakimov@math.lsu.edu}
\begin{abstract}
We prove a general theorem for constructing integral quantum cluster algebras over $\Zset[q^{\pm 1/2}]$, namely
that under mild conditions the integral forms of quantum nilpotent algebras always possess integral quantum cluster algebra structures. 
These algebras are then shown to be isomorphic to the corresponding upper 
quantum cluster algebras, again defined over $\Zset[q^{\pm 1/2}]$. Previously, this was only known for 
acyclic quantum cluster algebras. The theorem is applied to prove that for
every symmetrizable Kac--Moody algebra $\g$ and Weyl group element $w$, the dual canonical form 
$A_q(\n_+(w))_{\Zset[q^{\pm 1}]}$ of the corresponding quantum unipotent cell has the property that 
$A_q(\n_+(w))_{\Zset[q^{\pm 1}]} \otimes_{\Zset[q^{ \pm 1}]} \Zset [ q^{\pm 1/2}]$ 
is isomorphic to a quantum cluster algebra over $\Zset[q^{\pm 1/2}]$ and to the corresponding 
upper quantum cluster algebra over $\Zset[q^{\pm 1/2}]$.
\end{abstract}
\date{May 6, 2019}
\maketitle

\sectionnew{Introduction} \label{intro}
\subsection{Problems for integral quantum custer algebras}
Cluster algebras were introduced by Fomin and Zelevinsky in \cite{FZ1} and have been applied to a number of diverse areas such 
as representation theory, combinatorics, Poisson and algebraic geometry, mathematical physics and others. Their quantum 
counterparts, introduced by Berenstein and Zelevinsky \cite{BZ}, are similarly the topic of 
intensive research from various standpoints. In the uniparameter quantum case it is desirable to work over the minimal ring of definition, namely over
\begin{equation}
\label{iAA}
\iAA := \Zset[q^{\pm 1/2}],
\end{equation}
where $q$ is the quantum parameter.
We will refer to such structures as  {\em{integral quantum cluster algebras}}. 
Two fundamental problems that are being investigated are:
\begin{enumerate}
\item Given an algebra $R$ over the rational function field $\iFF:= \Qset(q^{1/2})$ and an integral form 
$R_{\iAA}$ of $R$ over $\iAA$ (i.e., $R \cong R_{\iAA} \otimes_{\iAA} \iFF$), when is $R_{\iAA}$
isomorphic to an integral quantum cluster algebra? 
\item When is the quantum cluster algebra $\AA$ in Problem (1) equal to the corresponding upper 
quantum cluster algebra $\UU$ defined over $\iAA$?
\end{enumerate}
The best known result on Problem (1) is a theorem of Kang, Kashiwara, Kim and Oh \cite{KKKO}
that the dual canonical forms (over $\iAA$) of the quantum unipotent cells for all symmetric 
Kac--Moody algebras possess integral quantum cluster algebra structures. Berenstein and Zelevinsky 
\cite{BZ} proved the equality $\AA = \UU$ in the acyclic case. Such an equality was proved by Muller \cite{Mu}
for (quantum) cluster algebras that are source-sink decomposable in the case when all 
frozen variables are inverted.
We are not aware of any affirmative solutions of Problem (2) for non-acyclic quantum cluster algebras
when frozen variables are not inverted.
A recent result of Gei\ss, Leclerc and Schr\"oer \cite{GLS18} establishes an equality of the form
\[
\AA \otimes_{\iAA} \Qset[q^{\pm 1/2}] = \UU \otimes_{\iAA} \Qset[q^{\pm 1/2}]
\] 
under the assumptions that $\AA$ is connected $\Zset_{\geq 0}$-graded with homogeneous cluster variables and that 
such an equality holds on the classical level.

\subsection{Main results} 
In this paper we provide affirmative answers to both Problems (1) and (2) in wide generality. As an application, 
affirmative answers to Problems (1) and (2) are obtained for the dual canonical forms of the quantum unipotent cells for all 
symmetrizable Kac--Moody algebras. 

For an iterated skew polynomial extension
\[
R := \iFF[x_1][x_2; \theta_2, \delta_2] \cdots [x_N; \theta_N, \delta_N]
\]
and $1 \leq j \leq k \leq N$,  
denote by $R_{[j,k]}$ the $\iFF$-subalgebra generated by $x_j, \ldots, x_k$ and set $R_k := R_{[1,k]}$. 

\begin{definition*}
An iterated skew polynomial extension $R$ is called a \emph{quantum nilpotent algebra} or a \emph{CGL extension} 
if it is equipped with a rational action of an $\iFF$-torus $\HH$ by $\iFF$-algebra automorphisms such that:
\begin{enumerate}
\item[(i)] The elements $x_1, \ldots, x_N$ are $\HH$-eigenvectors.
\item[(ii)] For every $k \in [2,N]$, $\de_k$ is a locally nilpotent 
$\theta_k$-derivation of the algebra $R_{k-1}$. 
\item[(iii)] For every $k \in [1,N]$, there exists $h_k \in \HH$ such that 
$\theta_k = (h_k \cdot)|_{R_{k-1}}$ and the $h_k$-eigenvalue of $x_k$, to be denoted by $\la_k$, is not a root of unity.
\end{enumerate}
\end{definition*}

A CGL extension is called {\em{symmetric}} if it has the same properties when its generators are adjoined in the opposite order. 
We will assume throughout Sections \ref{intro}, \ref{qSchubert}, \ref{Hprime-sec}, \ref{intAqn+w} that the $\theta_k$-eigenvalues of $x_j$ belong to $q^{\Zset/2}$ for $j \leq k$, where we abbreviate $\Zset/2 := \Zset\frac12$. 
Recall that a nonzero element $p \in R$ is called \emph{prime} 
if $Rp = pR$ and the ring $R/Rp$ is a domain.

\begin{theorem*} 
\cite[Theorem 4.3]{GY1} For each CGL extension $R$ and $k \in [1,N]$, 
the algebra $R_k$ has a unique {\em{(}}up to rescaling{\em{)}} homogeneous prime 
element $y_k$ which does not belong to $R_{k-1}$. It either equals $x_k$ or has the 
property that
\[
y_k - y_{p(k)} x_k \in R_{k-1}
\]
for some $p(k) \in [1,k-1]$. 
\end{theorem*}

In the following we will work with this choice of sets of homogeneous prime elements (and not with arbitrary $\iFF$-rescalings of them).
For a symmetric CGL extension the theorem can be applied to the interval subalgebras $R_{[p(k), k]}$ to obtain that each of them has a 
unique (up to rescaling) homogeneous prime element $y_{[p(k),k]}$ which does not belong to the smaller interval subalgebras. An $\iFF$-rescaling of the generators of a CGL 
extension $R$ leads to another CGL extension presentation of $R$. The generators $x_k$ can be always rescaled so that
\begin{equation}
\label{normaliz}
y_{[p(k), k]} = q^{m} x_{p(k)} x_k - q^{m'} \prod_i p_i^{n_i}
\end{equation}
for some $m, m' \in \Zset/2$ and $n_i \in \Zset_{\geq 0}$ where the product is over all homogeneous prime elements of $R_{[p(k), k]}$ 
from the theorem that are different from $y_{[p(k),k]}$
(see \S \ref{rescalegens}-\ref{normalizations} and \cite{GY1}). In the following we will assume that this 
normalization is made. Denote
\[
R_{\iAA} := \iAA \lcor x_1, \ldots, x_N \rcor \subseteq R.
\]

\begin{theoremA*} 
Let $R$ be a symmetric CGL extension for which $R_{\iAA}$ is an $\iAA$-form of $R$, that is, 
$R_{\iAA} \otimes_{\iAA} \iFF \cong R$. If the sequence of homogeneous prime elements $y_1, \ldots, y_N$ lies in $R_{\iAA}$, 
then there exists a quantum cluster algebra $\AA$ over $\iAA$ such that
\[
R \cong \AA = \UU
\]
where $\UU$ is the corresponding upper quantum cluster algebra over $\iAA$. For all $k \in [1,N]$ 
and $n \in \Zset_{>0}$ for which $p^n(k)$ is well defined {\em{(}}as in the previous theorem{\em{)}}, $q^{m} x_k$ and $q^{m'} y_{[p^n(k), k]}$ 
are cluster variables of $\AA$ for some $m, m' \in \Zset/2$.
\end{theoremA*}

We prove a more general result in \thref{qcl.intform} which deals with integral forms of multiparameter and arbitrary 
characteristic CGL extensions and quantum cluster algebras. 
In \S \ref{exam} we illustrate the theorem with various examples which are not connected $\Zset_{\geq 0}$-graded, 
including all quantized Weyl algebras, and with quantum cluster algebras over ${\mathbb{F}}_p[q^{\pm 1/2}]$. 

For each symmetrizable Kac--Moody algebra $\g$ and a Weyl group element $w$, De Concini--Kac--Procesi \cite{DKP} and Lusztig \cite{Lusztig2}
defined a quantum Schubert cell algebra $U^-[w]$ which is a subalgebra 
of the quantized universal enveloping algebra $U_q(\g)$ defined over $\Qset(q)$. 
The quantum unipotent cells of Gei\ss--Leclerc--Schr\"oer \cite{GLS} are $\Qset(q)$-algebras
$A_q(\n_+(w))$ which are antiisomorphic to  $U^-[w]$. Denote 
\begin{equation}
\label{iA}
\iA := \Zset[q^{\pm 1}].
\end{equation}
The dual canonical forms $A_q(\n_+(w))_{\iA}$ are $\iA$-forms of  $A_q(\n_+(w))$ which are obtained 
by transporting the Kashiwara--Lusztig dual canonical forms $U^-[w]\spcheck_{\iA}$ of $U^-[w]$. 

\begin{theoremB*}
Let $\g$ be an arbitrary symmetrizable Kac--Moody algebra and $w$ a Weyl group element. 
For the dual canonical form $A_q(\n_+(w))_{\iA}$ of the corresponding quantum unipotent cell,
there exists a quantum cluster algebra $\AA$ over $\iAA$ such that
\[
A_q(\n_+(w))_{\iA} \otimes_{\iA} \iAA \cong \AA = \UU
\]
where $\UU$ is the associated upper quantum cluster algebra defined over $\iAA$. 
\end{theoremB*}

Further details about the structure of the quantum cluster algebra $\AA$ are given in \thref{main2}.

The following special cases of parts of the theorem were previously proved: Qin \cite{Q} proved that $A_q(\n_+(w))_{\iA} \otimes_{\iA}\, \iAA \cong \AA$
for symmetric Kac--Moody algebras $\g$ and adaptable Weyl group elements $w$. Kang, Kashiwara, Kim and Oh \cite{KKKO} proved this 
isomorphism for symmetric Kac--Moody algebras $\g$ and all Weyl group elements $w$. Gei\ss, Leclerc and Schr\"oer \cite{GLS18} proved that
\[
\AA \otimes_{\iAA} \Qset[q^{\pm 1/2}] = \UU \otimes_{\iAA}  \Qset[q^{\pm 1/2}] 
\]
for symmetric Kac--Moody algebras $\g$ and all Weyl group elements $w$; however, the fact that $\AA= \UU$ is new even for simple cases 
like $\g = \sl_n$. For nonsymmetric Kac--Moody algebras $\g$ the results in the theorem are all new, including the  existence of a 
non-integral quantum cluster structure on $A_q(\n_+(w))_{\iA} \otimes_{\iA} \Qset(q^{1/2})$.

The previous approaches to integral quantum cluster structures \cite{CW,HL,KKKO,KKOP,N,Q} 
obtained monoidal categorifications of quantum cluster algebras. At the same time they also relied on extensive knowledge 
of categorifications which are available for concrete families of algebras.
The power of Theorem A for the construction of integral quantum cluster structures
lies in its flexibility to adjust to different situations and in the mild assumptions in it: one needs to 
verify the normalization condition \eqref{normaliz}, that $R_{\iAA}$ is an $\iAA$-form of $R$, and that 
the sequence of homogeneous prime elements $y_1, \ldots, y_N$ belongs to $R_{\iAA}$. 

\subsection{Notation and conventions}  
Throughout, $\KK$ denotes an infinite field of arbitrary characteristic. For integers $j \leq k$, set $[j,k]:=\{j, \ldots, k\}$. As above, $\Zset/2 := \Zset\frac12$.

An $N\times N$ matrix $\tb = (t_{kj})$ over a commutative ring $\DD$ is \emph{multiplicatively skew-symmetric} 
if $t_{jk} t_{kj} = t_{kk} = 1$ for all $j,k \in [1,N]$. Such a matrix gives rise to a skew-symmetric bicharacter $\Om_\tb : \Zset^N \times \Zset^N \rightarrow \DD^*$ 
for which
\begin{equation}  \label{Omt}
\Om_\tb(e_j,e_k) = t_{jk} \,, \quad \forall\, j,k \in [1,N],
\end{equation}
where $e_1,\dots,e_N$ are the standard basis vectors for $\Zset^N$. (We denote the group of units of $\DD$ by $\DD^*$.)
When we have need for formulas involving $\Zset^N$, we view its elements as column vectors. The transpose of an $N$-tuple 
$\mathbf{m} = (m_1,\dots,m_N)$ is denoted $\mathbf{m}^T$.

Given an algebra $A$ over a commutative ring $\DD$ and elements $a_1,\dots,a_k \in A$, we write $\DD \langle a_1,\dots,a_k \rangle$ to denote the unital $\DD$-subalgebra of $A$ generated by $\{a_1,\dots,a_k\}$.
\medskip
\\
\noindent
{\bf Acknowledgements.} We would like to thank Bernhard Keller for valuable suggestions. 
\red{We would also like to thank the anonymous referee whose suggestions were very helpful to us in improving the paper.}
\sectionnew{Quantum cluster algebras}  \label{qcluster-sec}
\label{q-cl}

We outline notation and conventions for quantum cluster algebras. To connect with the results of \cite{GYbig}, we describe a multiparameter setting which extends the uniparameter case originally developed by Berenstein and Zelevinsky \cite{BZ}. To allow for integral forms, we work over a commutative domain rather than over a field.

Fix a commutative domain $\DD$ contained in $\KK$ and a positive integer $N$.

Let $\FF$ be a division algebra over $\DD$. A {\em{toric frame}} (of rank $N$) for $\FF$ (over $\DD$) is a mapping 
$$
M : \Zset^N \longrightarrow \FF
$$
such that 
\begin{equation}
\label{MOm}
M(f) M(g) = \Om_\rb(f,g) M(f+g), \quad
\forall\,  f,g \in \Zset^N,
\end{equation}
where
\begin{itemize}
\item $\Om_\rb$ is a $\DD^*$-valued skew-symmetric bicharacter on $\Zset^N$ arising from a multiplicatively skew-symmetric matrix $\rb \in M_N(\DD)$ as in \eqref{Omt},
\item the elements in the image of $M$ are linearly independent over $\DD$, and
\item $\Fract \DD\langle M(\Zset^N) \rangle = \FF$.
\end{itemize}
The matrix $\rb$ is uniquely reconstructed from the toric frame $M$, 
and will be denoted by $\rb(M)$. The elements $M(e_1), \dots, M(e_N)$ are called \emph{cluster variables}. 
Fix a subset $\ex \subset [1,N]$, to be called the set of {\em{exchangeable  
indices}}; the remaining indices, those in $[1,N] \backslash \ex$, will be called {\em{frozen}}. 

An integral $N \times \ex$ matrix $\wt{B}$ will be called an {\em{exchange matrix}} if its principal part (the $\ex\times\ex$ submatrix) is 
skew-symmetrizable. \red{If the principal part of $\wt{B}$ is skew-symmetric, then it is represented by a quiver 
whose vertices are labelled by the integers in $[1,N]$. For $j, k \in [1,N]$, there is a directed edge from the vertex $j$ to the vertex 
$k$ if and only if $(\wt{B})_{jk}>0$ and the number of such directed edges equals $(\wt{B})_{jk}$. In particular, the quiver has no edges 
between any pair of vertices in $[1,N] \backslash \ex$. 
}

A \emph{quantum seed} for $\FF$ (over $\DD$) is a pair $(M, \wt{B})$ consisting of a toric frame $M$ for $\FF$ 
and an exchange matrix $\wt{B}$ compatible with $\rb(M)$ in the sense that
\begin{align*}
&
\Om_{\rb(M)}(b^k, e_j) = 1, \; \; \forall\,  k \in \ex, \; j \in [1,N], \; k \neq j
\quad \mbox{and}
\\
&\Om_{\rb(M)} (b^k, e_k) \; \; 
\mbox{are not roots of unity}, \; \; \forall\,  k \in \ex,
\end{align*}
where $b^k$ denotes the $k$-th column of $\wt{B}$.  

The \emph{mutation in direction $k \in \ex$} of a quantum seed $(M,\wt{B})$ is the quantum seed $(\mu_k(M), \mu_k(\wt{B}))$ where $\mu_k(M)$ is described below and $\mu_k(\wt{B})$ is the $N\times\ex$ matrix $(b'_{ij})$ with entries
$$
b'_{ij} :=
\begin{cases}
- b_{ij} \,, &\mbox{if} \; \;
i=k \; \; \mbox{or} \; \; j=k
\\
b_{ij} + \frac{ |b_{ik}| b_{kj} + b_{ik} | b_{kj}|}{2} \,,
&\mbox{otherwise},
\end{cases}
$$
\cite{FZ1}. \red{If the principal part of $\wt{B}$ is skew-symmetric, then $\mu_k(\wt{B})$ has the same property and the pair of quivers 
corresponding to $\wt{B}$ and $\mu_k(\wt{B})$ are obtained from each other by quiver mutation at the vertex $k$, see \cite[\S\S 2.1 and 2.7]{FWZ}  
for details.}
Corresponding to the column $b^k$ of $\wt{B}$ are $\DD$-algebra automorphisms $\rho_{b^k,\pm}$ of $\FF$ such that
$$
\rho_{b^k,\ep}(ME_\ep(e_j)) = 
\begin{cases}
ME_\ep(e_k) + ME_\ep(e_k + \ep b^k), \quad &\text{if} \ j = k  \\
ME_\ep(e_j), \quad &\text{if} \ j \ne k,
\end{cases}
$$
\cite[Proposition 4.2]{BZ} and 
\cite[Lemma 2.8]{GYbig}, where $E_\ep = E^{\wt{B}}_\ep$ is the $N\times N$ matrix with entries
$$
(E_\ep)_{ij} =
\begin{cases}
\delta_{ij} \,, & \mbox{if} \; j \neq k \\
-1, & \mbox{if} \; i=j=k \\
\max(0, - \ep b_{ik}), & \mbox{if} \; i \neq j = k.
\end{cases}
$$
The toric frame $\mu_k(M)$ is defined as
$$
\mu_k(M) := \rho_{b^k,\ep} M E_\ep : \Zset^N \longrightarrow \FF.
$$
It is independent of the choice of $\ep$, and, paired with $\mu_k(\wt{B})$, forms a quantum seed over $\KK$ \cite[Proposition 2.9]{GYbig}. (See also \cite[Corollary 2.11]{GYbig}, and compare with \cite[Proposition 4.9]{BZ} for the uniparameter case.) By \cite[Proposition 2.9 and Eq.~(2.22)]{GYbig}, the entries of $\rb(\mu_k(M)) = \mu_k(\rb(M))$ are products of powers of the entries of $\rb(M)$, so $\rb(\mu_k(M)) \in M_N(\DD)$. It follows that $\mu_k(M)$ is a toric frame for $\FF$ over $\DD$, so that $(\mu_k(M), \mu_k(\wt{B}))$ is a quantum seed over $\DD$.

Fix a subset $\inv$ of the set $[1,N] \backslash \ex$ of frozen indices -- the corresponding 
cluster variables will be inverted. The {\em{quantum cluster algebra}} 
$\AA(M, \wt{B}, \inv)_\DD$ is the unital $\DD$-subalgebra of $\FF$ generated by the cluster
variables of all seeds obtained from $(M, \wt{B})$ by iterated mutations and by 
$\{M(e_k)^{-1} \mid k \in \inv \}$. To each quantum seed $(M, \wt{B})$ and choice of $\inv$, one associates the mixed quantum 
torus/quantum affine space algebra
\begin{equation}  \label{TMBtil}
\DD\TT_{(M, \wt{B},\inv)}:= \DD \lcor M(e_k)^{\pm 1}, \;  M(e_j) \mid k \in \ex \cup \inv, \; j \in [1,N] \backslash (\ex \cup \inv) \rcor
\subset \FF.
\end{equation}
The intersection of all such subalgebras of $\FF$ associated to all seeds that are obtained by iterated mutation 
from the seed $(M,\wt{B})$ is called the {\em{upper quantum cluster algebra}} of 
$(M,\wt{B})$ and is denoted by $\UU(M, \wt{B}, \inv)_\DD$. The corresponding \emph{Laurent Phenomenon} \cite[Theorem 2.15]{GYbig} says that 
\begin{equation}  \label{Lphenom}
\AA(M, \wt{B}, \inv)_\DD \subseteq \UU(M, \wt{B}, \inv)_\DD \,.
\end{equation}

If $\KK$ is the quotient field of $\DD$, then $\FF$ is also a division algebra over $\KK$, and the above constructions may be performed over $\KK$. The corresponding quantum cluster algebras over $\KK$ are just the $\KK$-subalgebras of $\FF$ generated by the quantum cluster algebras over $\DD$:
$$
\AA(M, \wt{B}, \inv)_\KK = \KK\cdot \AA(M, \wt{B}, \inv)_\DD \,.
$$

The uniparameter quantum cluster algebras of Berenstein and Zelevinsky \cite{BZ} come from the above axiomatics when the following two conditions 
are imposed:
\begin{enumerate}
\item The base ring is taken to be 
\[
\DD = \iAA = \Zset[q^{\pm 1/2}].
\]
So, $\DD^*= (\iAA)^* = \pm q^{\Zset/2}$.

\item
The toric frame of one seed (and thus of any seed) has a multiplicatively skew-symmetric matrix $\rb \in M_N(\DD)$ of the form
\[
\rb = (q^{m_{ij}/2})_{i,j=1}^N \quad \mbox{for some} \quad m_{ij} \in \Zset. 
\]
\end{enumerate}

\sectionnew{Quantum nilpotent algebras}
\label{qNilp}

Quantum nilpotent algebras are iterated skew polynomial algebras over a base field, which we take to be $\KK$ in this section. We use the standard notation $S[x;\theta,\de]$ for a \emph{skew polynomial ring}, or \emph{Ore extension}; it denotes a ring generated by a subring $S$ and an element $x$ satisfying $xs = \theta(s) x + \de(s)$ for all $s \in S$, where $\theta$ is a ring endomorphism of $S$ and $\de$ is a {\rm(}left\/{\rm)} $\theta$-derivation of $S$. The ring $S[x;\theta,\de]$ is a free left $S$-module, with the nonnegative powers of $x$ forming a basis. For all skew polynomial rings $S[x;\theta,\de]$ considered in this paper, we assume that $\theta$ is an \emph{automorphism} of $S$. Moreover, we work in the context of algebras over a commutative ring $\DD$, so our coefficient rings $S$ will be $\DD$-algebras, the maps $\theta$ will be $\DD$-algebra automorphisms, and the maps $\de$ will be $\DD$-linear $\theta$-derivations. Under these assumptions, $S[x;\theta,\de]$ is naturally a $\DD$-algebra. Throughout the present section, $\DD = \KK$.

\subsection{CGL extensions}
\label{CGLext}

We focus on iterated skew polynomial extensions
\begin{equation} 
\label{itOre}
R := \KK[x_1][x_2; \theta_2, \delta_2] \cdots [x_N; \theta_N, \delta_N],
\end{equation}
where $\KK[x_1] = \KK[x_1; \id_\KK, 0]$. Set
$$
R_k := \KK \langle x_1, \dots, x_k \rangle = \KK[x_1][x_2; \theta_2, \delta_2] \cdots [x_k; \theta_k, \delta_k] \qquad \text{for}\ k \in [0,N];
$$
in particular, $R_0 = \KK$.

\bde{CGL} An iterated skew polynomial extension \eqref{itOre}
is called a \emph{quantum nilpotent algebra} or a \emph{CGL extension} 
\cite[Definition 3.1]{LLR} if it is equipped with a rational action of a $\KK$-torus $\HH$ 
by $\KK$-algebra automorphisms such that:
\begin{enumerate}
\item[(i)] The elements $x_1, \ldots, x_N$ are $\HH$-eigenvectors.
\item[(ii)] For every $k \in [2,N]$, $\de_k$ is a locally nilpotent 
$\theta_k$-derivation of the algebra $R_{k-1}$. 
\item[(iii)] For every $k \in [1,N]$, there exists $h_k \in \HH$ such that 
$\theta_k = (h_k \cdot)|_{R_{k-1}}$ and the $h_k$-eigenvalue of $x_k$, to be denoted by $\la_k$, is not a root of unity.
\end{enumerate}

Conditions (i) and (iii) imply that 
$$
\theta_k(x_j) = \la_{kj} x_j \; \; \mbox{for some} \; \la_{kj} \in \kx, \; \; \forall\,  1 \le j < k \le N.
$$
We then set $\la_{kk} :=1$ and $\la_{jk} := \la_{kj}^{-1}$ for $j< k$.
This gives rise to a multiplicatively skew-symmetric 
matrix $\lab := (\la_{kj}) \in M_N(\kx)$ and the corresponding skew-symmetric bicharacter $\Om_\lab$ from \eqref{Omt}. The elements $h_k \in \HH$ interact with the skew derivations $\de_k$ as follows:
\begin{equation}  \label{hde}
(h_k\cdot) \circ \de_k = \la_k \de_k \circ (h_k\cdot), \quad \forall\, k \in [1,N],
\end{equation}
see \cite[Eq.~(3.15)]{GYbig}.

The \emph{length} of $R$ is $N$ and its  \emph{rank} is given by by 
\begin{equation}
\label{rk}
\rk(R) := \{ k \in [1,N] \mid \de_k = 0 \} \in \Zset_{ > 0}
\end{equation}
(cf.~\cite[Eq. (4.3)]{GY1}).
Denote the character group of the torus $\HH$ by $\xh$.  The action of 
$\HH$ on $R$ gives rise to an $\xh$-grading of $R$, and the $\HH$-eigenvectors 
are precisely the nonzero homogeneous elements with respect to this grading. 
The $\HH$-eigenvalue of a nonzero homogeneous element $u \in R$ will be denoted by $\chi_u$; this equals its $\xh$-degree relative to the $\xh$-grading.
\ede 

By \cite[Proposition 3.2, Theorem 3.7]{LLR}, every CGL extension $R$ is an 
$\HH$-UFD, meaning that each nonzero $\HH$-prime ideal of $R$ contains a homogeneous prime element. (A \emph{prime element} of a domain $R$ is a nonzero element $p \in R$ such that $Rp = pR$ -- i.e., $p$ is a \emph{normal element} of $R$ -- and the ring $R/Rp$ is a domain.) A recursive description of the sets of homogeneous prime elements 
of the intermediate algebras $R_k$ of a CGL extension $R$ was obtained in  
\cite{GY1}. To state the result, we require the standard predecessor and successor functions, $p = p_\eta$ and $s = s_\eta$, of a function $\eta : [1,N] \to \Zset$, defined as follows:
\begin{equation}
\label{pred.succ}
\begin{aligned}
p(k) &:= \max \{ j <k \mid \eta(j) = \eta(k) \}, \\
s(k) &:= \min \{ j > k \mid \eta(j) = \eta(k) \}, 
\end{aligned}
\end{equation}
where $\max \varnothing = -\infty$ and $\min \varnothing = +\infty$. Corresponding order functions $O_\pm : [1,N] \rightarrow \Znn$ are defined by 
\begin{equation}
\label{O-+}
\begin{aligned}
O_-(k) &:= \max \{ m \in \Znn \mid p^m(k) \ne -\infty \},  \\
O_+(k) &:= \max \{ m \in \Znn \mid s^m(k) \ne +\infty \}.
\end{aligned}
\end{equation}

\bth{1} \cite[Theorem 4.3]{GY1} Let $R$ be a CGL extension of length $N$ and rank $\rk(R)$ as 
in \eqref{itOre}. There exist a function $\eta : [1,N] \to \Zset$ 
whose range has cardinality $\rk(R)$ and elements
$$
c_k \in R_{k-1} \; \; \mbox{for all} \; \; k \in [2,N] \; \; 
\mbox{with} \; \; p(k) \neq - \infty
$$
such that the elements $y_1, \ldots, y_N \in R$, recursively defined by 
\begin{equation}
\label{y}
y_k := 
\begin{cases}
y_{p(k)} x_k - c_k \,, &\mbox{if} \; \;  p(k) \neq - \infty \\
x_k \,, & \mbox{if} \; \; p(k) = - \infty,  
\end{cases}
\end{equation}
are homogeneous and have the property that for every $k \in [1,N]$,
\begin{equation}
\label{prime-elem}
\{y_j \mid j \in [1,k] , \, s(j) > k \}
\end{equation}
is a list of the homogeneous prime elements of $R_k$ up to scalar multiples.

The elements $y_1, \ldots, y_N \in R$ with these properties are unique.
The function $\eta$ satisfying the above 
conditions is not unique, but the partition of $[1,N]$ into a disjoint 
union of the level sets of $\eta$ is uniquely determined by the presentation \eqref{itOre} of $R$, as are the predecessor and successor functions $p$ and $s$.
The function $p$ has the property that $p(k) = - \infty$
if and only if $\de_k =0$.
\eth

The statement of \thref{1} is upgraded as in \cite[Theorem 3.6 and following comments]{GYbig}.
In the setting of the theorem, the rank of $R$ is also given by
\begin{equation}
\label{rankR}
\rk(R) = |\{ j \in [1,N] \mid s(j) > N \}|
\end{equation}
\cite[Eq. (4.3)]{GY1}.

\bde{lterm}
Denote by $\prec$ the reverse lexicographic order on $\Znn^N$:
\begin{multline}
\label{prec}
(m'_1, \ldots, m'_N) \prec (m_1, \ldots, m_N)
\quad \mbox{iff there exists} \\
\mbox{$n \in [1,N]$ with $m'_n < m_n \,, \, m'_{n+1} = m_{n+1} \,, \, \ldots, \, m'_N = m_N\,$.}
\end{multline}

A CGL extension $R$ as in \eqref{itOre} has the $\KK$-basis 
$$
\{ x^f := x_1^{m_1} \cdots x_N^{m_N} \mid 
f = (m_1, \ldots, m_N)^T \in \Znn^N \}.
$$
We say that a nonzero element
$b \in R$ 
has \emph{leading term} $t x^f$ and \emph{leading coefficient} $t$ where $t \in \kx$ and $f \in \Znn^N$ 
if 
$$
b = t x^f + \sum_{g \in \Znn^N,\; g \prec f} t_g x^g
$$
for some $t_g \in \KK$, and we set $\lc(b) := t$ and $\lt(b) := t x^f$.
\ede

The leading terms of the prime elements $y_k$ in \thref{1} are given by 
\begin{equation}
\label{lt-y}
\lt(y_k) = x_{p^{O_-(k)}(k)} \cdots x_{p(k)} x_k, \; \; \forall\,  k \in [1,N].
\end{equation}

The leading terms of reverse-order monomials $x_N^{m_N} \cdots x_1^{m_1}$ involve symmetrization scalars in $\kx$ defined by
\begin{equation}
\label{Scrlab}
\Scr_\lab(f) := \prod_{1 \le j< k \le N} \la_{jk}^{- m_j m_k}, \; \; \forall\,  f = (m_1,\dots,m_N)^T \in \Zset^N.
\end{equation}
Namely,
\begin{equation}
\label{x-comm}
\lt( x_N^{m_N} \cdots x_1^{m_1} ) = \Scr_\lab( (m_1,\dots,m_N)^T ) x_1^{m_1} \cdots x_N^{m_N}, \; \; \forall\,  (m_1,\dots,m_N)^T \in \Zset^N.
\end{equation}

\subsection{Symmetric CGL extensions}
\label{symmCGLext}

Given an iterated skew polynomial extension $R$ as in \eqref{itOre}, denote its interval subalgebras
$$
R_{[j,k]} := \KK \langle x_i \mid j \le  i \le k \rangle, \; \; \forall\,  j,k \in [1,N];
$$
in particular, $R_{[j,k]} = \KK$ if $j \nleq k$.

\bde{symCGL} A CGL extension $R$ as in 
\deref{CGL} is called {\em{symmetric}} provided
\begin{enumerate}
\item[(i)] For all $1 \leq j < k \leq N$,
$$
\de_k(x_j) \in R_{[j+1, k-1]}.
$$
\item[(ii)] For all $j \in [1,N]$, there exists $h^\sy_j \in \HH$ 
such that 
$$
h^\sy_j \cdot x_k = \la_{kj}^{-1} x_k = \la_{jk} x_k, \; \; \forall\,  
k \in [j+1, N]
$$
and $h^\sy_j \cdot x_j = \la^\sy_j x_j$ for some $\la^\sy_j \in \kx$ which is not a root of unity.
\end{enumerate}

Under these conditions, $R$ has a CGL extension presentation with the variables $x_k$ in descending order:
\begin{equation}  \label{reverseitOre}
R = \KK[x_N] [x_{N-1}; \theta^*_{N-1}, \de^*_{N-1}] \cdots [x_1; \theta^*_1, \de^*_1],
\end{equation}
see \cite[Corollary 6.4]{GY1}.
\ede

\bpr{la} {\rm\cite[Proposition 5.8]{GYbig}}
Let $R$ be a symmetric CGL extension of length $N$.
If $l \in [1,N]$ with $O_+(l) = m > 0$, then
\begin{equation}
\label{la-eq}
\la^*_l = \la^*_{s(l)} = \cdots = \la^*_{s^{m-1}(l)} = 
\la^{-1}_{s(l)} = \la^{-1}_{s^2(l)}= \cdots = \la^{-1}_{s^m(l)} \,.
\end{equation}
\epr

\bde{XiN}
Define the following subset of the symmetric group $S_N$:
\begin{equation}
\label{tau}
\begin{aligned}
\Xi_N := \{ \sig \in S_N \mid 
\sig(k) &= \max \, \sig( [1,k-1]) +1 \; \;
\mbox{or} 
\\
\sig(k) &= \min \, \sig( [1,k-1]) - 1, 
\; \; \forall\,  k \in [2,N] \}.
\end{aligned}
\end{equation}
In other words, $\Xi_N$ consists of those $\sig \in S_N$ 
such that $\sig([1,k])$ is an interval for all $k \in [2,N]$. The following subset of $\Xi_N$ will also be needed:
\begin{equation}
\label{GaN}  
\begin{aligned}
\Ga_N &:= \{ \sig_{i,j} \mid 1 \leq i \leq j \leq N \}, \; \; \text{where}  \\
\sig_{i,j} &:= [i+1, \ldots, j, i, j+1, \ldots, N, i-1, i-2, \ldots, 1] . 
\end{aligned}
\end{equation}
\ede

If $R$ is a symmetric CGL extension of length $N$, then for each $\sig \in \Xi_N$ there is a CGL extension presentation
\begin{equation}
\label{tauOre}
R = \KK [x_{\sig(1)}] [x_{\sig(2)}; \theta''_{\sig(2)}, \de''_{\sig(2)}] 
\cdots [x_{\sig(N)}; \theta''_{\sig(N)}, \de''_{\sig(N)}],
\end{equation}
see \cite[Remark 6.5]{GY1}, \cite[Proposition 3.9]{GYbig}.
Moreover, if $1 \le i \le k \le N$, then the subalgebra $R_{[i,k]}$ of $R$ is a symmetric CGL extension, to which \thref{1} applies.  In the case $k = s^m(i)$ we have

\bpr{yjk} {\rm\cite[Theorem 5.1]{GYbig}} Assume that $R$ is a symmetric CGL extension of length $N$, and $i \in [1,N]$ and $m \in \Znn$ are such that $s^m(i) \in [1,N]$. Then there is a unique homogeneous prime element $y_{[i,s^m(i)]} \in R_{[i,s^m(i)]}$ such that

{\rm{(i)}} $y_{[i,s^m(i)]} \notin R_{[i,s^m(i) -1]}$ and $y_{[i,s^m(i)]} \notin R_{[i+1, s^m(i) ]}$.

{\rm{(ii)}} $\lt( y_{[i,s^m(i)]} ) = x_i x_{s(i)} \cdots x_{s^m(i)}$.
\epr
The elements $y_{[i,s^m(i)]} \in R$ will be called {\em{interval prime elements}}. 
Certain combinations of the homogeneous prime elements from \prref{yjk} play an important role in the mutation formulas for quantum cluster variables of symmetric CGL extensions. They are given  in the following theorem, where we denote
\begin{equation}
\label{e-interval}
\begin{aligned}
e_{[j, s^l(j)]} := \; &e_j + e_{s(j)} + \cdots + e_{s^l(j)} \in \Zset^N,  \\
 &\forall\,  j \in [1,N], \; l \in \Znn \; \; \text{such that} \; \; s^l(j) \in [1,N],
 \end{aligned}
\end{equation}
and where we set $y_{[s(i),i]} := 1$.

\bth{uismi}
{\rm\cite[Corollary 5.11]{GYbig}}
Assume that $R$ is a symmetric CGL extension of length $N$,
and $i \in [1,N]$ and $m \in \Zset_{>0}$ are such that $s^m(i) \in [1,N]$.
Then
\begin{equation}
\label{uuu}
u_{[i,s^m(i)]} := 
y_{[i, s^{m-1}(i)]} y_{[s(i), s^m(i)]} - \Om_\lab(e_i, e_{[s(i), s^{m-1}(i)]}) 
y_{[s(i), s^{m-1}(i)]} y_{[i,s^m(i)]} 
\end{equation}
is a nonzero homogeneous normal element of $R_{[i+1,s^m(i)-1]}$ 
which is not a multiple of $y_{[s(i), s^{m-1}(i)]}$ if $m \geq 2$. 
\eth

The form and properties of the elements $u_{[i,s^m(i)]}$ mainly enter into the proofs of the mutation formulas for symmetric CGL extensions. However, an explicit normalization condition for the leading coefficients of these elements is required; see \eqref{cond} and \prref{rescalecond}.

\subsection{Rescaling of generators}
\label{rescalegens}

Assume $R$ is a CGL extension of length $N$ as in \eqref{itOre}. Given scalars $t_1,\dots, t_N \in \kx$, one can rescale the generators $x_j$ of $R$ in the fashion
\begin{equation}
\label{rescalexj}
x_j \longmapsto t_j x_j, \; \; \forall\,  j \in [1,N],
\end{equation}
meaning that $R$ is also an iterated Ore extension with generators $t_j x_j$; in fact,
\begin{equation}
\label{rescaleitOre}
R = \KK[t_1 x_1][t_2 x_2; \theta_2, t_2 \delta_2] \cdots [t_N x_N; \theta_N, t_N \delta_N].
\end{equation}
This is also a CGL extension presentation of $R$, and if \eqref{itOre} is a symmetric CGL extension, then so is \eqref{rescaleitOre}. 

Rescaling as in \eqref{rescalexj}, \eqref{rescaleitOre} does not affect the $\HH$-action or the matrix $\lab$, but various elements computed in terms of the new generators are correspondingly rescaled, such as the homogeneous prime elements from \thref{1} and \prref{yjk}. These transform as follows:
\begin{equation}
\label{rescaleyelems}
y_k \longmapsto \biggl( \, \prod_{l = 0}^{O_-(k)} t_{p^l(k)} \biggr) y_k \quad \text{and} \quad y_{[i,s^m(i)]} \longmapsto \biggl( \, \prod_{l=0}^m t_{s^l(i)} \biggr) y_{[i,s^m(i)]} \,.
\end{equation}
Consequently, the homogeneous normal elements \eqref{uuu} transform via
\begin{equation}
\label{rescaleuelems}
u_{[i,s^m(i)]} \longmapsto \bigl( t_i t_{s(i)}^2 \cdots t_{s^{m-1}(i)}^2 t_{s^m(i)} \bigr) u_{[i,s^m(i)]} \,.
\end{equation}

\subsection{Normalization conditions}
\label{normalizations}

In order for the homogeneous prime elements $y_k$ from \thref{1} to function as quantum cluster variables,  some normalizations are required. Throughout this subsection, assume that $R$ is a symmetric CGL extension of length $N$ as in Definitions \ref{dCGL} and \ref{dsymCGL}. Assume also that the following mild conditions on scalars are satisfied:

{\bf{Condition (A).}} The base field $\KK$ contains square roots $\nu_{kl} = \sqrt{\la_{kl}}$ of the scalars
$\la_{kl}$ for $1 \leq l < k \leq N$ such that the subgroup of $\kx$ generated 
by the $\nu_{kl}$ contains no elements of order $2$. Then set $\nu_{kk} := 1$ and $\nu_{kl} := \nu_{lk}^{-1}$ for $k < l$, so that $\nub := (\nu_{kl})$ is a multiplicatively skew-symmetric matrix.

{\bf{Condition (B).}} There exist positive integers $d_n$, for $n \in \eta([1,N])$,
such that 
$$
\la_k^{d_{\eta(l)}} = \la_l^{d_{\eta(k)}} , \quad \forall\, k, l \in [1,N]\; \text{with}\; p(k), p(l) \neq - \infty.
$$
In view of \prref{la}, this is equivalent to the condition
$$
(\la^*_k)^{d_{\eta(l)}} = (\la^*_l)^{d_{\eta(k)}}, \quad \forall\,  k,l \in [1,N] \; \text{with} \; s(k), s(l) \ne +\infty.
$$
\bre{sat-B} Note that Condition (B) is always satisfied if all $\la_k = q^{m_k}$ for some $m_k \in \Zset$ and 
$q \in \KK$ (which has to be a non-root of unity due to assumption (iii) in \deref{CGL}). This is the setting of Theorem A in the introduction.
\ere

In parallel with \eqref{Scrlab}, define
\begin{equation}
\label{Scrnub}
\Scr_\nub(f) := \prod_{1\le j < k \le N} \nu_{jk}^{-m_j m_k}, \; \; \forall\,  f = (m_1,\dots, m_N)^T \in \Zset^N.
\end{equation}
Then set
\begin{equation}
\label{eybarj}
\ol{e}_j := e_j + e_{p(j)} + \cdots + e_{p^{O_-(j)}(j)} \quad \text{and} \quad  \ol{y}_j := \Scr_\nub(\ol{e}_j) y_j \,, \; \; \forall\,  j \in [1,N].
\end{equation}

We analogously normalize the homogeneous prime elements described in \prref{yjk}:
\begin{equation}
\label{ybarismi}
\begin{aligned}
\ol{y}_{[i,s^m(i)]} := \; &\Scr_\nub( e_{[i,s^m(i)]} ) y_{[i,s^m(i)]} \,, \\ &\forall\,  i \in [1,N], \; m \in \Znn \; \; \text{such that} \; \; s^m(i) \in [1,N].
\end{aligned}
\end{equation}

A final normalization, for the leading coefficients of the elements $u_{[i,s^m(i)]}$, is needed in order to establish mutation formulas for the quantum cluster variables $\ol{y}_k$. For $i \in [1,N]$ and $m \in \Znn$ such that $s^m(i) \in [1,N]$, write
\begin{equation}
\label{ltu}
\begin{aligned}
\lt( u_{[i,s^m(i)]} ) =\;  &\pi_{[i,s^m(i)]} x^{f_{[i,s^m(i)]}} \,,  \\
  &\pi_{[i,s^m(i)]} \in \kx, \; \; f_{[i,s^m(i)]} \in \sum_{j=1+1}^{s^m(i)-1} \Znn \, e_j \subset \Znn^N \,.
  \end{aligned}
\end{equation}
We will require the condition
\begin{equation}
\label{cond}
\pi_{[i,s(i)]} = \Scr_\nub( -e_i + f_{[i,s(i)]} ), \; \; \forall\,  i \in [1,N] \; \; \text{such that} \; \; s(i) \ne +\infty.
\end{equation} 
This can always be achieved after a suitable rescaling of the $x_j$, as follows.

\bpr{rescalecond}
{\rm\cite[Propositions 6.3, 6.1]{GYbig}}
Let $R$ be a symmetric CGL extension of length $N$, satisfying condition {\rm(A)}. 

{\rm (i)} There exist $N$-tuples $(t_1,\dots, t_N) \in (\kx)^N$ such that after the rescaling \eqref{rescalexj}, condition \eqref{cond} holds.

{\rm (ii)} If \eqref{cond} holds, then
$$
\pi_{[i,s^m(i)]} = \;  \Scr_\nub( e_{[s(i),s^m(i)]} )^{-2} \Scr_\nub( -e_i + f_{[i,s^m(i)]} )
$$
for all $i \in [1,N]$, $m \in \Znn$ with $s^m(i) \in [1,N]$.
\epr

\subsection{Quantum cluster algebra structures for symmetric CGL extensions}
\label{thm8.2GYbig}

We present in this subsection the main theorem from \cite{GYbig}. 

Recall the notation on quantum cluster algebras from Section \ref{q-cl}. There is a right action 
of $S_N$ on the set of toric frames for a division algebra $\FF$, given by re-enumeration,
\begin{equation}
\label{r-act}
(M \cdot \tau) (e_k) := M(e_{\tau(k)}), \; \;  \rb(M \cdot \tau)_{jk} = \rb(M)_{\tau(j), \tau(k)} \,, \quad \tau \in S_N \,, \;\; j,k \in [1,N].
\end{equation}

Fix a symmetric CGL extension $R$ of length $N$ such that Conditions (A) and (B) hold. Define the multiplicatively skew-symmetric matrix $\nub$ as in Condition (A),  with associated bicharacter $\Om_\nub$ as in \eqref{Omt}, and define a second multiplicatively skew-symmetric matrix $\rb = (r_{kj})$ by
\begin{equation}
\label{defrb}
 r_{kj} := \Om_\nub( \ol{e}_k, \ol{e}_j ), \; \; \forall\,  k,j \in [1,N].
\end{equation}
Let $\ol{y}_1, \ldots, \ol{y}_N$ be the sequence of normalized homogeneous prime elements 
given in \eqref{eybarj}. (We recall that 
each of these is a prime element in some of 
the subalgebras $R_l$, but not necessarily in the full algebra $R=R_N$.) There is  a unique toric frame $M : \Zset^N \to \Fract(R)$ whose matrix is $\rb(M) := \rb$ 
and such that $M(e_k) := \ol{y}_k$ for all $k \in [1,N]$ \cite[Proposition 4.6]{GYbig}.

Next, consider an arbitrary element $\sig \in \Xi_N \subset S_N$, recall \eqref{tau}. For any $k \in [1,N]$, we see that
\begin{equation}
\label{tau-int}
\eta^{-1} \eta \sig(k) \cap \sig([1,k])
= \begin{cases}
\{p^n(\sig(k)), \ldots, p(\sig(k)), \sig(k)\}, & \mbox{if} \; \; \sig(1) \le \sig(k)
\\ 
\{\sig(k), s(\sig(k)), \ldots, s^n(\sig(k))\}, & \mbox{if} \; \; \sig(1) \ge \sig(k)
\end{cases}
\end{equation}
for some $n \in \Znn$. Corresponding to $\sig$, we have the CGL extension presentation \eqref{tauOre}, whose $\lab$-matrix is  the matrix $\lab_\sig$ with entries $(\lab_\sig)_{ij} := \lab_{\sig(i) \sig(j)}$. Analogously we define the 
matrix $\nub_\sig$, and denote by $\rb_\sig$ the corresponding multiplicatively skew-symmetric 
matrix derived from $\nub_\sig$ by applying \eqref{defrb} to the presentation  \eqref{tauOre}. Explicitly,
\begin{equation}
\label{tau-frame1}
(\rb_\sig)_{kj} = \prod \{ \nu_{il} \mid i \in \sig([1,k]),\; \eta(i) = \eta \sig(k), \;
l \in \sig([1,j]), \; \eta(l) = \eta \sig (j) \},
\end{equation}
cf.~\eqref{tau-int}.
Let $\ol{y}_{\sig,1}, \ldots, \ol{y}_{\sig,N}$ 
be the sequence of normalized prime elements given by applying \eqref{eybarj}
 to the presentation \eqref{tauOre}. By \cite[Proposition 4.6]{GYbig}, 
there is a unique toric frame $M_\sig : \Zset^N \to \Fract(R)$ whose matrix 
is $\rb(M_\sig) := \rb_\sig$ and such that for all $k \in [1,N]$,
\begin{equation}
\label{tau-frame2}
M_\sig(e_k) := \ol{y}_{\sig, k} =
\begin{cases} 
\ol{y}_{[p^n(\sig(k)), \sig(k)]}, & \mbox{if} \; \; \sig(1) \le \sig(k) 
\\ 
\ol{y}_{[\sig(k), s^n(\sig(k))]}, & \mbox{if} \; \; \sig(1) \ge \sig(k)
\end{cases}
\end{equation}
in the two cases of \eqref{tau-int}, respectively. The last equality 
is proved in \cite[Theorem 5.2]{GYbig}.

Recall that the set $P(N) := \{ k \in [1,N] \mid s(k) = + \infty \}$
parametrizes the set of homogeneous prime elements of $R$,
i.e.,  
$$
\{y_k \mid k \in P(N) \} \; \; \mbox{is a list of the homogeneous prime
elements of $R$}
$$ 
up to scalar multiples (\thref{1}). Define
$$
\ex := [1,N] \setminus P(N) = \{ l \in [1,N] \mid s(l) \neq + \infty \}.
$$
Since $|P(N)| = \rk(R)$,
the cardinality of the set $\ex$ is $N - \rk(R)$. 
For $\sig \in \Xi_N$, define the set
$$
\ex_\sig = \{ l \in [1, N] \mid \exists\, k > l \;\; \text{with} \;\; \eta \sig(k) = \eta \sig (l) \}
$$
of the same cardinality. Finally, recall the notation $\chi_u$ from \deref{CGL}.

In \cite[Theorem 8.2]{GYbig} we re-indexed all toric frames $M_\sig$ in such a way that 
the right action in \thref{maint} (c) was trivialized and the exchangeable variables in all 
such seeds were parametrized just by $\ex$, rather than by $\ex_\sig$. We omit the re-indexing here,
to simply the exposition. This affects the upper cluster algebra in the following way: When considering the quantum seed $(M_\sig,\wt{B}_\sig)$, the set $\ex$ must be replaced by $\ex_\sig$ in relations such as \eqref{TMBtil}.

\bth{maint} {\rm\cite[Theorem 8.2]{GYbig}}
Let $R$ be a symmetric CGL extension of length $N$ and rank $\rk(R)$
as in Definitions {\rm\ref{dCGL} \red{\it and\/} \ref{dsymCGL}}. Assume that Conditions {\rm (A), (B)} hold, and that the sequence of 
generators $x_1, \ldots, x_N$ of $R$ is normalized {\rm(}rescaled\/{\rm)}
so that condition \eqref{cond} is satisfied. 
Then the following hold:

{\rm(a)} For all $\sig \in \Xi_N$ {\rm(}see \eqref{tau}{\rm)} and $l \in \ex_\sig$, 
there exists a unique vector $b_\sig^l \in \Zset^N$ such that $\chi_{M_\sig(b_\sig^l)}=1$ and 
\begin{equation}
\label{linear-eq}
\Om_{\rb_\sig} ( b^l_\sig, e_j) = 1, \; \; \forall\,  j \in [1,N], \; j \neq l
\quad \mbox{and} \quad
\Om_{\rb_\sig} (b^l_\sig, e_l)^2 = \la^*_{ \min \eta^{-1} \eta (\sig(l))} \,.
\end{equation}
Denote by $\wt{B}_\sig \in M_{N \times |\ex| }(\Zset)$ the matrix with columns 
$b^l_\sig$, $l \in \ex_\sig$. Let $\wt{B}:= \wt{B}_\id$.
 
{\rm(b)} For all $\sig \in \Xi_N$, the pair $(M_\sig, \wt{B}_\sig)$ is a 
quantum seed for $\Fract(R)$. The principal part of $\wt{B}_\sig$ 
is skew-symmetrizable via the integers $d_{\eta(k)}$, $k \in \ex_\sig$ from Condition {\rm (B)}.

{\rm(c)} All such quantum seeds are mutation-equivalent to each other up to the $S_N$ action.
They are linked by the following one-step mutations.
Let $\sig, \sig' \in \Xi_N$ be such that
$$
\sig' = ( \sig(k), \sig(k+1)) \circ \sig = \sig \circ (k, k+1)
$$
for some $k \in [1,N-1]$.
If $\eta(\sig(k)) \neq \eta (\sig(k+1))$, then $M_{\sig'} = M_\sig \cdot (k,k+1)$
in terms of the action \eqref{r-act}.
If $\eta(\sig(k)) = \eta (\sig(k+1))$, 
then $M_{\sig'} = \mu_k (M_\sig)  \cdot (k,k+1)$.

{\rm(d)} The CGL extension $R$ equals the quantum cluster and upper cluster algebras associated to $M$, $\wt{B}$, $\varnothing$:
$$
R = \AA(M, \wt{B}, \varnothing)_\KK = \UU(M, \wt{B}, \varnothing)_\KK \,.
$$
In particular, $\AA(M, \wt{B}, \varnothing)_\KK$ and $\UU(M, \wt{B}, \varnothing)_\KK$ 
are affine and noetherian, and more precisely $\AA(M, \wt{B}, \varnothing)_\KK$ is 
generated by the cluster variables in the seeds parametrized by the finite subset 
$\Ga_N$ of $\Xi_N$, recall \eqref{GaN}.

{\rm(e)} Let $\inv$ be any subset of the set $P(N)$ of frozen variables.
Then 
$$
R[y_k^{-1} \mid k \in \inv] = \AA(M, \wt{B}, \inv)_\KK = \UU(M, \wt{B}, \inv)_\KK \,.
$$
\eth 

\sectionnew{Integral quantum cluster structures on quantum nilpotent algebras}
\label{Int-qNilp}

We introduce integral forms of CGL extensions and show that the quantum cluster algebra structure on a symmetric CGL extension $R$ satisfying the hypotheses of \thref{maint} passes to appropriate integral forms of $R$.

Throughout the section, let $R$ be a CGL extension of length $N$ as in \deref{CGL}, with associated torus $\HH$, scalars $\la_{kj}$ and $\la_k$, and other notation as in Section \ref{qNilp}. Let $\DD \subseteq \KK$ be a unital subring of $\KK$, and write $\DD^*$ for the group of units of $\DD$.

\subsection{Integral forms of CGL extensions}

\bde{DformCGL}
We say that the $\DD$-subalgebra $\DD \langle x_1, \dots, x_N \rangle$ of $R$ is a \emph{$\DD$-form} of the CGL presentation \eqref{itOre} -- and therefore that \eqref{itOre} \emph{has a $\DD$-form} -- provided this subalgebra is an iterated skew polynomial extension of the form
\begin{equation}  \label{DitOre}
\DD \langle x_1, \dots, x_N \rangle = \DD[x_1][x_2; \theta_2, \delta_2] \cdots [x_N; \theta_N, \delta_N],
\end{equation}
where we let $\theta_k$ (resp., $\de_k$) also denote the restriction of the original $\theta_k$ (resp., $\de_k$) to a $\DD$-algebra automorphism (resp., $\theta_k$-derivation) of $\DD\langle x_1,\dots,x_{k-1} \rangle$.
\ede

\bre{Dformcomments}
(a) The CGL presentation \eqref{itOre} has a $\DD$-form if and only if
\begin{itemize}
\item $\la_{kj} \in \DD^*$ for $1\le j< k\le N$;
\item $\de_k$ maps $\DD\langle x_1,\dots,x_{k-1} \rangle$ into itself for each $k \in [2,N]$.
\end{itemize}

(b) Whether \eqref{itOre} has a $\DD$-form depends on the choice of $\DD$ as well as the choice of CGL presentation \eqref{itOre}. For instance, if $N=2$ and $\de_2(x_1) \in \KK \setminus \DD$, then \eqref{itOre} does not have a $\DD$-form. However, if $\ga = \de_2(x_1)$, then $R$ has the CGL presentation $\KK[x_1] [\ga^{-1}x_2; \sig_2, \ga^{-1}\de_2]$, which does have a $\DD$-form.

(c) Even if \eqref{itOre} has a $\DD$-form,  the homogeneous prime elements $y_1,\dots,y_N$ from \thref{1} need not belong to $\DD \langle x_1, \dots, x_N \rangle$. For instance, if $R$ is the quantized Weyl algebra 
$$
A_1^q(\KK) = \KK \langle x_1,x_2 \mid x_1 x_2 = q x_2 x_1 + 1 \rangle
$$
with $q\in \kx$ transcendental over the prime field of $\KK$ and $\DD = (\Zset {\cdot} 1_\KK)[q^{\pm1}]$, then the above CGL presentation has a $\DD$-form, but $\DD\langle x_1,x_2 \rangle$ does not contain the element $y_2 = x_1 x_2 + (q-1)^{-1}$.
\ere

The problems indicated in \reref{Dformcomments} can typically be corrected by rescaling the generators $x_k$ as in \S\ref{rescalegens}, as we now show.

When working with a $\DD$-subalgebra $R' = \DD \langle x_1, \dots, x_N \rangle$ of $R$, we adapt previous notation and write
$$
R'_k := \DD \langle x_1,\dots,x_k \rangle \quad\text{and}\quad R'_{[j,k]} := \DD \langle x_j,\dots,x_k \rangle, \quad \forall\, j \le k \in [1,N].
$$

\bpr{rescale.Dform}
Assume that $\KK = \Fract \DD$ and that $\la_{kj} \in \DD^*$ for $1\le j< k\le N$. Then there exist $t_1,\dots,t_N \in \DD\setminus\{0\}$ such that

{\rm(a)} $R_\DD := \DD\langle t_1x_1, \dots, t_Nx_N \rangle$ is a $\DD$-form of the CGL presentation \eqref{rescaleitOre}.

{\rm(b)} The elements $y_1,\dots,y_N$ from Theorem {\rm\ref{t1}} for the presentation \eqref{rescaleitOre} all lie in $R_\DD$.
\epr

\begin{proof} Set $R' := R_\DD$ for the proof. We induct on $N$. The case $N=1$ holds trivially by taking $t_1=1$.

Now assume that $N>1$ and that there exist $t_1,\dots,t_{N-1} \in \DD\setminus\{0\}$ such that the algebra $R'_{N-1} :=  \DD\langle t_1x_1, \dots, t_{N-1}x_{N-1} \rangle$ satisfies condtions (a), (b). In particular, $R'_{N-1}$ is a $\DD$-form of the CGL presentation
\begin{equation}  \label{rescaleRN-1}
R_{N-1} =  \KK[t_1 x_1][t_2 x_2; \theta_2, t_2 \delta_2] \cdots [t_{N-1} x_{N-1}; \theta_{N-1}, t_{N-1} \delta_{N-1}].
\end{equation}
 Since $\la_{Nj}^{\pm1} \in \DD$ for all $j \in [1,N-1]$, the automorphism $\theta_N$ restricts to an automorphism of $R'_{N-1}$.

Write $\de_N(x_1), \dots, \de_N(x_{N-1})$ as $\KK$-linear combinations of monomials 
$$
(t_1x_1)^{m_1} \cdots (t_{N-1} x_{N-1})^{m_{N-1}}
$$
in the standard PBW basis for the presentation \eqref{rescaleRN-1}, and let $\kappa_i$ for $i\in I$ be a list of the nonzero coefficients that appear. Choose a nonzero element $b\in \DD$ such that $b \kappa_i \in \DD$ for all $i$. Set
$$
t_N := \begin{cases}
b  &(\text{if}\ p(N) = -\infty) \\
(\la_N-1)b  &(\text{if}\ p(N) \ne -\infty).
\end{cases}
$$
Since $b \kappa_i \in \DD$ for all $i$, we have $b \de_N(x_j) \in R'_{N-1}$ for all $j \in [1,N-1]$, and so $b \de_N$ maps $R'_{N-1}$ into itself. Then also $t_N \de_N$ maps $R'_{N-1}$ into itself. Therefore $R' = R'_{N-1} \langle t_N x_N \rangle$ is an Ore extension $R'_{N-1}[ t_N x_N; \theta_N, t_N \de_N]$ and (a) holds.

It remains to show that the element $y_N$ for the CGL presentation \eqref{rescaleitOre} lies in $R'$. If $p(N) = -\infty$, then $y_N = t_N x_N$ and we are done. Now assume that $p(N) \ne -\infty$. Then $y_N = y_{p(N)} x_N - c_N$ where $c_N \in R_N$ and $y_{p(N)}$ is the $p(N)$-th $y$-element for \eqref{rescaleitOre}. By our induction hypotheses, $y_{p(N)} \in R'_{N-1}$. From \cite[Proposition 4.7(b)]{GY1}, we have
$$
(\la_N-1)b \de_N( y_{p(N)} ) = t_N \de_N( y_{p(N)} ) = \prod_{m=1}^{O_-(N)} \la_{N,p^m(N)} (\la_N - 1) c_N \,.
$$
Since $\la_{N,p^m(N)} \in \DD^*$ for all $m \in [1, O_-(N)]$ and $b \de_N( y_{p(N)} ) \in R'_{N-1}$, we conclude that $c_N \in R'_{N-1}$. Therefore $y_N \in R'$, as required.
\end{proof}

\ble{Dforms1}
If the CGL presentation \eqref{itOre} has a $\DD$-form, then $\la_k \in \DD^*$ for all $k \in [2,N]$ such that $p(k) \ne -\infty$.
\ele

\begin{proof} If $k \in [2,N]$ and $p(k) \ne -\infty$, then $\de_k \ne 0$ (recall \thref{1}). Choose $i \in [1,k-1]$ such that $\de_k(x_i) \ne 0$, and choose a monomial $x^f$, for some $f = (m_1,\dots,m_{k-1})^T \in \Znn^{k-1}$, which appears with a nonzero coefficient in $\de_k(x_i)$. In view of \eqref{hde}, $h_k.\de_k(x_i) = \la_k \la_{ki} \de_k(x_i)$. Since all monomials in $x_1,\dots,x_N$ are $h_k$-eigenvectors, it follows that $h_k.x^f = \la_k \la_{ki} x^f$. On the other hand, $h_k.x^f = \theta_k(x^f) = \prod_{j=1}^{k-1} \la_{kj}^{m_j} x^f$, and consequently
$$
\la_k = \la_{ki}^{-1}  \prod_{j=1}^{k-1} \la_{kj}^{m_j} \in \DD^* \,,
$$
since all $\la_{kj} \in \DD^*$ (\reref{Dformcomments}(a)).
\end{proof}

 In case $R$ is symmetric and \eqref{itOre} has a $\DD$-form, the alternative CGL extension presentations of $R$ given in \eqref{tauOre} also have $\DD$-forms, as we now show. 

\ble{Dforms2}
Assume that $R_\DD = \DD \langle x_1,\dots, x_N \rangle$ is a $\DD$-form for \eqref{itOre}, and that $R$ is a symmetric CGL extension.

{\rm(a)} For $1\le j < k \le N$, the algebra $(R_\DD)_{[j,k]}$ is a $\DD$-form for the CGL presentation
\begin{equation}  \label{itOreRjk}
R_{[j,k]} = \KK [x_j] [x_{j+1}; \theta_{j+1}, \de_{j+1}] \cdots [x_k; \theta_k, \de_k].
\end{equation}

{\rm(b)} For each $\sig \in \Xi_N$, the algebra $R_\DD$ is a $\DD$-form for the CGL  presentation \eqref{tauOre} of $R$.
\ele

\begin{proof} Set $R' := R_\DD$.

(a) The symmetry assumption on $R$ implies that the $\KK$-subalgebra $R_{[j,k]}$ of $R$ is itself a CGL extension of the form \eqref{itOreRjk}, as noted following \deref{XiN}. 
For $l \in [j+1,k]$, closure of both $R_{[j,l-1]}$ and $R'_{l-1}$ under $\theta_l^{\pm1}$ and $\de_l$ implies that $R'_{[j,l-1]}$ is closed under $\theta_l^{\pm1}$ and $\de_l$. It follows that $R'_{[j,k]}$ is an iterated Ore extension of the form
$\DD [x_j] \cdots [x_k; \theta_k, \de_k]$, as required.

(b) We first consider the reverse CGL extension presentation \eqref{reverseitOre}. As shown in \cite[\S6.2]{GY1} (where $\theta_j^*$, $\de_j^*$ are denoted $\sig_j'$, $\de_j'$), we have
$$
\theta^*_j(x_k) = \la_{jk} x_k \quad \text{and} \quad \de^*_j(x_k) = - \la_{jk} \de_k(x_j), \quad \forall\, 1\le j < k \le N.
$$
Consequently, $\DD\langle x_{j+1}, \dots, x_N \rangle$ is stable under $(\theta^*_j)^{\pm1}$ and $\de^*_j$ for each $j \in [1,N-1]$. This allows us to write $R'$ as an iterated Ore extension in the form
\begin{equation}  \label{reverseRD}
R' = \DD[x_N] [x_{N-1}; \theta^*_{N-1}, \de^*_{N-1}] \cdots [x_1; \theta^*_1, \de^*_1],
\end{equation}
which shows that $R'$ is a $\DD$-form for \eqref{reverseitOre}.

Now let $\sig$ be an arbitrary element of $\Xi_N$ and consider the corresponding CGL extension presentation \eqref{tauOre} of $R$. As indicated in \cite[Remark 6.5]{GY1}, the automorphisms $\theta_j''$ and skew derivations $\de_j''$ appearing in \eqref{tauOre} are restrictions of either $\theta_j$, $\de_j$ or $\theta^*_j$, $\de^*_j$. Combined with the results of the previous paragraph, we conclude that $R'$
is an iterated Ore extension of the form
$$
\DD [x_{\sig(1)}] [x_{\sig(2)}; \theta''_{\sig(2)}, \de''_{\sig(2)}] 
\cdots [x_{\sig(N)}; \theta''_{\sig(N)}, \de''_{\sig(N)}].
$$
Therefore $R'$ is a $\DD$-form for \eqref{tauOre}. 
\end{proof}

\ble{Dforms3}
Assume that the CGL presentation \eqref{itOre} has a $\DD$-form $R_\DD = \DD \langle x_1,\dots, x_N \rangle$ which contains the elements $y_1,\dots,y_N$ from Theorem {\rm\ref{t1}}.

{\rm(a)} For each $k \in [1,N]$, the element $y_k$ is normal in $(R_\DD)_k$.

{\rm(b)} For any subset $I \subseteq [1,N]$, the multiplicative set generated by $\DD^* \cup \{ y_i \mid i \in I \}$ is a denominator set in $R_\DD$.
\ele

\begin{proof} Set $R' := R_\DD$.

(a) By \cite[Corollary 4.8]{GY1}, $y_k$ quasi-commutes with those $x_j$ such that $j < s(k)$ according to the rule
$$
y_k x_j = \biggl(\, \prod_{m=0}^{O_-(k)} \la_{j,p^m(k)} \biggr)^{-1} x_j y_k \,.
$$
Since the $\la_{j,p^m(k)}$ all lie in $\DD^*$, it follows that $y_k R'_k = R'_k y_k$.

(b) It suffices to show that $\DD^* y_k^\Nset$, the multiplicative set generated by $\DD^* \cup \{ y_k \}$, is a denominator set in $R'$ for each $k \in [1,N]$. By part (a), $\DD^* y_k^\Nset$ is a denominator set in $R'_k$.

Since $y_k$ is homogeneous (with respect to the $\xh$-grading on $R$), it is an eigenvector for each $h \in \HH$ and thus for $\theta_{k+1}, \dots, \theta_N$. The leading term of $y_k$ is $x_{p^{O_-(k)}(k)} \cdots x_{p(k)} x_k$, and so
$$
\theta_l(y_k) = \biggl(\, \prod_{m=0}^{O_-(k)} \la_{l,p^m(k)}\biggr) y_k \,, \quad \text{for}\ \ 1\le k < l \le N.
$$
Consequently, $\theta_l(\DD^* y_k^\Nset) = \DD^* y_k^\Nset$ for all $l > k$. It therefore follows from \cite[Lemma 1.4]{Go}, by induction on $l$, that $\DD^* y_k^\Nset$ is a denominator set in $R'_l$ for $l=k+1,\dots,N$.
\end{proof}

\bpr{tauDform}
Assume that $R$ is a symmetric CGL extension and that the CGL presentation \eqref{itOre} has a $\DD$-form $R_\DD = \DD \langle x_1,\dots, x_N \rangle$ which contains the elements $y_1,\dots,y_N$.

{\rm(a)} The elements $y_{[i,s^m(i)]}$ of Proposition {\rm\ref{pyjk}} all belong to $R_\DD$.

{\rm(b)} The elements $u_{[i,s^m(i)]}$ of Theorem {\rm\ref{tuismi}} all belong to $R_\DD$, and their leading coefficients $\pi_{[i,s^m(i)]}$ belong to $\DD$.

{\rm(c)} The elements $y_{\sig,k}$, for $\sig \in \Xi_N$ and $k \in [1,N]$, all belong to $R_\DD$.
\epr

\begin{proof} Set $R' := R_\DD$.

(a) We first recall that by the case $\tau=\id$ of \cite[Theorem 5.3]{GYbig}, $y_k$ is a scalar multiple of $y_{[p^{O_-(k)}(k),k]}$ for all $k \in [1,N]$. However, these elements both have leading coefficient $1$, so they are equal. Taking $k = s^m(i)$, we obtain
\begin{equation}  \label{yint=y}
y_{[i,s^m(i)]} = y_{s^m(i)} \,, \quad \forall\, i \in [1,N] \ \; \text{with} \ \; p(i) = -\infty.
\end{equation}
This verifies that $y_{[i,s^m(i)]} \in R'$ whenever $p(i) = -\infty$.

We next show, by induction on $i$, that all $y_{[i,s^m(i)]} \in R'$. The case $i=1$ follows from the previous result, since $p(1) = -\infty$. Now assume that $i>1$ and that $y_{[j,s^m(j)]} \in R'$ for all $j \in [1,i-1]$ and $m \in [0,O_+(j)]$. If $p(i) = -\infty$, we are done by the previous result, so we may assume that $p(i)  = j \in [1,i-1]$. Set $k = s^m(i) = s^{m+1}(j)$. By the induction hypothesis, $y_{[j,k]} \in R'$. According to \cite[Theorem 5.1(d)]{GYbig},
$$
y_{[j,k]} = x_j y_{[i,k]} - c'
$$
for some $c' \in R_{[j+1,k]}$. Since $R_{[j,k]}$ (resp., $R'_{[j,k]}$) is a free right module over $R_{[j+1,k]}$ (resp., $R'_{[j+1,k]}$) with basis $\{1,x_j,x_j^2, \dots\}$, the assumption $y_{[j,k]} \in R'$ implies $y_{[i,k]} \in R'$. This concludes the induction step. 

(b) Since all values of the bicharacter $\Om_\lab$ lie in $\DD^*$, the formula \eqref{uuu} together with part (a) yields $u_{[i,s^m(i)]} \in R'$. Consequently, its leading coefficient, $\pi_{[i,s^m(i)]}$, must lie in $\DD$.

(c) Fix $\sig \in \Xi_N$. We proceed by induction on $k \in [1,N]$ to show that $y_{\sig,k} \in R'$. The case $k=1$ holds trivially, since $y_{\sig,1} = x_{\sig(1)}$. Now let $k>1$ and assume that $y_{\sig,j} \in R'$ for all $j \in [1,k-1]$.

If $p(\sig(k)) \notin \sig([1,k-1])$, then $y_{\sig,k} = x_{\sig(k)}$ and we are done. Assume now that $p(\sig(k)) = \sig(l)$ for some $l \in [1,k-1]$. Then $y_{\sig,l} \in R'$ by induction, and
\begin{equation}  \label{ysigk}
y_{\sig,k} = y_{\sig,l} x_{\sig(k)} - c, \quad 0 \ne c \in R_{\sig([1,k-1])} \,.
\end{equation}

By \cite[Theorem 5.3]{GYbig}, one of the following cases holds:
\begin{enumerate}
\item[(i)] $\sig(k) > \sig(1)$, $y_{\sig,k} = \la y_{[p^m(\sig(k)), \sig(k)]}$, $m = \max \{ n \in \Znn \mid p^n(\sig(k)) \in \sig([1,k]) \}$,
\item[(ii)] $\sig(k) < \sig(1)$, $y_{\sig,k} = \la y_{[\sig(k), s^m(\sig(k))]}$, $m = \max \{ n \in \Znn \mid s^n(\sig(k)) \in \sig([1,k]) \}$,
\end{enumerate}
for some $\la \in \kx$.

Case (i). By the definition \eqref{tau} of $\Xi_N$, 
$$
\sig(k) = \max \sig([1,k]) \quad \text{and} \quad \sig([1,k-1]) \subseteq [1,\sig(k)-1].
$$
As $p(\sig(k)) = \sig(l) \in \sig([1,k])$, we also have $m \ge 1$, and so
$$
y_{[p^m(\sig(k)), \sig(k)]} = y_{[p^m(\sig(k)), p(\sig(k))]} x_{\sig(k)} - c', \quad 0 \ne c' \in R_{[p^m(\sig(k)), \sig(k)-1]} \,.
$$
Comparing terms in $R_{[1,\sig(k)]} = R_{[1,\sig(k)-1]} x_{\sig(k)} + R_{[1,\sig(k)-1]}$, we find that
$$
y_{\sig,l} = \la y_{[p^m(\sig(k)), p(\sig(k))]} \,.
$$
Since $\lc( y_{[p^m(\sig(k)), p(\sig(k))]} ) = 1$ and $y_{\sig,l} \in R'$, we find that $\la \in \DD$. In view of part (a), we conclude that
$$
 y_{\sig,k} = \la y_{[p^m(\sig(k)), \sig(k)]} \in R' \,.
 $$
 
Case (ii). Now $\sig(k) = \min \sig([1,k])$ and $\sig([1,k-1]) \subseteq [\sig(k)+1, N]$. If $m=0$, we would have $y_{\sig,k} = \la y_{\sig(k),\sig(k)]} = \la x_{\sig(k)}$, contradicting \eqref{ysigk}. Thus, $m>0$. By \cite[Theorem 5.1(d)]{GYbig},
$$
y_{[\sig(k), s^m(\sig(k))]} = x_{\sig(k)} y_{[s(\sig(k)), s^m(\sig(k))]} - c', \quad 0 \ne c' \in R_{[\sig(k)+1, s^m(\sig(k))]} \,.
$$
We may rewrite $y_{\sig,k}$ in the form
$$
y_{\sig,k} = \mu^{-1} x_{\sig(k)} y_{\sig,l} - \wt{c},
$$
where $\mu \in \kx$ arises from $\theta''_{\sig(k)}(y_{\sig,l}) = \mu y_{\sig,l}$ and
$$
\wt{c} = \mu^{-1} \de''_{\sig(k)}(y_{\sig,l}) + c \in R_{\sig([1,k-1])} \subseteq R_{[\sig(k)+1,N]} \,.
$$
Now $y_{\sig,l}$ is a homogeneous element of $R'$ and $R'$ is a $\DD$-form for \eqref{tauOre}. Moreover, $y_{\sig,l}$ has leading coefficient $1$ with respect to the presentation \eqref{tauOre}, so $\theta''_{\sig(k)}(y_{\sig,l})$ must be a $\DD^*$-multiple of $y_{\sig,l}$. Hence, $\mu \in \DD^*$.

Comparing terms in $R_{[\sig(k),N]} = x_{\sig(k)} R_{[\sig(k)+1,N]} + R_{[\sig(k)+1,N]}$,
we find that
$$
\mu^{-1} y_{\sig,l} = \la y_{[s(\sig(k)), s^m(\sig(k))]} \,.
$$
Since $\lc( y_{[s(\sig(k)), s^m(\sig(k))]} ) = 1$ while $y_{\sig,l} \in R'$ and $\mu \in \DD^*$, we obtain $\la \in \DD$, and therefore $y_{\sig,k} = \la y_{[\sig(k), s^m(\sig(k))]} \in R'$ in view of part (a). This concludes the second case of the inductive step.
\end{proof}

\subsection{Quantum cluster algebra structures on integral forms}
For integral forms of appropriately normalized symmetric CGL extensions, we have the following exact analog of \thref{maint}. Fix a symmetric CGL extension $R$ of length $N$ such that Conditions (A) and (B) hold. Set $\FF := \Fract(R)$, and let $\DD$ be a commutative domain whose field of fractions is $\KK$. Define toric frames $M_\sig : \Zset^N \rightarrow \FF$, multiplicatively skew-symmetric matrices $\rb_\sig \in M_N(\KK)$, and sets $\ex_\sig \subseteq [1,N]$ as in Subsection \ref{thm8.2GYbig}. (Recall the notation $M = M_\id$, $\rb = \rb_\id$, $\ex = \ex_\id$.) Provided the matrices $\rb_\sig$ have entries from $\DD^*$, the frames $M_\sig$ also qualify as toric frames over $\DD$, and we shall view them as such.

\bth{qcl.intform}
Let $R$ be a symmetric CGL extension of length $N$ as in Definitions {\rm\ref{dCGL}, \ref{dsymCGL}}, and assume that Conditions {\rm (A), (B)}, and \eqref{cond} hold. Let $\DD$ be a {\rm(}commutative{\rm)} domain with quotient field $\Fract (\DD) = \KK$, such that the scalars $\nu_{kl}$ in Condition {\rm(A)} all lie in $\DD^*$. Assume that the CGL presentation \eqref{itOre} has a $\DD$-form $R_\DD = \DD \langle x_1,\dots, x_N \rangle$ which contains the homogeneous prime elements $y_1,\dots,y_N$ from Theorem {\rm\ref{tmaint}}.

{\rm(a)} For each $\sig \in \Xi_N$, let $\wt{B}_\sig$ be the $N\times |\ex|$ integer matrix determined as in Theorem {\rm\ref{tmaint}(a)}. Then the pair $(M_\sig, \wt{B}_\sig)$ is a quantum seed for $\FF := \Fract(R) = \Fract(R_\DD)$ over $\DD$, and the principal part of $\wt{B}_\sig$ 
is skew-symmetrizable via the integers $d_{\eta(k)}$, $k \in \ex_\sig$ from Condition {\rm (B)}.

{\rm(b)} All the quantum seeds $(M_\sig, \wt{B}_\sig)$ from part {\rm(a)} are mutation-equivalent to each other up to the $S_N$ action. They are linked by  sequences of one-step mutations of the following kind.
Suppose $\sig, \sig' \in \Xi_N$ are such that
$$
\sig' = ( \sig(k), \sig(k+1)) \circ \sig = \sig \circ (k, k+1)
$$
for some $k \in [1,N-1]$.
If $\eta(\sig(k)) \neq \eta (\sig(k+1))$, then $M_{\sig'} = M_\sig \cdot (k,k+1)$
in terms of the action \eqref{r-act}.
If $\eta(\sig(k)) = \eta (\sig(k+1))$, 
then $M_{\sig'} = \mu_k (M_\sig)  \cdot (k,k+1)$.

{\rm(c)} The algebra $R_\DD$ equals the quantum cluster and upper cluster algebras over $\DD$ associated to $M$, $\wt{B}$, $\varnothing$:
$$
R_\DD = \AA(M, \wt{B}, \varnothing)_\DD = \UU(M, \wt{B}, \varnothing)_\DD \,.
$$
In particular, $\AA(M, \wt{B}, \varnothing)_\DD$ is a finitely generated $\DD$-algebra, and it is noetherian if $\DD$ is noetherian. In fact, $\AA(M, \wt{B}, \varnothing)_\DD$ is generated by the cluster variables in the seeds parametrized by the finite subset 
$\Ga_N$ of $\Xi_N$, recall \eqref{GaN}.

{\rm(d)} For any subset $\inv$ of the set $P(N)$ of frozen variables, there are equalities 
$$
R_\DD[y_k^{-1} \mid k \in \inv] = \AA(M, \wt{B}, \inv)_\DD = \UU(M, \wt{B}, \inv)_\DD \,.
$$
\eth 

\begin{proof}

(a) We already have from \thref{maint}(a) that $(M_\sig, \wt{B}_\sig)$ is a quantum seed for $\FF$ over $\KK$ and that the principal part of $\wt{B}_\sig$ is skew-symmetrizable via the $d_{\eta(k)}$, $k \in \ex_\sig$. The entries of $\rb(M_\sig) = \rb_\sig$, given in \eqref{tau-frame1}, lie in $\DD^*$ due to the assumption that all $\nu_{kl} \in \DD^*$. Since $\KK = \Fract(\DD)$, we have $\Fract \DD \langle M_\sig(\Zset^N) \rangle = \Fract \KK \langle M_\sig(\Zset^N) \rangle = \FF$, and so $(M_\sig, \wt{B}_\sig)$ is also a quantum seed for $\FF$ over $\DD$.

(b) This is immediate from \thref{maint}(c).

(c) and (d) are proved below. 
\end{proof}
\subsection{Examples}
\label{exam}
\bex{nongr&nonuni} Consider a uniparameter quantized Weyl algebra $R = A_n^{q,\boldsymbol\al}(\KK)$,
for a non-root of unity $q \in \kx$ and a skewsymmetric matrix $\boldsymbol\al = (a_{ij}) \in M_n(\Zset)$. This algebra is presented by generators $v_1,w_1,\dots, v_n,w_n$ and relations
\begin{equation} \label{Anqpsi.rels}
\begin{aligned}
w_iw_j &=q^{a_{ij}} w_jw_i \,, &\qquad &(\text{all\ } i,j), \\
v_iv_j &=q^{1+a_{ij}} v_jv_i \,, &\qquad &(i<j), \\ 
v_iw_j &= q^{-a_{ij}} w_j v_i \,, &\qquad &(i<j), \\ 
v_iw_j &=q^{1-a_{ij}} w_j v_i \,, &\qquad &(i>j), \\ 
v_jw_j &=1+qw_j v_j
+ (q-1) \sum_{l<j} w_l v_l \,, &\qquad &(\text{all\ }j). 
\end{aligned}
\end{equation}
The torus $\HH = (\kx)^n$ acts rationally on $R$ with
$$
(\al_1,\dots,\al_n).v_i = \al_i v_i \qquad\text{and}\qquad (\al_1,\dots,\al_n).w_i = \al_i^{-1} w_i
$$
for $(\al_1,\dots,\al_n) \in \HH$ and $i \in [1,n]$. With the variables $v_1,w_1,\dots, v_n,w_n$ in the listed order,  $R$ is a CGL extension, but that presentation is not symmetric. There is a symmetric CGL extension presentation with the variables in the order $w_n,\dots,w_1,v_1,\dots,v_n$, and $\DD\langle w_n,\dots,w_1,v_1,\dots,v_n \rangle$ is a $\DD$-form for this presentation, where $\DD = (\Zset{\cdot} 1_\KK) [q^{\pm1}]$. However, the homogeneous prime elements $y_1,\dots,y_{2n}$ from \thref{1} do not lie in this $\DD$-form; see \reref{Dformcomments}(c). This can be rectified by rescaling the generators as in \prref{rescale.Dform}. One such rescaling leads to the CGL presentation
\begin{equation}  \label{Anq.pres2}
R = \KK[(q-1) w_n] \cdots [(q-1) w_1; \theta_n] [v_1; \theta_{n+1}, \de_{n+1}] \cdots [v_n; \theta_{2n}, \de_{2n}],
\end{equation}
and $\DD\langle (q-1) w_n,\dots, (q-1) w_1,v_1,\dots,v_n \rangle$ is a $\DD$-form for this presentation which contains all the $y_k$.

If either $R$ or a $\DD$-form of $R$ is $\Znn$-graded, then in view of the final relations in \eqref{Anqpsi.rels} all the generators $v_i$, $w_j$ must be homogeneous of degree $0$. Thus, $R$ and its $\DD$-forms have no nontrivial $\Znn$-gradings. 
\qed
\eex

\bex{uniqWeyl}
Let $R = A_n^{q,\boldsymbol\al}(\KK)$ and $\HH = (\kx)^n$ as in \exref{nongr&nonuni}, and take the symmetric CGL presentation \eqref{Anq.pres2}. Set $\DD = (\Zset{\cdot} 1_\KK) [q^{\pm1}]$. Then 
$$
\DD\langle (q-1) w_n,\dots, (q-1) w_1,v_1,\dots,v_n \rangle
$$
is a $\DD$-form for the presentation \eqref{Anq.pres2} which contains the homogeneous prime elements $y_1,\dots,y_{2n}$ from \thref{1}.

The CGL presentation \eqref{Anq.pres2} satisfies Condition (B) with all $d_i=1$, and to obtain Condition (A) we just need to assume that $\KK$ contains a square root of $q$. Choose one, and label it $q^{1/2}$. The condition \eqref{cond}, however, only holds after a further rescaling of the generators. Namely, write $R$ as an iterated Ore extension with variables $x_1,\dots, x_{2n}$ where
$$
x_i := \begin{cases}
(q-1) w_{n+1-i}  &\quad (\text{if}\ i \in [1,n]),  \\
(-1)^{i-n} q^{(i-n-1)/2} v_{i-n}  &\quad (\text{if}\ i \in [n+1,2n]).
\end{cases}
$$
In order to express the relations among these $x_i$ in a convenient form, we use the following notation:
$$
l' := 2n+1-l \ \ (\text{for}\ l \in [1,2n]) \qquad\text{and}\qquad c_{ij} := a_{n+1-i, n+1-j} \ \ (\text{for}\ i,j \in [1,n]).
$$
Then $R$ has the presentation with generators $x_1,\dots, x_{2n}$ and defining relations
\begin{equation}  \label{Anq.pres3}
\begin{gathered}
\begin{aligned}
x_i x_j &= q^{c_{ij}} x_j x_i \,,  &\qquad (i,j \in [1,n])  \\
x_i x_j &= q^{1+c_{i'j'}} x_j x_i \,,  &\qquad (n< i< j\le 2n)  \\
x_i x_j &= q^{-c_{i'j}} x_j x_i \,,  &\qquad (j\le n< i < j'\le 2n)  \\
x_i x_j &= q^{1 -c_{i'j}} x_j x_i \,,  &\qquad (j \le n< j'< i\le 2n)  
\end{aligned}  \\
\begin{aligned}
x_{j'} x_j &= (-1)^{n+1-j} q^{(n-j)/2} (q-1) + q x_j x_{j'}  \\
&\qquad + (q-1) \sum_{1\le j<l \le n} (-1)^{l-j} q^{(l-j)/2} x_l x_{l'} \,,  \qquad (j \in [1,n]).  
\end{aligned}
\end{gathered}
\end{equation}

With the presentation \eqref{Anq.pres3}, $R$ is a symmetric CGL extension satisfying the required hypotheses (A), (B), and \eqref{cond} of \thref{qcl.intform}. 
It has a $\DD$-form 
\begin{equation}
\label{qWeyl}
A_n^{q,\boldsymbol\al}(\DD)  := \DD \langle x_1,\dots, x_{2n} \rangle
\end{equation}
where we now take $\DD = (\Zset{\cdot} 1_\KK) [q^{\pm1/2}]$. There are two possibilities for $\DD$:
\begin{align*}
&\DD \cong \iAA = \Zset[ q^{1/2} ], \quad &&\mbox{if} \quad \charr \KK = 0, 
\\
&\DD \cong {\mathbb{F}}_p [ q^{1/2} ], \quad &&\mbox{if} \quad \charr \KK = p
\end{align*} 
(recall \eqref{iAA}).
The scalars $\nu_{kl}$ from Condition (A) all lie in $\DD^*$, as do \red{the nonzero coefficients 
of} the homogeneous prime elements $y_1, \dots, y_{2n}$ from \thref{maint}. Therefore
$$
A_n^{q,\boldsymbol\al}(\DD) =  \AA(M, \wt{B}, \varnothing)_\DD =  \UU(M, \wt{B}, \varnothing)_\DD
$$
by \thref{qcl.intform}. The matrix $\rb = \rb(M)$ of the initial toric frame has the form
$$
\rb = \begin{bmatrix}
1 &s^{c_{12}} &\cdots &s^{c_{1,n-1}} &s^{c_{1n}}  &1 &1 &\cdots &1 &s^{-1}  \\
s^{c_{21}} &1 &\cdots &s^{c_{2,n-1}} &s^{c_{2n}}  &1 &1 &\cdots &s^{-1} &s^{-1}  \\
\vdots &\vdots &&\vdots &\vdots  &\vdots &\vdots &&\vdots &\vdots  \\
s^{c_{n1}} &s^{c_{n2}} &\cdots &s^{c_{n,n-1}} &1  &s^{-1} &s^{-1} &\cdots &s^{-1} &s^{-1}  \\
1 &1 &\cdots &1 &s  &1 &1 &\cdots &1 &1  \\
1 &1 &\cdots &s &s  &1 &1 &\cdots &1 &1  \\
\vdots &\vdots &&\vdots &\vdots  &\vdots &\vdots &&\vdots &\vdots  \\
s &s &\cdots &s &s  &1 &1 &\cdots &1 &1 
\end{bmatrix} ,
$$
where $s := q^{1/2}$.
The \red{quiver} of the initial seed is acyclic, namely it equals
\[
\begin{tikzpicture}
   \node (1') at (0,0) {$2n$};
    \node (2') at (2,0) {$2n-1$};
    \node (3') at (3.2,0) {$\cdots$}; 
    \node (4') at (4.4,0) {$n+2$};
    \node (5') at (6.4,0) {$n+1$};
    \node (1) at (0,1.2) {$1$};
    \node (2) at (2,1.2) {$2$};
    \node (3) at (3.2,1.2) {$\cdots$}; 
    \node (4) at (4.4,1.2) {$n-1$};
    \node (5) at (6.4,1.2) {$n$};
    \draw
    (1') edge[->,font=\scriptsize,>=angle 90] (1)
     (2') edge[->,font=\scriptsize,>=angle 90] (2)
     (4') edge[->,font=\scriptsize,>=angle 90] (4)
     (5') edge[->,font=\scriptsize,>=angle 90] (5)
     (1) edge[->,font=\scriptsize,>=angle 90] (2')
     (4) edge[->,font=\scriptsize,>=angle 90] (5');
\end{tikzpicture}
\]
where the top vertices are mutable and the bottom ones are frozen.
\qed
\eex

Next, we illustrate \thref{qcl.intform} with a CGL extension which is not $\Zset_{\geq 0}$-graded connected and whose quiver is not acyclic.

\bex{last-ex} Recall the notation $\iFF := \Qset(q^{1/2})$. Let $R$ be the $\iFF$-algebra with generators $x_1, \ldots, x_6$ and relations
$$
\begin{gathered}
\begin{aligned}
x_2 x_1 &= q x_1x_2 \,,  &\quad x_3x_1 &= q x_1x_3 +(1-q)  x_2^2 \,,  &\quad x_3x_2 &= q x_2x_3 \,,  \\
x_4 x_1 &= q x_1 x_4 + (1-q^2) x_2 x_3 \,,  &x_4 x_2 &= q x_2 x_4 + (q-1) x_3^2 \,,  &x_4 x_3 &= q x_3 x_4 \,,  \\
x_5 x_1 &= q^{-1} x_1 x_5 \,,  &x_5 x_2 &= q^{-1} x_2 x_5 \,,  &\quad x_5 x_3 &= q^{-1} x_3 x_5 \,,  \\
&&x_5 x_4 &= q^{-1} x_4 x_5 + (q-1)^3 \,,
\end{aligned}  \\
\begin{aligned}
x_6 x_1 &= q^{-1} x_1 x_6 + (q^{-1}-q) x_2 y x_5 + (1-q) x_3^2 x_5^2 \,,  \\
x_6 x_2 &= q^{-1} x_2 x_6 + (q-q^{-1}) x_3 y x_5 \,,
\end{aligned}  \\
\begin{aligned}
x_6 x_3 &= q^{-1} x_3 x_6 + (q-1) y^2 \,,  &\qquad\quad x_6 x_4 &= q^{-1} x_4 x_6 \,,  &\qquad\quad x_6 x_5 &= q x_5 x_6 \,, 
\end{aligned}
\end{gathered}$$
where 
$$
y := x_4 x_5 - q(1-q)^2 \,.
$$
The algebra $R$ is a symmetric CGL extension for the torus $\HH := ((\iFF)^{\times})^2$ acting so that for the corresponding grading by $\xh \cong \Zset^2$, the variables $x_1,\dots,x_6$ have degrees
$$
(4,3),\ (3,2),\ (2,1),\ (1,0),\ (-1,0),\ (-2,-1).
$$
The $h$-elements for this CGL extension are
$$
h_3 = h_4 = (q,q^{-1}), \; \;  h_5 =  h_6 = (q^{-1},q) \in \HH.
$$
Consequently, $\la_k =  q$ for $k \in [3,6]$. The (nonunique) elements $h_1, h_2 \in \HH$ can be also chosen so that 
$\la_k =  q$ for $k =1,2$. Obviously Conditions (A) and (B) hold.

Denote by $R_{\iAA}$ the $\iAA$-subalgebra of $R$ generated by $x_1, \ldots, x_6$.
The homogeneous prime elements $y_1,\dots, y_6$ belong to $R_{\iAA}$ and are given by 
$$
\begin{gathered}
y_1 = x_1 \,,  \qquad \quad y_2 = x_2 \,,  \qquad \quad y_3 = x_1 x_3 + q^{-1} x_2^2 \,,  \\
y_4 = x_2 x_4 - q^{-1} x_3^2 \,,  \qquad \quad y_5 = x_2 x_4 x_5 - q^{-1} x_3^2 x_5 - q (1-q)^2 x_2 \,,  \\
y_6 = x_1 x_3 x_6 + q^{-1} x_2^2 x_6 - q x_1 y^2 - (1+q^{-1}) x_2 x_3 y x_5 + q^{-2} x_3^3 x_5^2 \,.
\end{gathered}
$$
(The element $y$ is precisely the interval prime element $y_{[3,5]}$.) Consequently,  the $\eta$-function
from \thref{1}
is given by $\eta(1)= \eta(3) = \eta(6) = 1$ and $\eta(2) = \eta(4) = \eta(5) = 2$. 
Hence, the predecessor function $p$ maps $6 \mapsto 3 \mapsto 1$ and $5 \mapsto 4 \mapsto 2$. So, $\ex = [1,4]$.
One easily verifies that the condition \eqref{cond} is satisfied. 
The matrix of the initial toric frame for $R_{\iAA}$ from \thref{qcl.intform} is given by
$$
\rb = \begin{bmatrix}
1 &s^{-1} &s^{-1} &s^{-2} &s^{-1} &1  \\
s &1 &1 &s^{-1} &1 &s  \\
s &1 &1 &s^{-2} &1 &s^2  \\
s^2 &s &s^2 &1 &s^2 &s^4  \\
s &1 &1 &s^{-2} &1 &s  \\
1 &s^{-1} &s^{-2} &s^{-4} &s^{-1} &1
\end{bmatrix} ,
$$
where $s:= q^{1/2}$. The \red{quiver} of the initial quantum seed of $R_{\iAA}$ is
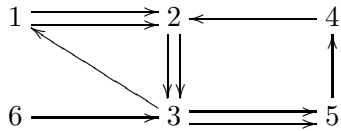
$$
{\xymatrixrowsep{2pc}  \xymatrixcolsep{4pc}
\xymatrix{
1 \ar[r]<0.5ex> \ar[r]<-0.5ex>  &2 \ar[d]<0.5ex> \ar[d]<-0.5ex>  &4 \ar[l]  \\
6 \ar[r]  &3 \ar[ul] \ar[r]<0.5ex> \ar[r]<-0.5ex>  &5 \ar[u]
}}
$$
where the vertices $5,6$ are frozen and the rest are mutable. \thref{qcl.intform} implies that 
$R_{\iAA}$ is isomorphic to the corresponding cluster and upper cluster algebras 
over $\iAA$ where the two frozen variables are not inverted. 

All statements in the example hold if $\iAA$ and $\iFF$ are replaced by ${\mathbb{F}}_p[q^{\pm 1/2}]$ and ${\mathbb{F}}_p(q^{1/2})$, 
respectively.
\qed
\eex

\red{
\bre{not-qunip} The algebras in Examples \ref{enongr&nonuni}--\ref{elast-ex} do not come from quantum unipotent cells in any symmetrizable Kac--Moody algebra, 
because the algebras in those examples are $\Zset$-graded but they are not $\Zset_{\geq 0}$-graded connected algebras while all quantum unipotent cells
are $\Zset_{\geq 0}$-graded connected algebras. In particular, these examples concern applications of \thref{qcl.intform} that are not covered by 
\cite{KKKO} or the results in Sect. \ref{qunip} of this paper.
\ere
}
\bre{mult-par} There are also simple examples of symmetric CGL extensions $R$ which cannot be ``untwisted" into a uniparameter form. More precisely, there are such $R$ for which no twist of $R$ relative to a $\kx$-valued cocycle on a natural grading group turns $R$ into a uniparameter CGL extension. For instance, this is true of the multiparameter quantized Weyl algebra $A_n^{Q,P}(\KK)$ when the parameters in the vector $Q = (q_1,\dots,q_n)$ generate a non-cyclic subgroup of $\kx$ (see \cite[Example 5.10]{GYtwist}). 
One can show that the quantized Weyl algebras $A_n^{Q,P}(\KK)$ have integral forms over subrings $\Zset[q_1^{ \pm 1/2}, \dots, q_n^{\pm 1/2}]$ of $\KK$.
\thref{qcl.intform} can be applied to prove that the integral forms are isomorphic to quantum cluster algebras over $\Zset[q_1^{ \pm 1/2}, \dots, q_n^{\pm 1/2}]$.
\ere

\subsection{Proof of parts (c), (d) of \thref{qcl.intform}} For the first part of this subsection, we assume only that $\KK = \Fract(\DD)$. The normalization assumptions in \thref{qcl.intform} will be invoked only in the proof of parts (c), (d) of the theorem.

In the following lemma and proposition, divisibility refers to divisibility within the ring $R_\DD$.

\ble{d.uv}
Assume that \eqref{itOre} has a $\DD$-form $R_\DD = \DD \langle x_1,\dots, x_N \rangle$.
Let $d \in \DD {\setminus}\{0\}$ and $u,v \in R_\DD {\setminus}\{0\}$ such that $d \mid uv$. If $\lc(v) \in \DD^*$, then $d \mid u$.
\ele

\begin{proof}
Let $\lt(u) = b x^f$ and $\lt(v) = c x^g$ where $b,c \in \DD {\setminus} \{0\}$ and $f,g \in \Znn^N$. By assumption, $c \in \DD^*$ and $uv = dw$ for some $w \in R' {\setminus}\{0\}$. We proceed by induction on $f$ with respect to $\prec$. If $f=0$, we have $u = b$ and $bc x^g = \lt(uv) = d\, \lt(w)$. In this case, $d \mid bc$, whence $d$ divides $b=u$ because $c$ is a unit in $\DD$.

Now assume that $f \succ 0$. In view of \cite[Eq.~(3.20)]{GYbig}, we have
$$
\la bc = \lc(uv) = d\, \lc(w)
$$
for some $\la \in \DD$ which is a product of $\la_{k,j}\,$s. By assumption, $\la$ is a unit in $\DD$, whence $b = d e$ for some $e \in \DD$. Now $u = de x^f + u'$ where either $u' = 0$ or $\lt(u') = b' x^{f'}$ with $b' \in \DD$ and $f' \prec f$. In the second case,
$$
u'v = uv - de x^f v = d(w - e x^f v).
$$
By induction, $d \mid u'$, and thus $d \mid u$. This verifies the induction step.
\end{proof}

\bpr{intersecprop}
Assume that \eqref{itOre} has a $\DD$-form $R_\DD = \DD \langle x_1,\dots, x_N \rangle$ which contains $y_1,\dots,y_N$.
If $Y$ is the multiplicative set generated by $\DD^* \cup \{y_1,\dots,y_N\}$, then
\begin{equation}  \label{intersect}
R_\DD [Y^{-1}] \cap R = R_\DD \,.
\end{equation}
\epr

Recall from \leref{Dforms3}(b) that $Y$ is a denominator set in $R_\DD$.

\begin{proof} 
If $r \in R_\DD [Y^{-1}] \cap R$, then $r = ay^{-1}$ for some $a \in R_\DD$ and $y \in Y$. Since $r \in \KK R_\DD$, we also have $r = d^{-1}b$ for some $d \in \DD {\setminus}\{0\}$ and $b \in R_\DD$. Now $da = by$. Since $\lc(y_j) = 1$ for all $j \in [1,N]$, we see via \cite[Eq.~(3.20)]{GYbig} that $\lc(y) \in \DD^*$. By \leref{d.uv}, $b = db'$ for some $b' \in R_\DD$. Thus $a = b'y$ and therefore $r = ay^{-1} = b' \in R_\DD$.
\end{proof}

From now on, assume that $R$ is a symmetric CGL extension and  that \eqref{itOre} has a $\DD$-form $R_\DD = \DD \langle x_1,\dots, x_N \rangle$ which contains $y_1,\dots,y_N$.. For each $\sig \in \Xi_N$, we have the CGL  presentation \eqref{tauOre} for $R$, and $R_\DD$ is a $\DD$-form of this presentation by \leref{Dforms2}(b). Let $y_{\sig,1}, \dots, y_{\sig,N}$ be the (unnormalized) sequence of homogeneous prime elements from \thref{maint} for the presentation \eqref{tauOre}, and let $E_\sig$ denote the multiplicative set generated by 
$$
\DD^* \cup \{ y_{\sig,l} \mid l \in [1,N],\ s_\sig(l) \ne +\infty \} = \DD^* \cup \{ y_{\sig,l} \mid l \in \ex_\sig \},
$$
where $s_\sig$ is the successor function for the level sets of $\eta\sig$. (By \cite[Corollary 5.6(b)]{GYbig}, $\eta\sig$ can be chosen as the $\eta$-function for the presentation \eqref{tauOre}.)
By \prref{tauDform}(c) and \leref{Dforms3}(b), $E_\sig$ is a denominator set in $R_\DD$.

\bpr{intersecprop2}
The ring $R_\DD$ equals the following intersection of localizations:
\begin{equation}  \label{intersect2}
R_\DD = \bigcap_{\sig \in \Ga_N} R_\DD [E_\sig^{-1}].
\end{equation}
\epr

\begin{proof} Let $T$ denote the right hand side of \eqref{intersect2}. Since $\bigcap_{\, \sig \in \Ga_N} R[E_\sig^{-1}] = R$ by \cite[Theorem 8.19(d)]{GYbig}, we have $T \subseteq R$. On the other hand, $E_\id$ is contained in the denominator set $Y$ of \prref{intersecprop}, and so $T \subseteq R_\DD[Y^{-1}]$. \prref{intersecprop} thus implies $T \subseteq R_\DD$, yielding \eqref{intersect2}.
\end{proof}

\bco{interseccor}
If $\inv$ is any subset of $[1,N]{\setminus}\ex$, then
\begin{equation}  \label{intersect3}
R_\DD [y_k^{-1} \mid k \in \inv] = \bigcap_{\sig\in\Ga_N} R_\DD[E_\sig^{-1}] [y_k^{-1} \mid k \in \inv].
\end{equation}
\eco

\begin{proof} This follows from \prref{intersecprop2} in the same way that \cite[Theorem 8.19(e)]{GYbig} follows from \cite[Theorem 8.19(d)]{GYbig}.
\end{proof} 

\begin{proof}[Proof of Theorem {\rm\ref{tqcl.intform}(c)(d)}]
Note that the scalars $\Scr_\nub(f)$ from \eqref{Scrnub}, for $f\in \Zset^N$, lie in $\DD^*$ because of our assumption that all $\nu_{kl} \in \DD^*$. Hence, invoking \prref{tauDform}(a), the normalized elements $\ol{y}_j$ and $\ol{y}_{[i,s^m(i)]}$ from \eqref{eybarj} and \eqref{ybarismi} belong to $R_\DD$. By \eqref{tau-frame2}, we thus have $\ol{y}_{\sig,k} \in R_\DD$ for all $\sig \in \Xi_N$ and $k \in [1,N]$.

We next show that
\begin{equation}  \label{RDgenybars}
R_\DD = \DD \langle \ol{y}_{\sig,k} \mid \sig \in \Ga_N\,, \ k \in [1,N] \rangle.
\end{equation}
The proof is parallel to that for the corresponding statement in \cite[Theorem 8.2(b)]{GYbig}. For each $j \in [1,N]$, there is an element $\sig \in \Ga_N$ with $\sig(1) = j$. By \eqref{tau-frame2}, $\ol{y}_{\sig,1}$ is a $\DD^*$-multiple of $y_{[j,j]} = x_j$, and so $x_j \in \DD^* \ol{y}_{\sig,1}$. Therefore all $x_j$ lie in the right hand side of \eqref{RDgenybars}, and the equation is established. Since all the $\ol{y}_{\sig,k} = M_\sig(e_k)$ are cluster variables, it follows that $R_\DD \subseteq \AA(M,\wt{B},\varnothing)_\DD$.

We have $\AA(M,\wt{B},\varnothing)_\DD \subseteq \UU(M,\wt{B},\varnothing)_\DD$ by the Laurent Phenomenon \eqref{Lphenom}, and 
$$
\UU(M,\wt{B},\varnothing)_\DD \subseteq \bigcap_{\sig \in \Xi_N} \DD \TT_{(M_\sig,\wt{B}_\sig,\varnothing)} = \bigcap_{\sig \in \Xi_N} \DD \langle \ol{y}_{\sig,k}^{\pm1},\, \ol{y}_{\sig,j} \mid k \in \ex_\sig, \, j \in [1,N]{\setminus}\ex_\sig \rangle,
$$
where $\ex_\sig$ appears instead of $\ex$ for the indexing reasons explained ahead of \thref{maint}.
Since $\DD \langle \ol{y}_{\sig,k}^{\pm1},\, \ol{y}_{\sig,j} \mid k \in \ex_\sig, \, j \in [1,N]{\setminus}\ex_\sig \rangle \subseteq R_\DD[E_\sig^{-1}]$ for each $\sig \in \Xi_N$, we obtain 
$$
\UU(M,\wt{B},\varnothing)_\DD \subseteq \bigcap_{\sig \in \Xi_N} R_\DD [E_\sig^{-1}].
$$
In view of \prref{intersecprop2}, we have the following sequence of inclusions:
\begin{equation}  \label{inclseq}
R_\DD \subseteq \AA(M,\wt{B},\varnothing)_\DD \subseteq \UU(M,\wt{B},\varnothing)_\DD \subseteq \bigcap_{\sig \in \Ga_N} R_\DD [E_\sig^{-1}] = R_\DD \,.
\end{equation}
All the inclusions in \eqref{inclseq} must be equalities, which establishes the first part of \thref{qcl.intform}(c). The finite generation statements concerning $\AA(M,\wt{B},\varnothing)_\DD$ now follow from \eqref{RDgenybars}. If $\DD$ is noetherian, the iterated Ore extension $R_\DD$ is noetherian by standard skew polynomial ring results. This concludes the proof of part (c).

Part (d) is proved analogously, using \coref{interseccor} in place of \prref{intersecprop2}.
\end{proof}

\sectionnew{Quantum Schubert cell algebras, canonical bases and quantum function algebras}
\label{qSchubert}
\subsection{Quantized universal enveloping algebras}
\label{qaff}
Fix a (finite) index set $I =[1,r]$ and consider a {\em{Cartan datum}} $(A,P,\Pi, P\spcheck, \Pi\spcheck)$ consisting of the following:
\begin{enumerate}
\item[(i)] A {\em{generalized Cartan matrix}} $A = (a_{ij})_{i,j \in I}$ such that $a_{ii}=2$ for $i \in I$, $-a_{ij} \in \Zset_{\geq 0}$ for $i \neq j \in I$,
and there exists a diagonal matrix $D = \diag(d_i)_{i \in I}$ with relatively prime entries $d_i \in \Zset_{> 0}$ for which $DA$ is symmetric.
\item[(ii)] A free abelian group $P$ ({\em{weight lattice}}).
\item[(iii)] A subset $\Pi = \{ \al_i \mid i \in I \} \subset P$ ({\em{set of simple roots}}). 
\item[(iv)] The dual group $P\spcheck := \Hom_\Zset(P, \Zset)$ ({\em{coweight lattice}}).
\item[(v)] Two linearly independent subsets $\Pi\spcheck = \{ h_i \mid i \in I \} \subset P\spcheck$ ({\em{set of simple coroots}}) 
such that $\lcor h_i, \al_j \rcor = a_{ij}$ for $i,j \in I$, and $\{ \vpi_i \in P \mid i \in I \}$ ({\em{set of fundamental weights}}) 
such that $\lcor h_i, \vpi_j \rcor = \delta_{ij}$.  
\end{enumerate}

Let $\g$ be the {\em{symmetrizable Kac--Moody algebra}} over $\Qset$ corresponding to this Cartan datum. 
Denote 
\[
Q := \oplus_{i \in I} \Zset \al_i \subset P, \; \; Q_+ := \oplus_{i \in I} \Zset_{\geq 0} \al_i
\]
and 
\[
P_+ := \{ \ga \in P \mid \lcor h_i, \ga \rcor \in \Zset_{\geq 0}, \forall i \in I\}, \; \;  
P_{++} := \{ \ga \in P \mid \lcor h_i, \ga \rcor \in \Zset_{>0} , \forall i \in I \}.
\]
Set $\h := \Qset \otimes_{\Zset} P\spcheck$.
There exists a $\Qset$-valued nondegenerate symmetric bilinear form $(.,.)$ on $\h^*= \Qset \otimes_{\Zset} P$ such that  
\begin{equation}
\label{di}
\lcor h_i, \mu \rcor = \frac{ 2 (\al_i, \mu)}{(\al_i, \al_i)}  
\quad \mbox{and} \quad (\al_i, \al_i) = 2 d_i \quad 
\mbox{for} \; \; i \in I,\ \mu \in \h^*.
\end{equation}
Set $\| \ga \|^2 := ( \ga, \ga)$ for $\ga \in \h^*$.
Denote by $W$ the Weyl group of $\g$ acting by isometries on $(\h^*, (.,.))$. Denote by $s_i$ its 
generators, by $\ell : W \to \Zset_{\geq 0}$ the length function on $W$, and by $\geq$ the Bruhat order on $W$.
We will also denote by  $(.,.)$ the transfer of this bilinear form to $\h$, satisfying $(h_i, h_j) = (\al_i, \al_j)/ d_i d_j$ for all $i,j \in I$.

Let $U_q(\g)$ be the quantized universal enveloping algebra of $\g$ over the rational function field $\Qset(q)$. It has generators 
$q^h, e_i, f_i$ for $i \in I$, $h \in P\spcheck$ and the following relations for $h, h' \in P\spcheck$, $i, j \in I$: 
\begin{align*}
&q^0 = 1, \quad q^h q^{h'} = q^{h + h'}, \\
&q^h e_i q^{-h} = q^{\langle h, \al_i \rangle} e_i, \quad q^h f_i q^{-h} = q^{ - \langle h, \al_i \rangle }  f_i, \\
& e_i f_j - f_j e_i = \delta_{ij} \frac{q^{d_i h_i} - q^{- d_i h_i }}{q_i - q_i^{-1}}, \\
&\sum_{k=0}^{1-a_{ij}} (-1)^k
\begin{bmatrix} 
1 - a_{ij} \\
k
\end{bmatrix}_i
e_i^{1 - a_{ij} -k} e_j e_i^k =0, \quad i \neq j,
\\
&\sum_{k=0}^{1- a_{ij}} (-1)^k 
\begin{bmatrix}
1 - a_{ij} \\
k
\end{bmatrix}_i
f_i^{1- a_{ij} - k} f_j f_i^k =0, \quad i \neq j,
\end{align*} 
where 
\[
q_i := q^{d_i}, \quad [n]_i := \frac{q_i^n - q_i^{-n}}{q_i - q_i^{-1}}, \quad
[n]_i! := [1]_i \cdots [n]_i \quad \mbox{and} \quad 
\begin{bmatrix}
n \\
k
\end{bmatrix}_i
:= \frac{[n]_i}{ [k]_i [n-k]_i}
\]
for $k \leq n$ in $\Zset_{\geq 0}$ and $i \in I$. The algebra $U_q(\g)$ is a Hopf algebra with coproduct, antipode and counit such that
\begin{gather*}
\Delta(q^h) = q^h \otimes q^h, \quad 
\Delta(e_i) = e_i \otimes 1 + q^{d_i h_i} \otimes e_i, \quad
\Delta(f_i) = f_i \otimes q^{- d_i h_i} + 1 \otimes f_i, 
\\
S(q^h) = q^{-h}, \quad S(e_i) = - q^{- d_i h_i } e_i, \quad S(f_i) = - f_i q^{d_i h_i}, 
\\
\ep(q^h) =1, \quad \ep(e_i) = \ep(f_i) =0
\end{gather*}
for $h \in P\spcheck$, $i \in I$. The Hopf algebra $U_q(\g)$ is $Q$-graded with
\begin{equation}
\label{grading}
\deg e_i =  \al_i, \quad \deg f_i = - \al_i, \quad 
\deg q^h =0.
\end{equation}
For a $Q$-graded subalgebra $R$ of  $U_q(\g)$, its graded components will be denoted by $R_\ga$, where $\ga \in Q$.
For a homogeneous $x \in U_q(\g)_\ga$, set $\wtt x := \ga$. Define the torus 
\[
\HH := (\Qset(q)^\times)^I.
\]
For $\ga = \sum n_i \al_i \in Q$, let $t \mt t^\ga$ denote the character of $\HH$ given by 
$(r_i)_{i \in I} \mt \prod_i r_i^{n_i}$. This identifies the rational character lattice of $\HH$ with $Q$. 
The torus $\HH$ acts on $U_q(\g)$ by
\begin{equation}
\label{HH-act}
t \cdot x = t^\ga x \quad \mbox{for} \quad x \in U_q(\g)_\ga,\ \ga \in Q.
\end{equation}

Let $\De_+ \subset Q_+$ be the set of positive roots of $\g$.
For $w \in W$, denote the following Lie subalgebras of the Kac--Moody algebra $\g$,
\begin{equation}
\label{n-alg}
\n_\pm := \oplus_{\al \in \De_+} \g^{\pm \al}, \quad
\n_\pm(w) := \oplus_{\al \in \De_+ \cap w^{-1}(- \De_+)} \g^{ \pm \al},
\end{equation}
where for $\al \in \De_+$, $\g^{\pm \al}$ are the corresponding root spaces in $\g$. Let $\b_\pm$ be the 
corresponding Borel subalgebras of $\g$. 
Denote by $U_q(\n_\pm)$ and $U_q(\h)$ the unital subalgebras of $U_q(\g)$ respectively generated by $\{e_i \mid i \in I \}$, 
$\{f_i \mid i \in I \}$ and $\{ q^h \mid h \in P\spcheck \}$. Denote the Hopf subalgebras $U_q(\b_\pm) := U_q(\n_\pm) U_q(\h)$
of $U_q(\g)$. 

Consider the $\Qset(q)$-linear anti-automorphisms $*$ and $\varphi$ of $U_q(\g)$ defined by
\begin{align*}
&e_i^* := e_i, \quad f_i^* := f_i, \quad (q^h)^* := q^{-h}, \quad \mbox{and} \\
&\varphi(e_i) := f_i, \quad \varphi(f_i) := e_i, \quad \varphi(q^h) := q^h
\end{align*}
for $i \in I$, $h \in P\spcheck$. 
Their composition $\varphi^* := \varphi \circ * = * \circ \varphi$ is the $\Qset(q)$-linear automorphism of $U_q(\g)$
satisfying
\[
\varphi^*(e_i) = f_i, \quad \varphi^*(f_i) = e_i \quad \mbox{and} \quad \varphi^*(q^h) = q^{-h}.
\]
Denote by $c \mt \ol{c}$ the automorphism of the field $\Qset(q)$ given by $\ol{q} = q^{-1}$. The bar involution 
$x \mt \ol{x}$ of $U_q(\g)$ is its $\Qset(q)$-skewlinear automorphism  
such that $\ol{c x} = \ol{c} \, \ol{x}$ for $c \in \Qset(q)$, $x \in U_q(\g)$ and $\ol{f}_i = f_i$, $\ol{e}_i = e_i$, $\ol{q^h} = q^{-h}$ for $i \in I$, $h \in P\spcheck$.
Denote the $\Qset(q)$-skewlinear antiautomorphism $\ol{\vp}$ of $U_q(\g)$, 
\[
\ol{\vp}(x) := \vp(\ol{x})= \ol{\vp(x)}, \quad \forall x \in U_q(\g).
\]

A $U_q(\g)$-module $V$ is called {\em{integrable}} if $e_i$ and $f_i$ act locally nilpotently on $V$ and
\[
V = \oplus_{\mu \in P} V_\mu \; \; \mbox{with} \; \;  \dim V_\mu < \infty, \quad \mbox{where} \quad V_\mu = \{ v \in M \mid q^h \cdot v = 
q^{ \lcor h, \mu \rcor } v , \ \forall h \in P\spcheck \}.
\]
The category $\Ointg$ consists of the integrable $U_q(\g)$ modules whose nontrivial graded subspaces have weights in 
$\cup_j (\mu_j + Q)$ for finitely many $\mu_1, \ldots, \mu_n \in P$ (depending on the module). It is a semisimple monoidal category 
with respect to the tensor product of $U_q(\g)$-modules and 
with simple objects given by the irreducible highest weight modules $V(\mu)$ with highest weights $\mu \in P_+$.

For $V \in \Ointg$ its restricted dual module with respect to the antiautomorphism $\vp$ 
is a module in $\Ointg$ defined by 
\[
D_\vp V := \oplus_{\mu \in P} V_\mu^*,  \quad \mbox{where $V_\mu^*$ is the dual $\Qset(q)$-vector space of $V_\mu$}.
\]
The $U_q(\g)$-action on $D_\vp V$ is given by $\lcor x \cdot \xi , v \rcor = \lcor \xi, \vp(x) \cdot v \rcor$ for $v \in V$, $\xi \in D_\vp V$. 

Denote by $\{T_i \mid i \in I \}$ the generators of the braid group of $W$. For $w \in W$, let $T_w := T_{i_1} \cdots T_{i_N}$ for a reduced expression $s_{i_1} \cdots s_{i_N}$ of $w$. 
We will denote by the same notation Lustig's braid group action \cite{Lusztig2} on $U_q(\g)$ and on the modules in $\Ointg$. 
We will follow the conventions of \cite{Jantzen1}.

\subsection{Two bilinear forms}
Consider the $\Qset(q)$-linear skew-derivations $e''_i$ of $U_q(\n_-)$,
\[
e''_i (f_j) = \delta_{ij} \quad \mbox{and} \quad
e''_i(x y) = e''_i(x) y + q_i^{- \lcor h_i , \gamma \rcor} x e''_i(y)
\]
for all $i,j \in I$, $x \in U_q(\n_-)_\ga$, $y \in U_q(\n_-)$.
The Kashiwara--Lusztig nondegenerate, symmetric bilinear form $(-, -)_{KL} : U_q(\n_-) \times U_q(\n_-) \to \Qset(q)$ is 
the unique bilinear form such that
\[
(1, 1)_{KL} =1 \quad \mbox{and} \quad (f_i x, y)_{KL} = (q_i^{-1} - q_i)^{-1} (x, e''_i(y))_{KL}, 
\quad \forall i \in I, \ x, y \in U_q(\n_-). 
\]
\bre{forms} The Lusztig form uses the scalars $(1-q_i^{-2})^{-1}$ instead of $(q_i^{-1} - q_i)^{-1}$, 
see \cite[Eq. (1.2.13)(a)]{Lusztig1}. For the Kashiwara form $(q_i^{-1} - q_i)^{-1}$ is replaced by $1$, 
and $e''_i$ are replaced by the skew-derivations $e'_i$ of $U_q(\n_-)$ satisfying  
$e'_i(x y) = e'_i(x) y + q^{(\al_i, \gamma)} x e'_i(y)$, 
see \cite[Eq. (3.4.4) and Proposition 3.4.4]{Kashiwara1}.

The use of the above form leads to minimal rescaling of dual PBW generators, quantum minors, and cluster variables.
\ere

Let $d \in \Zset_{> 0}$ be such that $(P\spcheck, P\spcheck) \subseteq \Zset/d$.
The Rosso--Tanisaki form \cite[\S 6.12]{Jantzen1}
\[
( -, -)_{RT} : U_q(\b_-) \times U_q(\b_+) \to \Qset(q^{1/d})
\]
is the Hopf algebra pairing satisfying
\begin{equation}
\label{RT-invar}
(x, y y')_{RT} = (\Delta(x), y' \otimes y)_{RT}, \quad (x x' , y)_{RT} = (x \otimes x', \Delta(y))_{RT}
\end{equation}
for $x,x' \in U_q(\b_-)$, $y, y' \in U_q(\b_+)$, and normalized by
\[
(f_i, e_j)_{RT} = \delta_{ij} (q_i^{-1} - q_i)^{-1}, 
\quad (q^h, q^{h'})_{RT} = q^{- (h , h') }, \quad (f_i, q^h)_{RT}= (q^h, e_i)_{RT}=0
\] 
for all $i, j \in I$, $h \in P\spcheck$. Its
restrictions to $U_q(\n_-) \times U_q(\b_+)$ and $U_q(\b_-) \times U_q(\n_+)$ take values in $\Qset(q)$. The above two forms are related by
\begin{equation}
\label{rel-forms}
(x, x')_{KL} = (x, \vp^*(x'))_{RT}, \quad \forall x, x' \in U_q(\n_-),
\end{equation}
see e.g. \cite[Lemma 3.8]{KO} or \cite[Proposition 8.3]{VY}.
\subsection{Integral forms and canonical bases}
Recall the notation \eqref{iA}. 
The (divided power) {\em{integral forms}} $U_q(\n_\pm)_\iA$ of $U_q(\n_\pm)$ are the $\iA$-subalgebras generated by $e_i^{(k)} := e_i^k/[k]_i!$
(resp. $f_i^{(k)} := f_i^k/[k]_i!$) for $i \in I$, $k \in \Zset_{>0}$. We have $\vp^*(U_q(\n_-)_\iA) = U_q(\n_+)_\iA$.
The {\em{dual integral form}} $U_q(\n_-)\spcheck_\iA$ of $U_q(\n_-)$ is the $\iA$-subalgebra
\begin{align}
\label{dual-U}
U_q(\n_-)\spcheck_\iA& = \{x \in U_q(\n_-) \mid ( x, U_q(\n_-)_\iA)_{KL} \subset \iA\} 
\\
&= \{x \in U_q(\n_-) \mid ( x, U_q(\n_+)_\iA)_{RT} \subset \iA\}.
\nn
\end{align}
 
Kashiwara \cite{Kashiwara1} defined a {\em{lower global basis}} ${\BB}^{\low}$ of $U_q(\n_-)_\iA$ 
and an {\em{upper global basis}} ${\BB}^{\up}$ of $U_q(\n_-)\spcheck_\iA$. 
The basis  ${\BB}^{\up}$ is defined from ${\BB}^{\low}$ as the dual basis with respect to the form $(-,-)_{KL}$. 
Lusztig \cite{Lusztig1} defined related {\em{canonical}} and {\em{dual canonical}} bases of $U_q(\n_+)_\iA$ 
and a dual integral form of $U_q(\n_+)$.  

\subsection{Quantum Schubert cell algebras, dual integral forms and CGL extensions}
\label{5.4}
To each $w \in W$, De Concini--Kac--Procesi \cite{DKP} and Lusztig \cite[\S 40.2]{Lusztig2}
associated {\em{quantum Schubert cell subalgebras}} of $U_q(\n_\pm)$. Given a reduced expression 
\begin{equation}
\label{reduced}
w = s_{i_1} \dots s_{i_N},
\end{equation}
define 
\[
w_{\leq k} := s_{i_1} \ldots s_{i_k}, \; \; w_{[j,k]} := s_{i_j} \ldots s_{i_k}, \; \; 
w_{\leq k}^{-1} := (w_{\leq k})^{-1}, \; \; w_{[j,k]}^{-1}:= (w_{[j,k]})^{-1}
\in W
\]
for $0 \le j \le k \le N$. Denote the roots and root vectors
\red{
\begin{equation}
\label{roots}
\be_k := w_{\leq k - 1} (\al_{i_k}), \; f_{\be_k} := T_{w_{\leq k-1}^{-1}}^{-1} (f_{i_k}) \in U_q(\n_-)_\iA, \;
e_{\be_k} := T_{w_{\leq k - 1}^{-1}}^{-1} (e_{i_k}) \in U_q(\n_+)_\iA
\end{equation}
}
for $k \in [1,N]$. 
The algebras $U_q(\n_\pm(w))$ are the unital $\Qset(q)$-subalgebras of $U_q(\n_\pm)$ generated by  $e_{\be_1}, \ldots, e_{\be_N}$ and
$f_{\be_1}, \ldots, f_{\be_N}$, respectively. These definitions are independent 
of the choice of reduced expression of $w$. Furthermore, 
\begin{equation}
\label{prod-Uw}
\begin{aligned}
U_q(\n_\pm(w)) &= U_q(\n_\pm) \cap T_{w^{-1}}^{-1} ( U_q(\n_\mp)),  \\
U_q(\n_\pm) &= \big( U_q(\n_\pm) \cap T_{w^{-1}}^{-1} ( U_q(\n_\pm)) \big) U_q(\n_\pm(w)).
\end{aligned}
\end{equation}
This was conjectured in \cite[Conjecture 5.3]{BerGreen} and proved in \cite{Kimura2,Tani}. 

Note that the algebras considered in \cite{DKP} (see also \cite{Jantzen1}) are 
\[
U_q^\pm[w] = * \big( U_q(\n_\pm(w)) \big).
\]
We use $U_q(\n_\pm(w))$ instead, to avoid making all algebras here antiisomorphic to the ones in \cite{GLS}.
The $\iA$-algebra
\[
U_q(\n_-(w))\spcheck_\iA := U_q(\n_-(w)) \cap U_q(\n_-)\spcheck_\iA 
\]
is called the {\em{dual integral form}} of $U_q(\n_-(w))$.
Define the {\em{dual PBW generators}} of $U_q(\n_-(w))$ 
\begin{equation}
\label{dual-PBWgen}
f_{\be_k}^* := \frac{1}{ (f_{\be_k}, e_{\be_k})_{RT}} f_{\be_k} = \frac{1}{ (\vp^*(e_{\be_k}), \vp^*(e_{\be_k}))_{KL}} \vp^*(e_{\be_k}) = 
(q_{i_k}^{-1}- q_{i_k}) f_{\be_k}
\end{equation}
for $k \in [1,N]$. Note that $\vp^*(e_{\be_k})$ differs from $f_{\be_k}$ by a unit of $\iA$, namely
$\vp^*(e_{\be_k}) = (- q_{i_k}) \prod_i (- q_i)^{n_i} f_{\be_k}$ 
where $n_i \in \Zset_{\geq 0}$ are such that $\be_k = \sum n_i \al_i$ (see e.g \cite[Eq. 8.14(9)]{Jantzen1}).
The inner products between the dual PBW monomials and the divided-power PBW monomials 
are given by 
\begin{equation}
\label{RT-mon}
\Big( (f_{\be_1}^*)^{m_1} \cdots (f_{\be_N}^*)^{m_N}, e_{\be_1}^{(l_1)} \cdots e_{\be_N}^{(l_N)} \Big)_{RT} = 
\prod_{k=1}^N \delta_{m_k l_k} q_{i_k}^{m_k(m_k-1)/2}, \quad \forall m_k, l_k \in \Zset_{\geq 0}
\end{equation}
(see e.g. \cite[\S 8.29-8.30]{Jantzen1}), where $e_{\be_k}^{(l_k)}:= e_{\be_k}^{l_k}/[l_k]_{i_k}$. 

\bth{basisUw} {\em{(Kimura)}} \cite[Prop. 4.26, Thms. 4.25 and 4.27]{Kimura1} The algebras $U_q(\n_-(w))\spcheck_\iA$ have the 
following decompositions as free $\iA$-modules:
\begin{align}
U_q(\n_-(w))\spcheck_\iA &= \bigoplus_{m_1, \ldots, m_N \in \Zset_{\geq 0}} \iA \cdot (f_{\be_1}^*)^{m_1} \cdots (f_{\be_N}^*)^{m_N}
\label{Int-PBW}
\\
&=\bigoplus_{d \in {\BB}^{\up} \cap U_q(\n_-(w))} \iA \cdot d.
\nn
\end{align}
The {\em{Levendorskii--Soibelman straightening law}} takes on the form
\begin{multline}
\label{LS}
f_{\be_k}^* f_{\be_j}^* - 
q^{ \red{( \be_k, \be_j )} }  
f_{\be_j}^* f_{\be_k}^*  \\
= \sum_{ {\bf{m}} = (m_{j+1}, \ldots, m_{k-1}) \, \in \, \Znn^{k-j-1} }
b_{\bf{m}} (f_{\be_{j+1}}^*)^{m_{j+1}} \ldots (f_{\be_{k-1}}^*)^{m_{k-1}},
\quad b_{\bf{m}} \in \iA
\end{multline}
for all $1 \le j < k \le N$. 
\eth
\bre{GLSrem1} Recall \eqref{HH-act} and denote
\begin{equation}
\label{t}
t :=( q_i^{-1} - q_i)_{i \in I} \in \HH.
\end{equation}
The objects associated to $U_q(\n_+)$ used by Gei\ss, Leclerc and Schr\"oer in \cite{GLS} are 
precisely the images under the isomorphism
\[
(t \cdot) \circ \vp^* : U_q(\n_-) \stackrel{\cong}{\lra} U_q(\n_+)
\]
of the objects associated to $U_q(\n_-)$ which we consider. Firstly, \cite{GLS} uses 
the canonical basis $\ol{\vp}( {\BB}^{\low}) = \vp^* ( {\BB}^{\low} )$ of $U_q(\n_+)$ and 
the PBW generators $e_{\be_k} = \ol{\vp}(f_{\be_k})$. They use the bilinear 
form $(-,-)$ on $U_q(\n_+)$ defined by 
\begin{equation}
\label{GLS-form}
(y, y'):= ( \vp^*(y), t^{-1} \vp^*(y'))_{KL}, \quad \forall y, y' \in U_q(\n_+)
\end{equation}
leading to the following:
\begin{enumerate}
\item The dual canonical basis of $U_q(\n_+)$ constructed from the canonical basis $\ol{\vp}( {\BB}^{\low}) = \vp^* ( {\BB}^{\low} )$ and the bilinear form \eqref{GLS-form}, 
thus giving the basis $(t \cdot) \circ \vp^*({\BB}^{\up})$;  
\item The dual PBW generators $e^*_{\be_k} := e_{\be_k} /(e_{\be_k}, e_{\be_k}) = (t \cdot) \circ \vp^*(f^*_{\be_k})$ of $U_q(\n_+(w))$; 
\item The dual integral forms $(t \cdot) \circ \vp^*(U_q(\n_-)\spcheck_\iA) = \{y \in U_q(\n_+) \mid ( x, U_q(\n_+)_\iA) \subset \iA\} $ of $U_q(\n_+)$
and $(t \cdot) \circ \vp^*(U_q(\n_-(w))\spcheck_\iA)$ of $U_q(\n_+(w))$. \qed
\end{enumerate}
\ere

For a reduced expression \eqref{reduced} of $w \in W$ and $k \in [1,N]$, fix elements $t_k, t^\sy_k \in \HH$ such that
\begin{equation}
\label{tk}
t_k^{\be_j} = 
q^{ \red{(\be_k, \be_j)} } \; \; 
\mbox{for} \; \;  j \in [1, k]
\quad \mbox{and} \quad
(t^\sy_k)^{\be_l} = 
q^{-(\be_k, \be_l) } \; \; 
\mbox{for} \; \;  l \in [k, N],
\end{equation}
cf. \eqref{HH-act}; such $t_k, t^\sy_k$ exist but are not unique since the restriction 
of the form $(.,.)$ to $Q$ is degenerate when $\g$ is not finite dimensional.
Note that the algebras $U_q(\n_-(w))\spcheck_\iA$ are preserved by the automorphisms 
$(t_k \cdot), (t_k^* \cdot)$. 

\ble{Uw} Let $w \in W$, \eqref{reduced} a reduced expression of $w$, and $t_k \in \HH$ satisfying \eqref{tk}.
\begin{enumerate}
\item[(a)] For $k \in [1,N]$, the algebra $U_q(\n_-(w_{\leq k}))\spcheck_\iA$ is an Ore extension 
\[
U_q(\n_-(w_{\leq k}))\spcheck_\iA \cong U_q(\n_-(w_{\leq k-1}))\spcheck_\iA [f_{\be_k}^*; (t_k \cdot), \delta_k],
\]
where $\delta_k$ is the locally nilpotent $(t_k \cdot)$-derivation
of $U_q(\n_-(w_{\leq k-1}))\spcheck_\iA$ given by
\[
\delta_k(x) := f_{\be_k}^* x - q^{\red{(\be_k, \wtt x)}} x f_{\be_k}^* \quad \mbox{for homogeneous} 
\quad
x \in U_q(\n_+(w_{\leq k-1}))\spcheck_\iA.
\]
The $t_k$-eigenvalue of $f_{\be_k}^*$ equals $q^2_{i_k}$, which is not a root of unity.

\item[(b)] The algebra 
\begin{equation}
\label{UwCGL}
U_q(\n_-(w)) \cong \Qset(q) [f_{\be_1}^*] [f_{\be_2}^*; (t_2 \cdot), \delta_2] \cdots [f_{\be_N}^*; (t_N \cdot), \delta_N]
\end{equation}
is a symmetric CGL extension. The algebra
\begin{equation}
\label{UwCGLint}
U_q(\n_-(w))\spcheck_\iA \cong \iA [f_{\be_1}^*] [f_{\be_2}^*; (t_2 \cdot), \delta_2] \cdots [f_{\be_N}^*; (t_N \cdot), \delta_N]
\end{equation}
with the generators $f_{\be_1}^*, \ldots, f_{\be_N}^*$ is an $\iA$-form of the CGL extension \eqref{UwCGL}.
\item[(c)] The interval subalgebras  of $U_q(\n_-(w))\spcheck_\iA$ are
\begin{equation}
\label{interval}
(U_q(\n_-(w))\spcheck_\iA)_{[j,k]} = T^{-1}_{w_{\leq j-1}^{-1}} \bigl( U_q(\n_-(w_{[j,k]}))\spcheck_\iA \bigr) \quad\mbox{for} \quad 1 \le j \le k \le N.
\end{equation}
\end{enumerate}
\ele

\begin{proof} Part (a) follows from \eqref{Int-PBW} and \eqref{LS}. 

(b) The facts that $U_q(\n_-(w))$ is a CGL extension and that $U_q(\n_-(w))\spcheck_\iA$ with the generators 
$f_{\be_1}^*, \ldots, f_{\be_N}^*$ is an $\iA$-form of it
follow by iterating (a). Its symmetricity is proved 
analogously to (a). 

(c) Applying twice \eqref{Int-PBW} and using \eqref{dual-PBWgen}, we obtain 
\begin{align*} 
&(U_q(\n_-(w))\spcheck_\iA)_{[j,k]} =  \oplus_{m_j, \ldots, m_k \in \Zset_{\geq 0}} \iA \cdot (f_{\be_j}^*)^{m_j} \cdots (f_{\be_k}^*)^{m_k} \\
&= T^{-1}_{w_{\leq j-1}^{-1}}  \big( \oplus_{m_j, \ldots, m_k \in \Zset_{\geq 0}} \iA \cdot ((q_{i_j}^{-1} -q_{i_j}) f_{i_j})^{m_j} \cdots 
((q_{i_k}^{-1} -q_{i_k})  T_{w_{[j,k-1]}^{-1}}^{-1} f_{i_k})^{m_k} \big) \\
&=  T_{w_{\leq j-1}^{-1}}^{-1} \bigl( U_q(\n_-(w_{[j,k]}))\spcheck_\iA \bigr),
\end{align*} 
which proves \eqref{interval}.
\end{proof}
An important feature of the normalization of $(-,-)_{KL}$ is that there are no additional scalars in \leref{Uw}(c) due to the 
braid group action.
\subsection{The quantum function algebra of $\g$}
Consider the full dual $\Qset(q)$-vector space $U_q(\g)^*$ which is canonically a unital algebra 
using the coproduct and counit of $U_q(\g)$. It is a $U_q(\g)$-bimodule by
\begin{equation}
\label{actions}
\lcor x \cdot c \cdot y, z \rcor := \lcor c , y z x \rcor \quad 
\mbox{for} \quad c \in U_q(\g)^*, \, x, y, z \in U_q(\g).
\end{equation}
For a right $U_q(\g)$-module $V$, let $V^\vp$ be the left $U_q(\g)$-module structure on the vector space $V$
such that
\[
x \cdot v = v \cdot \vp(x) \quad \mbox{for} \quad v \in V, \ x \in U_q(\g).
\]
For each $\mu \in P_+$, there exists a unique irreducible right $U_q(\g)$ module $V^{\mathrm{r}}(\mu)$ such that 
$V^{\mathrm{r}}(\mu)^\vp \cong V(\mu)$. Analogously to $\Ointg$, one defines an $\OO$-type category 
of integrable right $U_q(\g)$-modules; it is denoted by $\OO_{\mathrm{int}}(\g^{\mathrm{op}})$.

Kashiwara defined \cite[Sect. 7]{Kashiwara5} the quantized coordinate ring $A_q(\g)$ of the Kac--Moody group 
of $\g$ as the unital subalgebra of $U_q(\g)^*$ consisting of those $f \in U_q(\g)^*$ such that 
\[
U_q(\g) \cdot f \in \Ointg \quad \mbox{and} \quad
f \cdot U_q(\g) \in \OO_{\mathrm{int}}(\g^{\mathrm{op}}).
\]
Kashiwara also proved \cite[Proposition 7.2.2]{Kashiwara5} a quantum version of the Peter-Weyl theorem that there is an isomorphism
of $U_q(\g)$-bimodules
\begin{equation}
\label{PW}
A_q(\g) \cong \bigoplus_{\mu \in P_+} V^{\mathrm{r}}(\mu) \otimes V(\mu).
\end{equation}

For $M \in \Ointg$ and $v \in M$, $\xi \in D_\vp M$ define the matrix coefficient 
\begin{equation}
\label{matr-coeff}
c_{\xi v} \in U_q(\g)^* \quad \mbox{given by} \quad
\lcor c_{\xi v}, x \rcor := \lcor \xi, x \cdot v \rcor \; \; \forall x \in U_q(\g).
\end{equation}
It follows from \eqref{PW} that
\[
A_q(\g) = \{ c_{\xi v} \mid M \in \Ointg, \, v \in M, \, \xi \in D_\vp M \}
= \oplus_{\mu \in P_+} \{ c_{\xi v} \mid v \in V(\mu), \, \xi \in D_\vp V(\mu) \}.
\]
This is the form in which quantum function algebras were defined in the finite dimensional case \cite{LS}. 
The algebra $A_q(\g)$ is $P \times P$-graded by
\begin{equation}
\label{Aq-grad}
A_q(\g)_{\mu, \nu} = \{c_{\xi v} \mid \xi \in (V_\mu)^* \subset D_\vp V, \, v \in V_\nu, \, V \in \Ointg \}, 
\quad \forall \mu, \nu \in P.
\end{equation}

\sectionnew{Homogeneous prime ideals of $A_q(\n_+(w))$}  \label{Hprime-sec}
\subsection{The algebras $A_q(\n_+)$ and $A_q(\n_+(w))$}
\label{6.1}
It follows from the first identity in \eqref{RT-invar} and the nondegeneracy of the form 
that  the map
\begin{equation}
\label{iota}
\iota : U_q(\n_-) \to U_q(\b_+)^*  \quad \mbox{given by} \quad  \lcor \iota(x), y \rcor = (x, y)_{RT}, \; \; \forall 
x \in  U_q(\n_-), \, y \in U_q(\b_+)
\end{equation}
is an injective algebra homomorphism. Here $U_q(\b_+)^*$ denotes the unital algebra which is the full dual 
of the Hopf algebra $U_q(\b_+)$ over $\Qset(q)$.  

Following Gei\ss--Leclerc--Schr\"oer \cite[\S 4.2]{GLS}, denote
the subalgebra $A_q(\n_+) \subset U_q(\b_+)^*$ consisting of those $f \in U_q(\b_+)^*$ such that 
\begin{enumerate}
\item[(i)] $f(x q^h) = f(x)$ for all $x \in U_q(\n_+)$, $h \in P\spcheck$ and 
\item[(ii)] $f(x)=0$ for all $x \in U_q(\n_+)_\ga$ and $\ga \in Q_+ \backslash S$ for a finite subset $S$ of $Q_+$. 
\end{enumerate}
The properties
\[
(x q^h, y q^{h'})_{RT} = (x, y)_{RT} q^{ - (h,h')}, \quad 
(U_q(\n_-)_{-\ga}, U_q(\n_+)_{\de} )_{RT} =0 
\]
for $x \in U_q(\n_-)$, $y \in U_q(\n_+)$, $h, h' \in P\spcheck$, $\ga \neq \de$ in $Q_+$
(see \cite[Eq 6.13(1)]{Jantzen1}) and the nondegeneracy of $(-,-)_{RT}$ imply that 
$A_q(\n_+):= \iota(U_q(\n_-))$. Thus
\begin{equation}
\label{iota-isom1}
\iota : U_q(\n_-) \stackrel{\cong}{\lra} A_q(\n_+)
\end{equation}
is an algebra isomorphism. Following \cite[\S 7.2]{GLS}, define $A_q(\n_+(w)):= \iota( U_q(\n_-(w)))$. Hence, $\iota$ restricts 
to the algebra isomorphsim
\begin{equation}
\label{iota-isom2}
\iota :  U_q(\n_-(w)) \stackrel{\cong}{\lra} A_q(\n_+(w)).
\end{equation}
Using the isomorphism $\iota$, transport the isomorphisms $T_w : U_q(\n_-) \cap T_w^{-1}( U_q(\n_-)) \to T_w( U_q(\n_-)) \cap U_q(\n_-)$ 
to such maps on $A_q(\n_+)$. Denote the integral forms over $\iA$
\[
A_q(\n_+)_\iA:= \iota( U_q(\n_-)\spcheck_\iA) \quad \mbox{and} \quad
A_q(\n_+(w))_\iA:= \iota( U_q(\n_-(w))\spcheck_\iA )
\]
of $A_q(\n_+)$ and $A_q(\n_+(w))$. The algebra $A_q(\n_+(w))$ is $Q_+$-graded by
\[
A_q(\n_+(w))_\ga := \iota ( U_q(\n_-(w))_{- \ga}), \quad \forall \ga \in Q_+.
\]
In other words, the isomorphism \eqref{iota-isom2} is not $\HH$-equivariant, but satisfies 
$\iota (t \cdot u) = t^{-1} \cdot \iota(u)$ for $t \in \HH$, $u \in U_q(\n_-(w))$. 
\bre{GLSrem2} Using the bilinear form \eqref{GLS-form}, in \cite{GLS} the algebra $A_q(\n_+)$ is identified with $U_q(\n_+)$ via the isomorphism
\[
\Psi \colon U_q(\n_+) \stackrel{\cong}{\lra} A_q(\n_+), \quad 
\lcor \Psi(y), y' q^h \rcor := (y, y')_{KL}, \quad \forall y,y' \in U_q(\n_+), \ h \in P\spcheck.
\]
$\Psi$ fits in the commutative diagram
\[
\begin{tikzpicture}[scale=1.6]
\node (A) at (0.75, 0.87) {$A_q(\n_+)$};
\node (C) at (0,0) {$U_q(\n_-)$};
\node (D) at (1.5,0) {$U_q(\n_+)$};
\path[->,font=\scriptsize,>=angle 90]
(C) edge node[left]{$\iota$} (A)
(D) edge node[right]{$\Psi$} (A)
(C) edge node[above]{$(t \cdot) \circ \vp^*$} (D);
\end{tikzpicture}
\]
in terms of $t \in \HH$ given by \eqref{t}.
This and \reref{GLSrem1} imply that $A_q(\n_+)_\iA$ and $\iota({\BB}^{\up})$ are precisely the integral form of $A_q(\n_+)$ and the 
dual canonical basis of $A_q(\n_+)$ considered in \cite{GLS}. However the braid group action of \cite{GLS} on $A_q(\n_+)$ 
is a conjugate of ours by an element of the torus $\HH$, and involves extra scalars compared to our formulas. 
\ere

\subsection{An algebra isomorphism}
For $\mu \in P_+$, fix a highest weight vector $v_\mu$ of $V(\mu)$. For $w \in W$, define 
the extremal weight vector 
\[
v_{w \mu} := T_{w^{-1}}^{-1} v_\mu \in V(\mu)_{w \mu}.
\]
Denote the associated {\em{Demazure modules}}
\[
V_w^\pm(\mu) := U_q(\b_\pm) v_{w \mu} \subseteq V(\mu).
\]
Let
\[
\xi_{w \mu} \in V(\mu)_{w \mu}^* \subset D_\vp (V(\mu)) \quad \mbox{be such that} \quad \lcor \xi_{w \mu}, v_{w \mu} \rcor =1.
\]
For $u, w \in W$ and $\mu \in P_+$, using the notation \eqref{matr-coeff}, define the quantum minors 
\[
\De_{u \mu, w \mu} := c_{\xi_{u \mu}, v_{w \mu}} \in A_q(\g),
\]
which are equivalently given by \cite[Eq. (9.10)]{BZ}, \cite[Eq. (3.5)]{GLS}.
It is well known that 
\begin{equation}
\label{prod-Tw}
T_{w^{-1}}^{-1} (v_\mu \otimes v_\nu) = T_{w^{-1}}^{-1} v_\mu \otimes T_{w^{-1}}^{-1} v_\nu
\end{equation}
for all $\mu, \nu \in P_+$. This implies that 
\begin{equation}
\label{DeDe}
\De_{u \mu, w \mu} \De_{u \nu, w \nu} = \De_{u (\mu + \nu), w (\mu + \nu)}, \quad 
\forall \mu, \nu \in P_+.
\end{equation}

Following Joseph \cite[\S 9.1.6]{Joseph-book}, denote the subalgebra
\[
A^+_q(\g) := \oplus_{\mu \in P_+} \{ c_{\xi v_\mu} \mid \xi \in D_\vp ( V(\mu) )\} 
\] 
of $A_q(\g)$. By \cite[Lemma 2.1(i)]{FY}, the multiplicative set
\[
E_w := \{ \De_{w \mu, \mu} \mid \mu \in P_+ \}  
\]
is a denominator set in $A_q^+(\g)$. Denote the subsets 
\[
J^\pm_w := \oplus_{\mu \in P_+} \{ c_{\xi v_\mu} \mid \xi \in D_\vp ( V(\mu) ), \ \xi \perp V^\pm_w(\mu) \} \subset A_q^+(\g).
\]
By the proofs of Theorems \ref{thom1} and \ref{thom2} below, they are completely prime ideals of $A_q^+(\g)$.

The $P \times P$-grading \eqref{Aq-grad} of $A_q(\g)$ extends to a $P \times P$-grading 
of the localization $A_q^+(\g)[E_w^{-1}]$. For a graded subalgebra $R \subseteq A_q^+(\g)[E_w^{-1}]$, denote the subalgebra
\[
R_0:= \oplus_{\nu \in P} R_{\nu, 0},
\]
noting that $R_0$ is naturally $P$-graded.
It is easy to show that every element of $(A_q^+(\g)[E_w^{-1}])_0$ has the form $c_{\xi, v_\mu} \De_{w \mu, \mu}^{-1}$ for some
$\mu \in P_+$, $\xi \in D_\vp ( V(\mu) )$; in particular, this algebra is $Q$-graded.
The following theorem was proved in the finite dimensional case in \cite{Y-plms}. 
 
\bth{hom1} For all symmetrizable Kac--Moody algebras $\g$ and $w \in W$, there exists a $Q$-graded surjective 
homomorphism $\psi_w : ( A_q^+(\g) [E_w^{-1}])_0 \to A_q(\n_+(w))$ such that
\begin{equation}
\label{psi-w}
\lcor \psi_w( c_{\xi, v_\mu} \De_{w \mu, \mu}^{-1}), y q^h \rcor = \lcor \xi, y v_{w \mu} \rcor
\end{equation}
for $\mu \in P_+$, $\xi \in D_\vp ( V(\mu) )$, $y \in U_q(\b_+)$, $h \in P\spcheck$.
Its kernel equals $(J_w^+[E_w^{-1}])_0$. 
\eth
We will need the following lemma.

\ble{hom} \cite[Lemma 3.2]{Y-plms} Let $H$ be a Hopf algebra over $\KK$ and $A$ be an $H$-module algebra equipped with a right $H$-action.
For every  algebra homomorphism $\theta : A \to \KK$, the map $\psi : A \to H^*$, given by
\[
\psi(a) (h) = \theta (a  \cdot h),
\]
is an algebra homomorphism.
\ele

\begin{proof}[Proof of Theorem {\rm\ref{thom1}}] Eq.~\eqref{prod-Tw} implies that 
\[
\theta_w : A_q^+(\g) \to \Qset(q) \quad \mbox{given by} \quad \theta_w( c_{\xi v_\mu}) := \lcor \xi, v_{w \mu} \rcor, \; \; 
\forall \mu \in P_+, \ \xi \in D_\vp ( V(\mu) )
\]
is an algebra homomorphism. We apply the lemma to it and to the right action \eqref{actions} of $U_q(\b_+)$ on $A_q^+(\g)$. 
It shows that the map $\psi_w : A_q^+(\g) \to U_q(\b_+)^*$, given by
\[
\lcor \psi_w( c_{\xi v_\mu}), y \rcor := \lcor \xi, y v_{w \mu} \rcor, \quad \forall \mu \in P_+, \ \xi \in D_\vp ( V(\mu) ), \ y \in U_q(\b_+),
\]
is an algebra homomorphism. The element $\psi_w( \De_{w \mu, \mu})$ is a unit of $U_q(\b_+)^*$ because
\[
\lcor \psi_w( \De_{w \mu, \mu}),y q^h \rcor = \epsilon(y) q^{ \lcor h, \mu \rcor}, \quad 
\forall y \in U_q(\n_+), \ h \in P\spcheck.
\]
Hence, $\psi_w$ extends to $A_q^+(\g)[E_w^{-1}]$, $\psi_w ( (A_q^+(\g)[E_w^{-1}])_0) \subset A_q(\n_+)$, and the restriction
of $\psi_w$ to $(A_q^+(\g)[E_w^{-1}])_0$ is given by \eqref{psi-w}. From now on we will denote by $\psi_w$ this restriction. 
The formula \eqref{psi-w} implies at once that the kernel of $\psi_w$ equals $(J_w^+[E_w^{-1}])_0$ and
\[
\lcor \Im\, \psi_w, U_q(\n_+(w)) y \rcor =0, \quad 
\forall y \in \big( U_q(\n_+) \cap T_{w^{-1}}^{-1} (U_q(\n_+)) \big)_\ga, \ \ga \in Q_+ \backslash \{0\}. 
\]
For each $\ga \in Q_+$ such that $U_q(\n_+(w))_\ga \neq 0$, there exists $\mu \in P_+$ such that the pairing 
\[
(V_w(\mu)_{\ga + w \mu})^* \times U_q(\n_+(w))_\ga \quad 
\mbox{given by} \quad
\xi, y \mt \lcor \xi, y v_{w\mu} \rcor
\]
is nondegenerate. This, the second equality in \eqref{prod-Uw} and the fact that
\[
(U_q(\n_-(w)), U_q(\n_+(w)) y)_{RT}= 0, \quad 
\forall y \in \big( U_q(\n_+) \cap T_{w^{-1}}^{-1} (U_q(\n_+)) \big)_\ga, \ \ga \in Q_+ \backslash \{0\}
\]
imply that $\Im\, \psi_w = A_q(\n_+(w))$. 
\end{proof}

\bth{hom2} In the setting of Theorem {\rm\ref{thom1}}, there exists a {\em{(}}$Q$-graded{\em{)}} homomorphism 
$\psi_w^- : ( A_q^+(\g) [E_w^{-1}])_0 \to U_q(\b_-)^*$ such that
\[
\lcor \psi_w(c_{\xi, v_\mu} \De_{w \mu, \mu}^{-1}), y q^h \rcor = \lcor \xi, y^* v_{w \mu} \rcor
\]
for $\mu \in P_+$, $\xi \in D_\vp ( V(\mu) )$, $y \in U_q(\b_-)$, $h \in P\spcheck$. 
Its kernel equals $(J_w^-[E_w^{-1}])_0$. Its image is
contained in the image of the antiembedding $U_q(\n_+(w)) \to (U_q(\b_-))^*$ 
coming from the second component of the Rosso-Tanisaki form.
\eth

The proof of the theorem is analogous to that of \thref{hom1}. 

\subsection{The prime spectrum of $A_q(\n_+(w))$} The fact that $\Ointg$ is a braided monoidal category gives rise 
to ${\mathcal{R}}$-matrix commutation relations in $A_q(\g)$, \cite[Proposition 9.1.5]{Joseph-book}. Particular cases of those are 
the relations
\begin{equation}
\label{De-Aq}
\De_{w \mu, \mu} x = q^{\pm( (w \mu, \nu)  - (\mu, \ga) )} x \De_{w \mu, \mu} \mod J_w^\pm, \quad \forall x \in A_q^+(\g)_{\nu, \ga}, \ \mu \in P_+, \ \nu, \ga \in P.
\end{equation}
For $u \in W$, $\mu \in P_+$, denote the unipotent quantum minors
\[
D_{u \mu, w \mu} := \psi_w( \De_{u \mu, \mu} \De_{w \mu, \mu}^{-1}) \in A_q(\n_+(w)).
\]
They are alternatively defined as the elements of $A_q(\n_+(w))_{(u-w)\mu} \subset A_q(\n)$ such that 
\begin{equation}
\label{D-unip}
\lcor D_{u \mu, w \mu}, x q^h \rcor = \lcor \xi_{u \mu}, x v_{w \mu} \rcor, 
\quad \forall x \in U_q^+(\g), \ h \in P\spcheck,
\end{equation}
which implies that they are precisely the elements of $A_q(\n_+(w))$ defined in \cite[Eqs. (5.3)-(5.4)]{GLS}. 
Set 
\[
W^{\leq w} = \{ u \in W \mid u \leq w \}.
\]
For $u \in W^{\leq w}$, denote the ideals
\[
I_w(u) := \psi_w \big( (J_u^-[E_w^{-1}])_0 \big) \quad \mbox{of} \quad
A_q(\n_+(w)).
\]
It follows from \eqref{DeDe} and \eqref{De-Aq} that
\[
D_{u \mu, w \mu} D_{u \nu, w \nu} = q^{(w \mu, u \nu)- (\mu,\nu)} D_{u (\mu+ \nu), w (\mu+ \nu)}, \quad
\forall \mu, \nu \in P_+
\]
and that
\[
D_{u \mu, w \mu} x = q^{((w+u)\mu, \wtt x)} x  D_{u \mu, w \mu} \mod I_w(u),
\quad
\forall \mu \in P_+, \ \text{homogeneous} \ x \in A_q(\n_+(w)). 
\]
We have $I_w(1)=0$, thus
\begin{equation}
\label{D-normal}
D_{\mu, w \mu} x = q^{((w+1)\mu, \wtt x)} x  D_{\mu, w \mu},
\quad
\forall \mu \in P_+, \ \text{homogeneous} \ x \in A_q(\n_+(w)).
\end{equation}
Denote the multiplicative sets
\[
E_w(u) := q^\Zset \{ D_{u \mu, w \mu} \mid \mu \in P_+ \} 
\quad \mbox{in} \quad
A_q(\n_+(w)).
\]
Analogously to \eqref{HH-act}, we use the $Q_+$-grading of 
$A_q(\n_+(w))$ to construct an action of the torus $\HH$ on it.
\bth{spec-Uw} For all symmetrizable Kac--Moody algebras $\g$ and $w \in W$, the following hold:
\begin{enumerate}
\item[(a)] The graded prime ideals of $A_q(\n_+(w))$ are the ideals $I_w(u)$ for $u \in W^{\leq w}$. 
The map $u \mt I_w(u)$ is an isomorphism of posets from $W^{\leq w}$ with the Bruhat order
to the set of graded prime ideals of $A_q(\n_+(w))$ with the inclusion order. 
\item[(b)] All prime ideals of $A_q(\n_+(w))$ are completely prime and 
\[
\Spec A_q(\n_+(w)) = \bigsqcup_{u \in W^{\leq w}} \Spec_u A_q(\n_+(w)),
\]
where
$\Spec_u A_q(\n_+(w)) := \{ \JJ \in \Spec A_q(\n_+(w)) \mid \cap_{t \in \HH} ( t \cdot \JJ ) = I_w(u) \}.$
\end{enumerate}
The following hold for $u \in W^{\leq w}$:
\begin{enumerate}
\item[(c)] $I_w(u) \cap E_w(u) = \varnothing$ and the localization 
$R_{u,w} = \big( A_q(\n_+(w))/I_w(u) \big)[E_w(u)^{-1}]$ is an $\HH$-simple domain. 
\item[(d)]  For $u \in W^{\leq w}$,  the center $Z(R_{u,w})$ is a Laurent polynomial ring over $\Qset(q)$ and
there is a homeomorphism 
\[
\eta_u \colon \Spec Z(R_{u,w}) \stackrel{\cong}{\lra} \Spec_u A_q(\n_+(w))
\]
where for ${\mathcal{J}} \in \Spec Z(R_{u,w})$, $\eta_u(\JJ)$ is the ideal of $A_q(\n_+(w))$ containing 
$I_w(u)$ such that 
$\eta_u(\JJ) /I_w(u) = \JJ R_{u,w} \cap (A_q(\n_+(w)) /I_w(u))$.
\end{enumerate}
\eth
Denote for brevity the algebra
\[
A_w^+:= ( A_q^+(\g) [E_w^{-1}])_0 \subseteq  A_q^+(\g) [E_w^{-1}].
\]
It is $Q$-graded by 
\[
(A_w^+)_\nu := ( A_q^+(\g) [E_w^{-1}])_{\nu, 0} \quad \mbox{for} \quad \nu \in Q
\]
in terms of the $P \times P$-grading \eqref{Aq-grad} of $A_q^+(\g) [E_w^{-1}]$.
Define the commuting (inner) automorphisms $\tau_w^\mu \in \Aut (A_w^+)$ for $\mu \in P_+$ by
\[
\tau_w^\mu( c ) : = \Delta_{w \mu, \mu}^{-1}  c  \Delta_{w \mu, \mu}. 
\]
For each $i \in I$, define the automorphism $\kappa_i \in \Aut (A_q^+(\g))$ by $\kappa_i(c) := c \cdot q^{d_i h_i}$ and 
the locally nilpotent (right skew) $\kappa_i$-derivation $\partial_i$ of $A_q^+(\g)$ by $\partial_i(c) := c \cdot f_i$ in terms of the second action in 
\eqref{actions}. It easy to check that $\kappa_i \partial_i \kappa_i^{-1} = q_i \partial_i$. 
Following Joseph \cite[\S A.2.9]{Joseph-book}, for $c \in A_q^+(\g) \backslash \{0 \}$ set
\[
\deg_i (c) := \max \{ n \in \Zset_{> 0} \mid \partial_i^n (c) \neq 0 \}
\]
and
\[
\partial^*_w(c) := \partial_{i_1}^{n_1} \ldots \partial_{i_N}^{n_N} (c) \neq 0
\]
where $n_N, \ldots, n_1 \in \Zset_{\geq 0}$ are recursively defined by $n_k := \deg_{i_k} ( \partial_{i_{k+1}}^{n_{k+1}} \ldots \partial_{i_N}^{n_N} (c))$ 
in terms of the reduced expression \eqref{reduced}. Set $\partial_w^*(0) := 0$.

\begin{proof} We carry out the proof in four steps as follows:

{\em{Step 1.}} {\em{For all $u \in W^{\leq w}$, the ideals $I_w(u)$ of $A_q(\n_+(w))$ are completely prime.}} 

The image of $\psi_w$ is an iterated skew polynomial extension, and thus is a domain. Similarly one shows that the 
image of $\psi^-_w$ is also a domain. Therefore $(J^\pm_w[E_w^{-1}])_0$ are completely prime ideals of $A_w^+$. 
By direct extension and contraction arguments one gets that $J^\pm_u$ are completely prime ideals of $A_q^+(\g)$ 
for $u \in W$, and that the same is true for the ideals $(J^\pm_u[E_w^{-1}])_0$ of $A_w^+$.
The remaining part of the proof of the statement of step 1 uses elements of Gorelik's and Joseph's proofs \cite{Gorelik,Joseph-book} 
of related facts in the finite dimensional case. 
We prove the stronger fact that there exists an embedding 
\[
A_q(\n_+(w))/ I_w(u) \hra A_w^+ /  (J^-_w[E_w^{-1}])_0
\]
which we construct next. For a linear map $\tau$ on a $\Qset(q)$-vector space $V$ and $t \in \ol{\Qset(q)}$, denote by 
$\EE_\tau(t)$ the generalized $t$-eigenspace of $\tau$. Using the first action \eqref{actions}, one shows that for all
$i \in I$, $w \in W$ such that $\ell(s_i w) < \ell(w)$ and $\nu \in P$, $\la \in P_+$, $t \in \ol{\Qset(q)}$:
\[
\mbox{If} \quad c \in \EE_{\tau_w^\mu}( t) \cap A_q^+(\g)_{\nu, \la}, \quad \mbox{then} \quad
\partial_i^n (c) \in \EE_{\tau_{s_i w}^\mu}( t q^{ (w \mu, \nu) - (s_i w \mu, \nu + n \al_i)} )
\]
where $n:= \deg_i (c)$; the proof of this is analogous to \cite[Lemma 6.3.1]{Gorelik}. By induction on the length of $w$, this 
implies that
\begin{multline}
\label{Adecomp}
A_w^+ = \oplus_{\ga \in Q_+} (A_w^+)[2 \ga] \quad \mbox{where}  \\
\quad (A_w^+)[2 \ga] : = \moplus_{\nu \in Q} \{ c \in (A_w^+)_\nu \mid c \in \EE_{\tau_w^\mu}( q^{ - (w^{-1} \nu + 2 \ga, \mu)}), \; \forall \mu \in P_+ \},
\end{multline}
and that for $\ga \in Q_+$, $\la \in P_+$, 
\begin{equation}
\label{Adecmp-emb}
c_{\xi, v_\la} \De_{w \la, \la}^{-1} \in (A_w^+)[2 \ga] \quad \Rightarrow \quad 
(\partial^*_{w^{-1}} (c_{\xi, v_\la})) \De_{\la, \la}^{-1} \in (A_1^+)[2 \ga].
\end{equation}
The base of the induction for $w=1$ follows from \eqref{De-Aq} applied to $J_1^-=0$, which gives that 
$\tau_1^\mu(c) = q^{(\mu, \nu)} c$ for all $c \in (A_q^+(\g)[E_1^{-1}])_{\nu, 0}$, $\nu \in - Q_+$; that is
\begin{equation}
\label{A1gamma}
(A_1^+)[2 \ga] = \moplus_{\la \in P_+} \{ c_{\xi, v_\la} \De_{\la, \la}^{-1} \mid \xi \in  ( V(\la)_{\la - \ga}) )^* \subset D_\varphi ( V(\la) ) \}, \quad \forall \ga \in Q_+.
\end{equation}

Furthermore, we have 
\begin{equation}
\label{Jw-decomp}
(J^+_w[E_w^{-1}])_0 =  \oplus_{\ga \in Q_+ \backslash \{0 \}} A_w^+ [2 \ga]. 
\end{equation}
By \eqref{De-Aq}, applied to $J_w^+$, the right hand side is contained in the left one. Because of \eqref{Adecomp} 
it remains to show that $(J^+_w[E_w^{-1}])_0 \cap A_w^+[0] =0$. Assume the opposite, that 
$(J^+_w[E_w^{-1}])_0 \cap A_w^+[0] \neq 0$. By putting elements over a common denominator, 
each element of $A_w^+$ can be represented in the form $c_{\xi, v_\la} \De_{w \la, \la}^{-1}$ for some $\la \in P_+$, $\xi \in D_\varphi (V(\la))$. Choose 
a nonzero element of this form in $(J^+_w[E_w^{-1}])_0 \cap A_w^+[0]$. By \eqref{Adecmp-emb}, 
$(\partial^*_{w^{-1}} (c_{\xi, v_\la})) \De_{\la, \la}^{-1} \in (A_1^+)[0]$. 
Hence, \eqref{A1gamma} implies that 
$\partial^*_{w^{-1}} (c_{\xi, v_\la}) = r c_{\xi_\la, v_\la}$ for some $r \in \Qset(t)^*$. The definition of $\partial_{w^{-1}}^*$ gives that 
\[
\lcor \xi, f_{i_1}^{n_1} \ldots f_{i_N}^{n_N} v_\la \rcor \neq 0 \quad \mbox{for some} \quad 
n_1, \ldots, n_N \in \Zset_{\geq 0}
\] 
in terms of the reduced expression \eqref{reduced}. However, $f_{i_1}^{n_1} \ldots f_{i_N}^{n_N} v_\la \in V_w^+(\la)$
by the standard presentation of Demazure modules \cite[Lemma 4.4.3(v)]{Joseph-book}. This contradicts with 
$c_{\xi, v_\la} \De_{\la, \la}^{-1} \in (J_w^+[E_w^{-1}])_0$ and proves \eqref{Jw-decomp}. 

Since $\tau_w^\mu \in \Aut (A_w^+)$, $A_w^+[0]$ is a subalgebra of $A_w^+$. \thref{hom1} and \eqref{Adecomp}, \eqref{Jw-decomp} imply
\begin{align*}
A_q(\n_+(w))/I_w(u) &\cong A_w^+/( (J_w^+[E_w^{-1}])_0 + (J_u^-[E_w^{-1}])_0) 
\\
&\cong A_w^+[0]/ (A_w^+[0] \cap (J_u^-[E_w^{-1}])_0) 
\hra A_w^+/ (J_u^-[E_w^{-1}])_0. 
\end{align*}

{\em{Step 2.}} {\em{For all $u \in W^{\leq w}$, $I_w(u) \cap E_w(u) = \varnothing$.}} 

Denote by $G^{\min}$ the minimal Kac--Moody group associated to $\g$, see \cite[\S 7.4]{Ku} for details. 
Let $H$ be the Cartan subgroup of $G^{\min}$, and $N^{\min}_+$ and $N_-$ the subgroups of $G^{\min}$ generated by 
its one-parameter unipotent subgroups for positive and negative roots, respectively. 
Denote by ${\mathcal{B}}^{\min}_+$ and $\mathcal{{B}}_-$ the associated Borel subgroups of $G^{\min}$.
Denote by $N_+(w)$ the unipotent subgroup of $N^{\min}_+$ corresponding to $\n_+(w)$. By \cite[Theorem 4.44]{Kimura1},
we have the specialization isomorphism
\begin{equation}
\label{special-1}
A_q(\n_+(w))_\iA \otimes \Cset \cong \Cset[N_+(w)]
\end{equation}
for the map $\iA \to \Cset$ given by $q \mt 1$.
By \cite[Proposition 9.7]{VY}, $I_w(u) \cap A_q(\n_+(w))_\iA$ is an $\iA$-form of $I_w(u)$. 
The definitions of $J_u^-$ and $I_w(u)$ in terms of Demazure modules imply 
that under the specialization isomorphism \eqref{special-1},
$I_w(u) \cap A_q(\n_+(w))_\iA$ is mapped to functions that 
vanish on the nonempty set 
\begin{equation}
\label{Richards}
N(w) \cap {\mathcal{B}}_- u {\mathcal{B}}^{\min}_+ w^{-1},
\end{equation}
which is isomorphic to the open Richardson variety in the flag scheme of $G^{\min}$ 
corresponding to the pair $u \leq w \in W$.
Let $\mu \in P_+$.  Analogously to the quantum situation, using special representatives 
of $w \in W$ in the normalizer of $H$ in $G^{\min}$,  
one defines the generalized minor ${\ol{\De}}_{u \mu, w \mu}$ which is a strongly regular function on $G^{\min}$. 
It is well known that under the specialization isomorphism \eqref{special-1}, 
the element $D_{u \mu, w \mu} \in E_w(u)$ corresponds to the restriction of ${\ol{\De}}_{u \mu, w \mu}$ to $N^{\min}_+$. This function is 
nowhere vanishing on the set \eqref{Richards}. Therefore, the 
specializations $I_w(u)$ and $E_w(u)$ are disjoint, so 
$I_w(u) \cap E_w(u) = \varnothing$.

For the next step, we denote for brevity
\[
c_\xi:= \psi_w( c_{\xi, v_\la} \De_{w \la, \la}^{-1}) \in A_q(\n_+(w)) 
\quad \mbox{for} \quad \xi \in D_\vp ( V(\la) ), \ \la \in P_+.
\]
For $\JJ \in \Spec A_q(\n_+(w))$ and $\la \in P_+$, denote
\[
C_{\JJ}(\la) = \{ \nu \in P \mid \exists \xi \in (V(\la)_\nu)^* \subset D_\vp (V(\la)) \; \; \mbox{such that} \; \; 
c_\xi \notin \JJ \}.
\]
Since $c_{\xi_{w \la}} = 1 \notin \JJ$, $w \la \in C_{\JJ}(\la)$. 
Denote by $M_{\JJ}(\la)$ the set of maximal elements of $C_{\JJ}(\la)$
with respect to the partial order $\nu \preceq \nu'$ if $\nu' - \nu \in Q_+$.

{\em{Step 3.}} {\em{For every $\JJ \in \Spec A_q(\n_+(w))$, there exists a unique $u \in W^{\leq w}$ such that $M_{\JJ}(\la) = \{ u \la \}$
for all $\la \in P_+$.}}

This step is similar to  \cite[Proposition 9.3.8]{Joseph-book}. 
Let  $\la \in P_+$ and $\nu \in M_{\JJ}(\la)$, so there exists $\xi \in (V(\la)^*)_\nu$ such that $c_\xi \notin \JJ$. The ${\mathcal{R}}$-matrix commutation relations in 
$A_q(\g)$ (see e.g. \cite[Proposition 9.1.5]{Joseph-book}) and the homomorphism from \thref{hom1} imply that
\[
c_\xi  x  \equiv q^{ - ( \nu+ w \la, \ga ) } 
x c_\xi \mod \JJ, \; \; 
\forall x \in A_q(\n_+(w))_\ga, \ \ga \in Q_+.
\] 
Take any other pair $\la' \in P_{++}$ and $\nu' \in M_{\JJ}(\la')$ going with $\xi' \in (V(\la')^*)_{\nu'}$ such that $c_{\xi'} \notin \JJ$. 
Applying the last relation twice gives 
\[
c_\xi c_{\xi'} \equiv 
q^{- ( \nu + w \la, \nu' - w \la')  
- (\nu- w \la, \nu' + w \la') } 
c_{\xi'}  c_\xi
\mod \JJ.
\]
Since $A_q(\n_+(w))/ \JJ$ is a prime ideal and the images of $c_\xi$, $c_{\xi'}$ are nonzero normal elements, they are regular.
Therefore the power of $q$ above must equal 0, and thus,
\begin{equation}
\label{lala'}
(\la, \la') - (\nu, \nu') = 0. 
\end{equation}
It follows from \cite[Lemma A.1.17]{Joseph-book} that $\nu = u_{\la} (\la)$ for some $u_\la \in W$; that is $M_{\JJ}(\la) = \{ u_\la \la \}$
(note that $u_\la$ is non-unique for $\la \in P_+ \backslash P_{++}$). 
It follows from the inclusion relations for Demazure modules \cite[Proposition 4.4.5]{Joseph-book} and the definition of $J_w^+$ 
that $u_\la \in W^{\leq w}$ for $\la \in P_{++}$. Applying one more time \eqref{lala'} gives that $u_\la = u_{\la'}$ for 
$\la, \la' \in P_{++}$ and that for $\la \in P_+$, $\la' \in P_{++}$, the element $u_\la$ can be chosen 
so that $u_\la = u_{\la'}$. 

{\em{Step 4. Completion of proof.}} By step 3, 
\begin{align}
&\Spec A_q(\n_+(w)) = \bigsqcup_{u \in W^{\leq w}} \Spec'_u A_q(\n_+(w)), \quad \mbox{where}
\label{de'}
\\
&\Spec'_u A_q(\n_+(w)) := \{ \JJ \in \Spec A_q(\n_+(w)) \mid M_{\JJ}(\la) = \{ u \la \}, \ 
\forall \la \in P_+\}.
\nn 
\end{align}
Steps 1, 2 and 3 and the fact that $\dim V(\la)_{w \la} =1$ imply the following: 

(*) For all $u \in W^{\leq w}$, we have $I_w(u) \in \Spec'_u A_q(\n_+(w))$, all ideals in $\Spec'_u A_q(\n_+(w))$ contain $I_w(u)$ and  
the stratum $\Spec'_u A_q(\n_+(w))$ contains no other $Q_+$-graded prime ideals. 

Therefore $\{ I_u(w) \mid u \in W^{\leq w} \}$ exhaust all $Q_+$-graded prime ideals of $A_q^+(\n_+(w))$. 
For $u_1 \leq u_2$ in $W^{\leq w}$, we have $I_{u_1}(w) \subseteq I_{u_2}(w)$ because 
$V^-_{u_1} (\la) \supseteq V^-_{u_2}(\la)$. Step 2 and the inclusion relations between Demazure modules \cite[Proposition 4.4.5]{Joseph-book}
imply that there are no other inclusions between these ideals. This proves part (a). 

All prime ideals of $A_q(\n_+(w))$ are completely prime by \cite[Theorem 2.3]{Go-Le-94}. 
It follows from (*) and the definition of $M_{\JJ}(\la)$ that the stratum
$\Spec_u A_q(\n_+(w))$, defined in part (b) of the theorem, coincides with $\Spec'_u A_q(\n_+(w))$ and equals 
\[
\{ \JJ \in \Spec A_q(n_+(w)) 
\mid \JJ \supseteq I_w(u), \ \JJ \cap E_w(u) = \varnothing \}. 
\]
The second statement in part (b) follows from \eqref{de'}, or equivalently, from \cite[\S II.2.1]{BG}. 

The properties (*) imply that the ring $(A_q(\n_+(w))/I_u(w))[E_u(w)^{-1}]$ is $\HH$-simple since the 
stratum $\Spec'_u A_q(\n_+(w))$ has a unique $Q_+$-graded ideal. This and step 2 prove part (c). Part (d) now
follows from \cite[Lemma II.3.7, Proposition II.3.8, Theorem II.6.4]{BG}. 
\end{proof}

\subsection{The homogeneous prime elements of $A_q(\n_+(w))$} 
Denote the support of $w$:
\[
\Sbb(w):= \{ i \in I \mid s_i \leq w \} = \{ i \in I \mid i= i_k \; \; \mbox{for some} \; \; k \in [1, N]\}
\]
where the second formula 
is in terms of a reduced expression \eqref{reduced}.
\bco{primeAn} The homogeneous prime elements of $A_q(\n_+(w))$ up to scalar multiples are 
\begin{equation}
\label{Uwprimes}
D_{\vpi_i, w \vpi_i} \quad \mbox{for} \quad i \in \Sbb(w).
\end{equation}
\eco

\begin{proof} \thref{spec-Uw}(i) implies that the height one $Q_+$-graded prime ideals of $A_q(\n_+(w))$ are $I_w(s_i)$ 
for $i \in \Sbb(w)$. Since $A_q(\n_+(w)) \cong U_q(\n_-(w))$ is a CGL extension (\leref{Uw}), it is an
$\HH$-UFD; thus, its height one $Q_+$-graded prime ideals are principal and their generators are 
precisely the homogeneous prime elements of $A_q(\n_+(w))$. Applying \thref{spec-Uw}(c) for $u =1$ 
and taking into account that $I_w(1) = 0$ gives that $I_w(s_i) \cap E_w(1) \neq \varnothing$ 
for $i \in \Sbb(w)$.
However, $E_w(1)$ consists of monomials in the elements \eqref{Uwprimes}. 
Hence each of the (completely prime) ideals $I_w(s_i)$, $i \in \Sbb(w)$ is generated by one of the elements in 
\eqref{Uwprimes}. The two sets have the same number of elements and 
$D_{\vpi_i, w \vpi_i} \in I_w(s_i)$. Hence, 
\[
I_w(s_i) = D_{\vpi_i, w \vpi_i} A_q(\n_+(w)), \quad \forall i \in \Sbb(w),
\]
and the set \eqref{Uwprimes} exhausts all homogeneous prime elements of $A_q(\n_+(w))$ up to scalar multiples.
\end{proof}

\sectionnew{Integral cluster structures on $A_q(\n_+(w))$}  \label{intAqn+w}
\label{qunip}
\subsection{Statements of main results}
\label{7.1}
Recall the notation \eqref{iAA}. 
Throughout the section, $\g$ denotes an arbitrary symmetrizable Kac--Moody algebra and $w$ a Weyl group element.
We fix a reduced expression \eqref{reduced}. 
Set
\[
U_q(\n_-(w))\spcheck_{\iAA} := U_q(\n_-(w))\spcheck_\iA \otimes_{\iA} \iAA, \quad
A_q(\n_+(w))_{\iAA} := A_q(\n_+(w))_{\iA} \otimes_{\iA} \iAA
\]
and extend $\iota$ to an algebra isomorphism
\begin{equation}
\label{iota1/2}
\iota : U_q(\n_-(w))\spcheck_{\iAA} \stackrel{\cong}{\lra} A_q(\n_+(w))_{\iAA}.
\end{equation}
For $k \in [1,N]$, denote
\begin{equation}
\label{x-gen}
x_k := q_{i_k}^{1/2} \iota(f_{\be_k}^*) = 
q_{i_k}^{1/2} (q_{i_k}^{-1} - q_{i_k}) \iota(f_{\be_k}) 
\in A_q(\n_+(w))_{\iAA},
\end{equation}
recall \eqref{dual-PBWgen}. For $j < k \in [1,N]$, set
\begin{equation}
\label{ajk}
a[j,k] := \| (w_{[j,k]} - 1) \vpi_{i_k} \|^2 /4 \in \Zset / 2.
\end{equation}

By applying $\iota$ to \eqref{UwCGL} and extending the scalars from $\Qset(q)$ to $\Qset(q^{1/2})$, we see that 
$A_q(\n_+(w)) \otimes_{\Qset(q)} \Qset(q^{1/2})$ is a symmetric CGL extension 
on the generators $\iota(f_{\be_1}^*), \ldots, \iota(f_{\be_N}^*)$. It follows from \leref{Uw}(a) that the scalars $\la_l, \la_k^*$ 
of the CGL extension are given by  
\begin{equation}
\label{lala}
\la_k = q^2_{i_k}, \quad \la_k^* = q_{i_k}^{-2}, \quad
\forall k \in [1,N].
\end{equation}
\leref{Uw}(b) implies that $A_q(\n_+(w))_{\iAA}$ with the generators $x_1, \ldots, x_N$ is an $\iAA$-form 
of the symmetric CGL extension $A_q(\n_+(w)) \otimes_{\Qset(q)} \Qset(q^{1/2})$. It follows from \eqref{LS} 
that the scalars 
\[
\nu_{kj} := q^{ ( \be_k, \be_j)/2}, \quad \forall 1 \leq j < k \leq N
\]
satisfy Condition (A) in \S \ref{normalizations}. 
\leref{Uw}(c) implies that the interval subalgebras of $A_q(\n_+(w))_{\iAA}$ are 
\begin{equation}
\label{interval-Anw}
\big( A_q(\n_+(w))_{\iAA} \big)_{[j,k]} = T_{w_{\leq j-1}^{-1}}^{-1}  \big( A_q(\n_+(w_{[j,k]}]))_{\iAA} \big), 
\quad \forall j \leq k \ \text{in} \ [1,N].
\end{equation}

Our first main theorem on quantum Schubert cells is:
\bth{main1} Let $\g$ be a symmetrizable Kac--Moody algebra and $w \in W$ with a reduced expression \eqref{reduced}. 
Consider the $\iAA$-form $A_q(\n_+(w))_{\iAA}$ of the symmetric CGL extension 
$A_q(\n_+(w)) \otimes_{\Qset(q)} \Qset(q^{1/2})$ with the generators $x_1, \ldots, x_N$
given by \eqref{x-gen}.
\begin{enumerate}
\item[(a)] The sequence of prime elements from Theorem {\rm\ref{t1}} of $A_q(\n_+(w)) \otimes_{\Qset(q)} \Qset(q^{1/2})$ 
with respect to the generators $x_1, \ldots, x_N$ is 
\[
y_k = q_{i_k}^{(O_-(k) + 1)/2} D_{\vpi_{i_k}, w_{\leq k} \vpi_{i_k}}, \quad k=1, \ldots, N. 
\]
The corresponding sequence of normalized prime elements is 
\[
\ol{y}_k = q^{a[1,k]} \, D_{ \vpi_{i_k}, w_{\leq k} \vpi_{i_k}}, \quad k=1, \ldots, N.
\]
Moreover, $y_1, \ldots, y_N, \ol{y}_1, \ldots, \ol{y}_N \in A_q(\n_+(w))_{\iAA}$.
\item[(b)] The $\eta$-function $\eta : [1,N] \to \Zset$ of $A_q(\n_+(w)) \otimes_{\Qset(q)} \Qset(q^{1/2})$ from \thref{1} is given by
\begin{equation}
\label{eta-Uw}
\eta(k) := i_k, \quad \forall k \in [1,N].
\end{equation}
\item[(c)] The normalized interval prime elements of $A_q(\n_+(w))_{\iAA}$ are
\[
\ol{y}_{[j,k]} = q^{a[j,k]}
\, D_{w_{\leq j -1} \vpi_{i_k}, w_{\leq k} \vpi_{i_k}} = \
q^{a[j,k]}
\,T_{w_{\leq k -1}} D_{\vpi_{i_k}, w_{[j,k]} \vpi_{i_k}}
\]
for all $j < k$ in $[1,N]$ such that $i_j = i_k$.
\end{enumerate}
\eth

In the rest of this section we will use the notation \eqref{pred.succ} for the predecessor and successor functions $p$, $s$ and 
the notation \eqref{O-+} for the functions $O_\pm$ associated to the $\eta$-function 
\eqref{eta-Uw}. Eq. \eqref{D-normal} and \thref{main1} imply that for $k >j$, 
\[
D_{\vpi_{i_k}, w_{\leq k} \vpi_{i_k}} D_{ \vpi_{i_j}, w_{\leq j} \vpi_{i_j}} = 
 q^{- ( (w_{\leq k} + 1)\vpi_{i_k}, (w_{\leq j} - 1) \vpi_{i_j} )}
 D_{ \vpi_{i_j}, w_{\leq j} \vpi_{i_j}} D_{\vpi_{i_k}, w_{\leq k} \vpi_{i_k}},
\]
and thus there is a unique toric frame 
$M^w \colon \Zset^N \to \Fract(A_q(\n_+(w))\spcheck_{\iAA})$ with cluster variables
\[
M^w(e_k) = q^{a[1,k]} \, D_{\vpi_{i_k}, w_{\leq k} \vpi_{i_k}}, 
\; \; \forall k \in [1,N]
\]
and matrix $\rb^w$ with
\begin{equation}
\label{rM-q}
(\rb^w)_{kj} := q^{- ( (w_{\leq k} + 1)\vpi_{i_k}, (w_{\leq j} - 1) \vpi_{i_j} ) /2},
\; \; \forall \, 1 \leq j < k \leq N. 
\end{equation}
We will use a quantum cluster algebra in which the exchangeable variables are
\begin{equation}
\label{exch-u}
\ex(w):=\{ k \in [1,N] \mid s(k) \neq \infty \}. 
\end{equation}
The number of elements of this set is $N- |\Sbb(w)|$. We will index the columns of the exchange matrices of this quantum cluster algebra 
(which have sizes $N \times (N - |\Sbb(w)|)$) by the elements of the set $\ex(w)$.

\bpr{compet-initial} The matrix $\wt{B}^w$ of size $N \times (N - |\Sbb(w)|)$ with entries
\[
(\wt{B}^w)_{jk} = 
\begin{cases}
1, &\mbox{if} \; \; j = p(k)
\\
-1, &\mbox{if} \; \; j = s(k)
\\
a_{i_j i_k}, & \mbox{if} \; \; 
j < k < s(j) < s(k)
\\ 
- a_{i_j i_k}, & \mbox{if} \; \; 
k < j < s(k) < s(j)
\\
0, & \mbox{otherwise}
\end{cases}
\]
is compatible with $\rb^w$, and more precisely, its columns $(\wt{B}^w)^k$, $k \in \ex(w)$ satisfy 
\begin{equation}
\label{comp}
\Om_{\rb^w}((\wt{B}^w)^k, e_l) = q_{i_k}^{-\delta_{kl}} = (\la^*_k)^{\delta_{kl}/2}\quad \mbox{and} \quad \sum_j (\wt{B}^w)_{jk}  (w_{\leq j}-1)  \vpi_{i_j} =0
\end{equation}
for all $k \in \ex(w)$, $l \in [1,N]$, recall \eqref{lala}.
\epr
The next theorem relates the integral quantum cluster algebra and upper quantum cluster algebra with initial seed $(M^w, \wt{B}^w)$
(both defined over $\iAA$) to the algebra $A_q(\n_+(w))_{\iAA}$.

\bth{main2} In the setting of Theorem {\rm\ref{tmain1}} the following hold:
\begin{enumerate}
\item[(a)] $A_q(\n_+(w))_{\iAA} = \AA(M^w, \wt{B}^w, \varnothing)_{\iAA} = \UU(M^w, \wt{B}^w, \varnothing)_{\iAA}$.
\item[(b)] For each $\sig \in \Xi_N \subset S_N$, the quantum cluster algebra $\AA(M^w, \wt{B}^w, \varnothing)_{\iAA}$ has a seed
with cluster variables
\[
M^w_\sig(e_l) = q^{a[j,k]}
\, D_{w_{\leq j-1} \vpi_{i_j}, w_{\leq k} \vpi_{i_j}} = 
q^{a[j,k]}
\,T_{w_{\leq j-1}} D_{\vpi_{i_j}, w_{[j,k]} \vpi_{i_j}}
\]
for $j := \min \sig([1,l])$ and $k := \max \{ m \in \sig([1,l]) \mid i_m = i_j \}$.
The initial seed $(M^w, \wt{B}^w)$ equals the seed corresponding to $\sig = \id_N \in \Xi_N$. 
\item[(c)] The seeds in {\rm(b)} are linked by sequences of one-step mutations of the following kind:

Suppose that $\sig, \sig' \in \Xi_N$ are such that $\sig' = ( \sig(k), \sig(k+1)) \circ \sig = \sig \circ (k, k+1)$
for some $k \in [1,N-1]$. If $\eta(\sig(k)) \neq \eta (\sig(k+1))$, then $M^w_{\sig'} = M^w_\sig \cdot (k,k+1)$
in terms of the action \eqref{r-act}.
If $\eta(\sig(k)) = \eta (\sig(k+1))$, then $M^w_{\sig'} = \mu_k (M^w_\sig)  \cdot (k,k+1)$.
\end{enumerate}
\eth
\red{We illustrate \thref{main2} and the constructions in \S\S \ref{5.4}, \ref{6.1} and \ref{7.1} 
with two examples of quantum unipotent cells in  non-symmetric Kac--Moody algebras $\g$:
a finite dimensional one and an affine one. 
}
\red{
\bex{B2} Let $\g$ be of type $B_2$ and $w$ be the longest Weyl group element $s_1 s_2 s_1 s_2$. 
The corresponding root sequence \eqref{roots} is
\[
\be_1 = \al_1, \; \; \be_2 = s_1 (\al_2) = \al_1 + \al_2, \; \; \be_3 = s_1 s_2 (\al_1) = \al_1 + 2 \al_2, \;  \; \mbox{and} \; \; \be_4 = s_1 s_2 s_1(\al_2) = \al_2.
\] 
The root vectors $f_{\be_k}$, $1 \leq k \leq 4$, satisfy
\begin{align*}
&f_{\be_2} f_{\be_1} = q^2  f_{\be_1} f_{\be_2}, &&f_{\be_3} f_{\be_1} = f_{\be_1} f_{\be_3} +\frac{1- q^{-2}}{q^{-1} + q} f_{\be_2}^2, 
&&&f_{\be_3} f_{\be_2} = q^2  f_{\be_2} f_{\be_3},  
\\
&f_{\be_4} f_{\be_1} = q^{-2} f_{\be_1} f_{\be_4} - q^{-2} f_{\be_2}, &&f_{\be_4} f_{\be_2} = f_{\be_2} f_{\be_4} - (q^{-1} + q) f_{\be_3},
&&&f_{\be_4} f_{\be_3} = q^2  f_{\be_3} f_{\be_4}.
\end{align*}
Note that the scalar in the right hand side of the second equation is not in $\iA$. 
The CGL extension $U_q(\n_-(w))=U_q(\n_-)$ is the $\Cset(q)$-algebra with these generators and relations. 
Its $\eta$-function from \thref{1} is given by $\eta(1)=\eta(3) =1$, $\eta(2)=\eta(4)=2$. 
The generators of the integral form $U_q(\n_-)\spcheck_{\iA}$ of the CGL extension  $U_q(\n_-)$
(cf. \eqref{dual-PBWgen} and \leref{Uw}(b)) are 
\[
f_{\be_k}^* = c_k f_{\be_k}, \quad \mbox{where} \; \; 
c_1 = c_3 = q^{-2} - q^2, \; \; 
c_2 = c_4 = q^{-1} - q.
\]
They satisfy 
\begin{align*}
&f_{\be_2}^* f_{\be_1}^* = q^2  f_{\be_1} f_{\be_2}, &&f_{\be_3}^* f_{\be_1}^* = f_{\be_1}^* f_{\be_3}^* - q^{-1} (q^{-2} - q^2) (f_{\be_2}^*)^2,
\\
&f_{\be_3}^* f_{\be_2}^* = q^2  f_{\be_2}^* f_{\be_3}^*,  
&&f_{\be_4}^* f_{\be_1}^* = q^{-2} f_{\be_1}^* f_{\be_4}^* - q^{-2} (q^{-2} - q^2)  f_{\be_2}^*, 
\\
&f_{\be_4}^* f_{\be_2}^* = f_{\be_2}^* f_{\be_4}^* - (q^{-1} - q) f_{\be_3}^*,
&&f_{\be_4}^* f_{\be_3}^* = q^2  f_{\be_3}^* f_{\be_4}^*.
\end{align*}
Recall the isomorphism \eqref{iota1/2}. The rescaled generators of $A_q(\n_+(w))_{\iAA}= A_q(\n_+)_{\iAA}$ are
\[
x_k = c'_k \iota(f_{\be_k}^*), \quad \mbox{where} \; \; 
c'_1 = c'_3 = q, \; \; 
c'_2 = c'_4 = q^{1/2}. 
\]
The algebra $A_q(\n_+)_{\iAA}$ is the $\iAA$-algebra with generators $x_1, \ldots, x_4$ and relations
\begin{align*}
&x_2 x_1 = q^2  x_1 x_2, &&x_3 x_1 = x_1 x_3 - (q^{-2}- q^2) x_2^2, 
&&&x_3 x_2 = q^2 x_2 x_3,  
\\
&x_4 x_1 = q^{-2} x_1 x_4 - q^{-1}(q^{-2} - q^2) x_2, &&x_4 x_2 = x_2 x_4 - (q^{-1} - q) x_3,
&&&x_4 x_3 = q^2  x_3 x_4.
\end{align*}
By  \thref{main2}, $A_q(\n_+)_{\iAA}$ has the structure of a quantum cluster algebra over $\iAA$ with initial cluster variables 
\[
\ol{y}_1 = x_1, \quad \ol{y}_2 =  x_2, \quad
\ol{y}_3 = x_1 x_3 - q^{-2} x_2^2, \quad \ol{y}_4 = x_2 x_4 - q^{-1} x_3
\]
(where the 3rd and 4th variables are frozen)
and mutation matrix 
\[
\wt{B} = 
\begin{bmatrix}
0 & -1 
\\
2 & 0
\\
-1 & 0
\\
0 & -1
\end{bmatrix}.
\]
Note that this is the quantum cluster algebra of type $B_2$ with principal coefficients. 
\qed
\eex
}
\red{
\bex{A22} Let $\g$ be the twisted affine Kac--Moody algebra of type $A_2^{(2)}$, whose Dynkin diagram is
\[
\begin{tikzpicture}[start chain]
\dnodeanj{}{0}
\dnodeanj{}{1}
\path (chain-1) -- node {\QLeftarrow} (chain-2);
\end{tikzpicture}
\]
Following the standard convention \cite[Ch. 6,8]{Kac}, we label its simple roots by $\{0,1\}$ instead of $\{1,2\}$.
We have $(\al_0, \al_0)=2$, $(\al_1, \al_1)=8$, $d_0=1$, $d_1=4$ and $q_0 = q$, $q_1 =q^4$. 
Consider the Weyl group element $w = s_0 s_1 s_0 s_1 s_0$. 
The corresponding root sequence \eqref{roots} is
\begin{align*}
\be_1 &= \al_0,  &\qquad \be_2 &= 4\al_0+\al_1,  &\qquad \be_3 &= 3\al_0+\al_1,  \\
\be_4 &= 8\al_0+3\al_1,  &\be_5 &= 5\al_0+2\al_1.
\end{align*}
The root vectors $f_{\be_k}$, $1 \leq k \leq 5$, given by \eqref{roots}, satisfy the relations 
\begin{equation}  \label{[15]rels}
\begin{gathered}
\begin{aligned}
z_2 z_1 &= q^4 z_1 z_2,  &\qquad z_3 z_1 &= q^2 z_1 z_3 + a z_2,  &\qquad z_3 z_2 &= q^4 z_2 z_3,  \\
z_4 z_1 &= q^4 z_1 z_4 + \frac{ab}{q^4-q^2} z_3^3 , &z_4 z_2 &= q^8 z_2 z_4 + b z_3^4  &z_4 z_3 &= q^4 z_3 z_4, 
\end{aligned}  \\
\begin{aligned}
&\!\!\!\!\!\!\!\!z_5 z_1 = q^2 z_1 z_5 + \frac{abc (q^6-1)} {q^2 (q^2-1)^2 (q^8-1)} z_3^2,  &\qquad z_5 z_2 &= q^4 z_2 z_5 + \frac{bc}{q^4-q^2} z_3^3,  \\
&\!\!\!\!\!\!\!\!z_5 z_3 = q^2 z_3 z_5 + c z_4,   &z_5 z_4 &= q^4 z_4 z_5
\end{aligned}
\end{gathered}
\end{equation}
with $a = c = -q^2 [4]_q$ and $b = -q^2 (q^{-1}-q)^3 / [4]_q$, where $[n]_q = (q^n - q^{-n})/(q- q^{-1})$. 
Note that $b \notin \iA$. 
The CGL extension $U_q(\n_-(w))$ is the $\Cset(q)$-algebra with these generators and relations. 
Its $\eta$-function from \thref{1} is given by $\eta(1)=\eta(3)= \eta(5) =0$, $\eta(2)=\eta(4)=1$. 
The generators of the integral form $U_q(\n_-(w))\spcheck_{\iA}$ of the CGL extension  $U_q(\n_-(w))$
(cf. \eqref{dual-PBWgen} and \leref{Uw}(b)) are 
\[
f_{\be_k}^* = c_k f_{\be_k}, \quad \mbox{where} \; \; 
c_1 = c_3 = c_5 = q^{-1} - q, \; \; 
c_2 = c_4 = q^{-4} - q^4.
\]
They satisfy the relations \eqref{[15]rels} for $a = c=q(q^2-1), \ b= q^{-2} (q^8-1) \in \iA$, and furthermore, $U_q(\n_-(w))\spcheck_{\iA}$ is 
the $\iA$-algebra with these generators and relations. Recall the isomorphism \eqref{iota1/2}. 
The rescaled generators of $A_q(\n_+(w))_{\iAA}$ are
\[
x_k = c'_k \iota(f_{\be_k}^*), \quad \mbox{where} \; \; 
c'_1 = c'_3 = c'_5 = q^{1/2}, \; \; 
c'_2 = c'_4 = q^2. 
\]
They satisfy the relations \eqref{[15]rels} for $a = c = q^2-1, \ b = q^8-1 \in \iA \subset \iAA$, and furthermore,
$A_q(\n_+)_{\iAA}$ is the $\iAA$-algebra  
with these generators and relations. 
By  \thref{main2}, $A_q(\n_+)_{\iAA}$ has the structure of a quantum cluster algebra over $\iAA$ with initial cluster variables
\begin{align*}
\ol{y}_1 &= x_1, \qquad \ol{y}_2 = x_2,  \qquad
 \ol{y}_3 = q x_1 x_3 -q^{-1} x_2,  \qquad
 \ol{y}_4 = q^4 x_2 x_4 -q^{-4} x_3^4, \\
\ol{y}_5 &= q^3 x_1 x_3 x_5 -q x_2 x_5 -q^{-1} [3]_q x_3^3 -q x_1 x_4. 
\end{align*}
(where the 4th and 5th variables are frozen)
and mutation matrix 
\[
\wt{B} = 
\begin{bmatrix}
0 & -4 & 1
\\
1 & 0 & -1
\\
-1 & 4 & 0
\\
0 & -1 & 1
\\
0 & 0 & -1
\end{bmatrix}.
\]
\qed
\eex
}
\subsection{Proof of \thref{main1}} If $u_1, u_2 \in W$ are such that $\ell(u_1 u_2) = \ell(u_1) + \ell(u_2)$, 
then we have the decomposition 
\[
A_q(\n_+(u_1 u_2))_\iA = A_q(\n_+(u_1))_{\iA} \, T_{u_1^{-1}}^{-1} ( A_q(\n_+(u_2)))_\iA.
\]
This follows by applying the isomorphism $\iota$ to the dual PBW basis \eqref{Int-PBW} 
of $U_q(\n_-(u_1 u_2))\spcheck_\iA$. The next lemma shows the equality of the unipotent quantum minors 
in \thref{main1}(c) and that they belong to the correct integral forms.  
\ble{D-int} If $u_1, u_2 \in W$ are such that $\ell(u_1 u_2) = \ell(u_1) + \ell(u_2)$, then
\begin{equation}
\label{equal-D}
D_{u_1 \mu, u_1 u_2 \mu} = T_{u_1^{-1}}^{-1} D_{\mu, u_2 \mu} \in T_{u_1^{-1}}^{-1} A_q(\n_+(u_2))_{\iA} \subset A_q(\n_+(u_1 u_2)), 
\quad \forall \mu \in P_+.
\end{equation}
\ele
\begin{proof} It was proved in \cite[Proposition 6.3]{GLS} that $\Psi^{-1} ( D_{\mu, u_2 \mu} ) \in t \cdot \vp^*({\BB}^{\up})$
in the notation of Remarks \ref{rGLSrem1} and \ref{rGLSrem2}. \thref{basisUw} and the commutative diagram in 
\reref{GLSrem2} imply that $D_{\mu, u_2 \mu} \in \iota ( {\BB}^{\up}) \subset  A_q(\n_+(u_1))_{\iA}$. 

The equality \eqref{equal-D} can be derived from \cite[Proposition 7.1]{GLS} and \reref{GLSrem2}, but it also has a direct proof as follows.
For all $y_k \in U_q(\n_+(u_k))$, $k =1,2$ and $h \in P\spcheck$, we have
\begin{align*}
&\lcor D_{u_1 \mu, u_1 u_2 \mu}, y_1 T_{u_1^{-1}}^{-1} (y_2) q^h \rcor = \lcor \xi_{u_1 \mu}, y_1 T_{u_1^{-1}}^{-1} (y_2) v_{u_1 u_2 \mu} \rcor 
= \lcor \xi_{\mu}, T_{u_1^{-1}}(y_1) y_2 v_{u_2 \mu} \rcor
\\
&= \lcor \xi_\mu, y_2 v_{u_2 \mu} \rcor \epsilon(y_1)
= \lcor D_{\mu, u_2 \mu}, y_2 \rcor \epsilon(y_1) = ( \iota^{-1} (D_{\mu, u_2 \mu}), y_2 )_{RT} \, \epsilon(y_1) 
\\
&= ( T_{u_1^{-1}}^{-1} \iota^{-1} (D_{\mu, u_2 \mu}), T_{u_1^{-1}}^{-1} y_2 )_{RT} \, \epsilon(y_1) = 
\lcor T_{u_1^{-1}}^{-1} D_{\mu, u_2 \mu}, y_1 T_{u_1^{-1}}^{-1} (y_2) q^h \rcor, 
\end{align*}
where the sixth equality uses \eqref{RT-mon}.
\end{proof}

\begin{proof}[Proof of Theorem {\rm\ref{tmain1}}] 
We have
\[
e_{\be_k} v_{w_{\leq k-1} \vpi_{i_k}} = T_{w_{\leq k-1}^{-1}}^{-1} (e_{i_k} T_{i_k}^{-1} v_{\vpi_{i_k}}) =  T_{w_{\leq k-1}^{-1}}^{-1} v_{\vpi_{i_k}} = 
v_{w_{\leq p(k)} \vpi_{i_k}}
\]
and $e_{\be_k}^m v_{w_{\leq k-1} \vpi_i} =0$ for $m>1$. Hence
\[
\lcor D_{\vpi_{i_k}, w_{\leq k} \vpi_{i_k}}, y_1 e_{\be_k}^m \rcor = \delta_{m1} \lcor D_{\vpi_{i_k}, w_{\leq p(k)} \vpi_{i_k}}, y_1 \rcor 
\]
for all $y_1 \in U_q(\n_+(w_{\leq k-1}))$, $m \in \Zset_{\geq 0}$. It follows from \eqref{RT-mon} and \eqref{interval-Anw} that in 
$A_q(\n_+(w_{\leq k}))_\iA = \big( A_q(\n_+(w))_{\iAA} \big)_{[1,k]}  \subset A_q(\n_+(w))_\iA$ we have 
\begin{equation}
\label{recur}
D_{\vpi_{i_k}, w_{\leq k} \vpi_{i_k}} \equiv D_{\vpi_{i_k}, w_{\leq p(k)} \vpi_{i_k}} \iota(f^*_{\be_k})  \mod A_q(\n_+(w_{\leq k-1}))_\iA. 
\end{equation}
Therefore, 
\[
q_{i_k}^{(O_-(k)+1)/2} D_{\vpi_{i_k}, w_{\leq k} \vpi_{i_k}} 
\equiv q_{i_k}^{(O_-(p(k) ) +1)/2}  D_{\vpi_{i_k}, w_{\leq p(k)} \vpi_{i_k}} x_k 
\mod A_q(\n_+(w_{\leq k-1}))_{\iAA}
\]
for all $k \in [1,N]$. Part (b) and the first statement in part (a) now follow from \coref{primeAn}.

We have $w_{\leq k} \vpi_{i_k} = w_{\leq k-1} ( \vpi_{i_k} - \al_{i_k}) = w_{\leq p(k)} \vpi_{i_k} - \be_k$. Iterating this gives
\begin{align*}
a[1,k] &= \| w_{\leq k} \vpi_{i_k} - \vpi_{i_k} \|^2 / 4 = 
\| \be_{p^{O_-(k)}\red{(k)}} + \cdots + \be_k \|^2 
\\
&= (O_-(k) + 1) \| \al_{i_k} \|^2/4 + \sum_{0 \leq l \leq m \leq O_-(k)} ( \be_{p^l(k)}, \be_{p^m(k)})/2.  
\end{align*}
Therefore,
\begin{align*}
\ol{y}_k &= \Big( \prod_{0 \leq l \leq m \leq O_-(k)} \nu_{p^l(k) p^m(k)}^{-1} \Big) y_k 
\\
&= \Big( \prod_{0 \leq l \leq m \leq O_-(k)} q^{( \be_{p^l(k)}, \be_{p^m(k)})/2} \Big) q^{(O_-(k) + 1) \| \al_{i_k} \|^2 /4} D_{\vpi_{i_k}, w_{\leq k}} 
= q^{a[1,k]} D_{\vpi_{i_k}, w_{\leq k}}
\end{align*} 
which proves the second statement in part (a).

It follows from \leref{D-int} that  $y_1, \ldots, y_N, \ol{y}_1, \ldots, \ol{y}_N \in A_q(\n_+(w))_{\iAA}$.
Part (c) follows from \eqref{interval-Anw} and part (b).
\end{proof}

\subsection{Proof of \thref{main2}}

\begin{proof}[Proof of Proposition {\rm\ref{pcompet-initial}}] Extend $\wt{B}^w$ to an $(N+r) \times (N - \Sbb(w))$ matrix whose rows are indexed 
by $[-r, -1] \sqcup [1,N]$ and columns by $\ex(w)$ by setting
\[
(\wt{B}^w)_{-i, k} := 
\begin{cases}
1, &\mbox{if} \; \; i_k =i \; \; \mbox{and} \; \; p(k) = - \infty
\\
0, & \mbox{otherwise}
\end{cases}
\]
for $i \in [1,r]$, $k \in \ex(w)$.

Denote for simplicity $b_{jk}:= (\wt{B})_{jk}$. 
We apply \cite[Theorem 8.3 and \S 10.1]{BZ} to the double word $1, \ldots, r, - i_1, \ldots, i_N$, which gives 
\begin{multline}
\sum_{j=1}^N b_{jk} \sign(j-l) \big( (w_{\leq j} \vpi_{i_j}, w_{\leq l} \vpi_{i_l} ) - ( \vpi_{i_j}, \vpi_{i_l} ) \big) 
\label{comp1}
\\
+ \sum_{i=1}^r b_{-i, k} ( (w_{\leq j} -1) \vpi_{i_j}, \vpi_i  ) = 2 \delta_{kl} d_k 
\end{multline}
for all $k \in \ex(w)$, $l \in [1,N]$.
The graded nature of the seed corresponding to the double word (cf. \cite[Definition 6.5]{BZ}) means that
\begin{align}
& \sum_{j=1}^N b_{jk} w_{\leq j} \vpi_{i_j} + \sum_{i=1}^r b_{-i, k} \vpi_i = 0,
\label{comp2}
\\
& \sum_{j=1}^N b_{jk} \vpi_{i_j} + \sum_{i=1}^r b_{-i, k} \vpi_i = 0
\label{comp3}
\end{align}
for all $k \in \ex(w)$. Subtracting \eqref{comp2} from \eqref{comp1} gives the second identity in \eqref{comp}. 
The linear combination $\mbox{\eqref{comp1}} + ( \mbox{\eqref{comp2}}, \vpi_{i_l} ) -  ( \mbox{\eqref{comp3}}, w_{\leq l} \vpi_{i_l} )$ 
yields the identity
\[
\sum_{j=1}^N b_{jk} \sign(j-l) \big( (w_{\leq j} +1) \vpi_{i_j}, (w_{\leq l} -1) \vpi_{i_l} )  = 2 \delta_{kl} d_k
\]
for all $k \in \ex(w)$, $l \in [1,N]$, which is precisely the first identity in \eqref{comp} in view of \eqref{rM-q}.
\end{proof} 

\bpr{verify-cond} In the setting of Theorem {\rm\ref{tmain1}}, the $\iAA$-form $A_q(\n_+(w))_{\iAA}$ of the symmetric CGL extension 
$A_q(\n_+(w)) \otimes_{\Qset(q)} \Qset(q^{1/2})$ with the generators $x_1, \ldots, x_N$ from \eqref{x-gen}
satisfies all conditions in Theorem {\rm\ref{tqcl.intform}}. 
\epr

\begin{proof} The scalars $\nu_{kl}$ are integral powers of $q^{1/2}$ and thus are units of $\iAA$. Obviously condition 
(A) is satisfied for the base field $\KK = \Qset(q^{1/2})$. Recall from \leref{Uw}(a) and \eqref{lala} that
$$
\la_k = q^2_{i_k} = q^{2 d_{i_k}} = q^{ \| \al_{i_k} \|^2 } \quad \mbox{for} \quad k \in [1,N],
$$
and from \thref{main1}(b) that $\eta(k) = i_k$ for $k \in [1,N]$. Therefore, Condition (B) is satisfied for the positive integers $\{ d_i \mid i \in I \}$ from \eqref{di}. 

The homogenous 
prime elements $y_1, \ldots, y_N$ belong to $A_q(\n_+(w))_{\iAA}$ by \thref{main1}(a). 

It remains to show that the condition \eqref{cond} holds. Because of \eqref{interval-Anw} and \leref{D-int} it is sufficient 
to consider the case when $i=1$ and $s(i) =N$. Since the $\eta$-function of the CGL extension $A_q(\n_+(w))$ is given by \thref{main1}(b), 
this means that $i_1 = i_N =i$ and $i_k \neq i$ for $k \in [2,N-1]$. It is well known that for $\g = \sl_2$ and 
$l \geq n \in \Zset_{> 0}$
\[
e_1^{(l)} \cdot T_1^{-1} v_{n \vpi_1} = \delta_{l n} v_{n \vpi_1}.
\]
For $k \in [2,N-1]$, $T_i^{-1} v_{\vpi_i}$ is a highest weight vector for the copy of $U_q(\sl_2)$ inside $U_q(\g)$ generated by 
$e_{i_k}, f_{i_k}, h_{i_k}$ of highest weight $\lcor s_i \vpi_i, h_{i_k} \rcor = - a_{i_k i} \vpi_{i_k}$. Hence, for $l \geq - a_{i_k i}$, 
\[
e_{\be_k}^{(l)} \cdot T_{w_{\leq k}^{-1}}^{-1} T_i^{-1} v_{\vpi_i}
= T_{w_{\leq k-1}^{-1}}^{-1} ( e_{i_k}^{(l)} \cdot T_{i_k}^{-1} T_i^{-1} v_{\vpi_i}) 
= 
\delta_{ l, -a_{i_k i}} T_i^{-1} v_{\vpi_i}. 
\]
Set $a := (- a_{i_1 i}, \ldots, -a_{i_{N-2} i}) \in \Zset_{\geq 0}^{N-2}$. 
Iterating this and using \eqref{recur} and the identity $T_i^{-2} v_{\vpi_i} = - q_i^{-1}  v_{\vpi_i}$ 
gives 
\[
\big{\lcor} D_{\vpi_i, w \vpi_i} - q_i^{-1} x_1 x_N, e_{\be_{2}}^{(l_{2})} \ldots e_{\be_{N-1}}^{(l_{N-1})} \big{\rcor} = 
\begin{cases}
- q_i^{-1}, &\mbox{if} \; \; a = (l_2, \ldots, l_{N-1})
\\
0, &\mbox{if} \; \; a \precneqq (l_2, \ldots, l_{N-1})
\end{cases}
\]
with respect to the the reverse lexicographic order \eqref{prec}.
It follows from \eqref{RT-mon} and \eqref{x-gen} that
\[
\lt (D_{\vpi_i, w \vpi_i} - q_i^{-1} x_1 x_N ) = - q_i^{-1} \Big( \prod_{k=2}^{N-1} q_{i_k}^{-(a_{i_k i}^2 + a_{i_k i} +1)/2} \Big) 
x_2^{-a_{i_2 i}} \ldots x_{N-1}^{-a_{i_{N-1} i}}.   
\]
By a straightforward calculation with powers of $q$, one obtains from this that the condition \eqref{cond} is satisfied.
\end{proof}
\thref{main2} follows by combining Theorems \ref{tqcl.intform} and \ref{tmain1} and Propositions \ref{pcompet-initial} and \ref{pverify-cond}.


\end{document}